\documentclass[11pt,a4paper,reqno]{amsart}  
\usepackage{fullpage}
\usepackage[utf8]{inputenc}
\usepackage[english]{babel}
\usepackage[T1]{fontenc}
\usepackage{microtype, fancyhdr, lmodern, xcolor}

\usepackage{amsfonts, amsmath,  amssymb, mathrsfs, amscd}
\usepackage{faktor}
\allowdisplaybreaks

\usepackage{url, enumitem}
\usepackage[normalem]{ulem} 
\usepackage[noadjust]{cite}
\usepackage{hyperref}

\usepackage{multirow,bigdelim, caption}

\usepackage[all]{xypic}
\xyoption{matrix,arrow}
\xyoption{graph}   

\usepackage{graphicx}
\graphicspath{{../}{../ppm/}}

\newcommand{\figdir}{xfig/}  

\newcommand{\module}[1]{}
\newcommand{\moduler}[1]{}

\newcommand{\ppm}[1]{{\textcolor{red}{#1}}}

\newcommand{\ppmx}[1]{}

\theoremstyle{plain}
\newtheorem{theorem}{Theorem}[section]
\newtheorem{proposition}[theorem]{Proposition}

\newtheorem{lemma}[theorem]{Lemma}

\theoremstyle{remark}
\newtheorem{remark}[theorem]{Remark}

\theoremstyle{definition}
\newtheorem{definition}[theorem]{Definition}

\newcommand{\beq}{\begin{equation}}
\newcommand{\eq}{\end{equation}} 
\newcommand{\mat}{\left( \begin{array}}
\newcommand{\tam}{\end{array} \right)} 

\newcommand{\ignore}[1]{}  

\newcommand\Nn{\mathbb{N}}
\newcommand\Zz{\mathbb{Z}}

\newcommand{\C}{\mathbb{C}}

\newcommand{\Sym}{\Sigma}        
\newcommand{\SymXZ}{\Zz_2\times\Sigma}        

\newcommand{\LBcat}{{\mathsf{ LB}}}  
\newcommand{\Bcat}{{\mathsf{ B}}}    
\newcommand{\Mat}{{\mathsf{ Mat}}}   
\newcommand{\Vect}{{\mathsf{ Vect}}} 

\newcommand{\FF}{{\mathsf F}} 
\newcommand{\FFF}{{\mathsf f}} 

\newcommand{\Scat}{{\mathsf S}} 

\newcommand{\ket}[1]{ | #1 \rangle } 
\newcommand{\Ccat}{{\mathbf C}}

\newcommand{\mate}{{\tiny \mat{cc} a&b\\ c&d\tam}}
\newcommand{\HH}{{\mathcal H}}
\newcommand{\integ}{{natural}}
\newcommand{\Or}[2]{{#1}_{|_{{#2}}}}  

\newcommand{\modelss}{\models^{\!\!\!\!\!\circ}} 
\newcommand{\gaml}{\gamma^l} 
\newcommand{\gamr}{\gamma^r}

\newcommand{\NL}{N}  
\newcommand{\braidrel}[2]{(#1\otimes#2)(#2\otimes#1)(#1\otimes#2)
  =(#2\otimes#1)(#1\otimes#2)(#2\otimes#1)}

\newcommand{\vs}{\varsigma}  
\newcommand{\vsigma}{\vs}  
\newcommand{\botimes}{\boxtimes}

\newcommand{\AAAA}[1]{\underline{A}(#1)}  
\newcommand{\off}{/}
\newcommand{\upp}{+}
\newcommand{\doo}{-}
\newcommand{\uppi}{+_i}

\newcommand{\uppii}{+_{ii}}

\newcommand{\II}{{\mathcal I}}

\newcommand{\vi}{v}
\newcommand{\vj}{w}
\newcommand{\offmat}{\mat{cc} 0&1\\1&0\tam}
\newcommand{\omat}{\mat{cc} 1&0\\0&1\tam}
\newcommand{\Rec}{{\mathsf R}}     
\newcommand{\Cpq}{\C^{(p,q)}}
\newcommand{\germ}{\SSS^\C}
\newcommand{\grm}{\TTT^\C}
\newcommand{\Cx}{\C^{\times}} 

\newcommand{\cc}{charge-conserving}
\newcommand{\N}{{\mathbb N}}
\newcommand{\ul}[1]{\underline{#1}}
\newcommand{\Match}{{\mathsf{Match}}}
\newcommand{\Braid}{{\mathsf B}}
\newcommand{\PP}{{\mathsf P}}  
\newcommand{\SSS}{{\mathfrak S}}  
\newcommand{\TTT}{{\mathfrak T}}  
\newcommand{\ppp}{{\mathsf p}_{\Lambda}}
\newcommand{\orr}{\triangleright}

\newcommand{\GG}{{\mathcal G}}
\newcommand{\ZZ}{{\mathbb Z}}

\newcommand{\stackx}{ \stackrel{X}{=} }

\newcommand{\chippy}{parted}
\newcommand{\choppy}{choppy}
\newcommand{\roo}{Rule-of-1}
\newcommand{\configuration}{$\II$-configuration}
\newcommand{\SET}{{\mathsf{ Set}}}  
\newcommand{\Hom}{Hom}
\newcommand{\Power}{{\mathcal P}}   
\newcommand{\pig}{\overline{\pi}}
\newcommand{\bb}{{\mathsf b}} 
\newcommand{\rrrr}{{\mathsf r}} 
\newcommand{\alphaa}{\pmb{\alpha}} 

\newcommand{\ova}{\prime\hspace{.051em}}  
\newcommand{\ovar}{\hspace{.05em}\prime\hspace{.05em}}  

\newcommand{\kay}[2]{\raisebox{.1421in}{\xymatrix@R=3pt{ #1 \ar@{-}[d] \\ #2}}}
\newcommand{\kfy}[2]{\raisebox{.1421in}{\xymatrix@R=3pt{ #1 \ar@{-}[d] \ar@/^1pc/ @{-} [d] \\ #2}}}
\newcommand{\kayy}[3]{\raisebox{.231in}{\xymatrix@R=3pt{
#1 \ar@{-}[d] \ar@/^1pc/ @{-} [d] \\ #2 \ar@{-}[d] \\ #3}}}
\newcommand{\kfyy}[3]{\raisebox{.231in}{\xymatrix@R=3pt{
#1 \ar@{-}[d] \ar@/^1pc/ @{-} [dd] \\ #2 \ar@{-}[d] \\ #3}}}
\newcommand{\kffy}[3]{\raisebox{.231in}{\xymatrix@R=3pt{
#1 \ar@{-}[d] \\ #2 \ar@{-}[d] \ar@/^1pc/ @{-} [d] \\ #3}}}
\newcommand{\kfffy}[3]{\raisebox{.231in}{\xymatrix@R=3pt{
#1 \ar@{-}[d] \ar@/^1pc/ @{-} [d] \\ #2 \ar@{-}[d] \ar@/^1pc/ @{-} [d] \\ #3}}}
\newcommand{\kayyy}[4]{\raisebox{.231in}{\xymatrix@R=3pt{
#1 \ar@{-}[d] \ar@/^1pc/ @{-} [d] \\ #2 \ar@{-}[d] \\ #3 \ar@{-}[d] \\ #4}}}
\newcommand{\kfyyy}[4]{\raisebox{.231in}{\xymatrix@R=3pt{
#1 \ar@{-}[d] \ar@/^1pc/ @{-} [dd] \\ #2 \ar@{-}[d] \\ #3 \ar@{-}[d] \\ #4}}}
\newcommand{\kffyyy}[4]{\raisebox{.231in}{\xymatrix@R=3pt{
#1 \ar@{-}[d] \ar@/^1pc/ @{-} [ddd] \\ #2 \ar@{-}[d] \\ #3 \ar@{-}[d] \\ #4}}}
\newcommand{\kfffyyy}[4]{\raisebox{.231in}{\xymatrix@R=3pt{
#1 \ar@{-}[d] \ar@/^1pc/ @{-} [d] \\
#2 \ar@{-}[d] \ar@/^1pc/ @{-} [d] \\ #3 \ar@{-}[d] \\ #4}}}
\newcommand{\kffffyyy}[4]{\raisebox{.231in}{\xymatrix@R=3pt{
#1 \ar@{-}[d] \ar@/^1pc/ @{-} [d] \\
#2 \ar@{-}[d] \ar@/^1pc/ @{-} [dd] \\ #3 \ar@{-}[d] \\ #4}}}
\newcommand{\kfffffyyy}[4]{\raisebox{.231in}{\xymatrix@R=3pt{
#1 \ar@{-}[d] \ar@/^1pc/ @{-} [dd] \\
#2 \ar@{-}[d] \\ #3 \ar@{-}[d] \ar@/^1pc/ @{-} [d] \\ #4}}}
\newcommand{\kffffffyyy}[4]{\raisebox{.231in}{\xymatrix@R=3pt{
#1 \ar@{-}[d] \\
#2 \ar@{-}[d] \ar@/^1pc/ @{-} [d] \\ #3 \ar@{-}[d] \ar@/^1pc/ @{-} [d] \\ #4}}}
\newcommand{\kfffffffyyy}[4]{\raisebox{.231in}{\xymatrix@R=3pt{
#1 \ar@{-}[d] \ar@/^1pc/ @{-} [d] \\
#2 \ar@{-}[d] \ar@/^1pc/ @{-} [d] \\ #3 \ar@{-}[d] \ar@/^1pc/ @{-} [d] \\ #4}}}

\newcommand{\whisper}[1]{{\textcolor{lightgray}{#1}}}
\newlength{\ccm}
\newcommand{\frametitle}[1]{}

\newcommand{\frx}[2]{}
\newcommand{\fr}[2]{   #2 }
\newcommand{\chapter}{\section}

\newcommand{\befff}[2]{\fr{}{#2}}
\newcommand{\beff}[2]{\fr{}{#2}}
\newcommand{\befs}[3]{\fr{}{#3}}
\newcommand{\prooff}{\proof}

\newcounter{minidef}[section]

\newcommand{\mdef}{\refstepcounter{theorem} 
\medskip \noindent ({\bf \thetheorem}) }
\newcounter{minicapt}

\newcommand{\fii}{\mathsf{f}} 
\newcommand{\ai}{\mathsf{a}}  
\newcommand{\pfii}{{\fii}} 
\newcommand{\uppfii}{\pfii}
\newcommand{\pai}{{\ai}}
\newcommand{\mfii}{\underline{\fii}} 
\newcommand{\mai}{\underline{\ai}}
\newcommand{\uppai}{\pai}
\newcommand{\ooo}[1]{\textcolor{orange}{\textcircled{\tiny #1}}}

\newcommand{\ber}[5]{ 
\xymatrix{
       &  {a_{#1}}_{\ooo{2}} \ar[dr]|{#4} \\
{{a_1}}^{\ooo{1}} \ar[ur]|{#3} \ar[rr]|{#5} &&  {}_{\ooo{3}}{a_{#2}}
}}

\newcommand{\bermss}[6]{ 
\xymatrix{
       &  {a_{#2}}_{\ooo{2}} \ar[dr]|{#5}="b" \\
{{a_{#1}}}^{\ooo{1}} \ar[ur]|{#4}="a" \ar[rr]|{#6} &&  {}_{\ooo{3}}{a_{#3}}
\ar@{=}_{} "a";"b"
}}

\newcommand{\bermpp}[6]{ 
\xymatrix{
       &  {a_{#2}}_{\ooo{2}} \ar[dr]|{#5}="b" \\
{{a_{#1}}}^{\ooo{1}} \ar[ur]|{#4}="a" \ar[rr]|{#6}="c" &&  {}_{\ooo{3}}{a_{#3}}
\ar@{=}_{} "a";"b"
\ar@{=}_{} "a";"c"
\ar@{=}_{} "b";"c"
}}

\newcommand{\mbermpp}[9]{ 
       &  {a_{#2}}_{\ooo{2}} \ar[dr]|{#5}="b" \\
{{a_{#1}}}^{\ooo{1}} \ar[ur]|{#4}="a" \ar[rr]|{#6}="c" &&  {}_{\ooo{3}}{a_{#3}}
\ar@{=}|{#7} "a";"b"
\ar@{=}|{#8} "a";"c"
\ar@{=}|{#9} "b";"c"
}

\newcommand{\mbermupp}[9]{ 
       &  {{#2}}_{\ooo{2}} \ar[dr]|{#5}="b" \\
{{{#1}}}^{\ooo{1}} \ar[ur]|{#4}="a" \ar[rr]|{#6}="c" &&  {}_{\ooo{3}}{{#3}}
\ar@{=}|{#7} "a";"b"
\ar@{=}|{#8} "a";"c"
\ar@{=}|{#9} "b";"c"
}

\newcommand{\abcab}{\ar@{=}_{} "a";"b"}
\newcommand{\abcac}{\ar@{=}_{} "a";"c"}

\newcommand{\bermxx}[7]{ 
\xymatrix{
       &  {a_{#2}}_{\ooo{2}} \ar[dr]|{#5}="b" \\
{{a_{#1}}}^{\ooo{1}} \ar[ur]|{#4}="a" \ar[rr]|{#6}="c" &&  {}_{\ooo{3}}{a_{#3}}
#7
}}

\newcommand{\bertx}[9]{  
\xymatrix@R=37pt{
& {}^{a_{#1}} \ooo{2}_{} \ar[rr]^{#5} \ar[drrr]^{#8}  && \ooo{3}^{a_{#2}} \ar[dr]^{#7} \\
{}_{a_1}\ooo{1} \ar[ur]^{#4} \ar[urrr]^{#6} \ar[rrrr]^{#9} &&&& \ooo{4}_{a_{#3}}
}}

\newcommand{\bert}[9]{ 
\xymatrix@R=37pt{
& {a_{#1}}_{\ooo{2}} \ar[rr]^{#5} \ar[drrr]^{#8}  && {}_{\ooo{3}}{a_{#2}} \ar[dr]^{#7} \\
{{a_1}}^{\ooo{1}} \ar[ur]^{#4} \ar[urrr]^{#6} \ar[rrrr]^{#9} &&&& {}^{\ooo{4}}{a_{#3}}
}}

\newcommand{\berto}[9]{ 
\xymatrix@C=10pt{
1 \ar[drr]^{#4}="a"
  \ar@/^1pc/@{=}[rrrrrr]
           &&                  &          &
                                           && 1
                                           \ar[dll]^{#5}_{\;\;\;}="b"
                                           \ar@/^2pc/@{=}[dddlll] \\
           && #1 \ar[rr]^{#6}="c" \ar[dr]^{#7} &          & #3   \\
           &&                  & #2 \ar[ur]^{#8} &   \\
           &&                  & 1 \ar[u]^{#9} \ar@/^2pc/@{=}[uuulll]
}}

\newcommand{\bertovold}[9]{ 
\xymatrix@C=10pt{
1 \ar[drr]^{#4}="a"
  \ar@/^1pc/@{=}[rrrrrr]
           &&                  &          &
                                           && 1
                                           \ar[dll]^{#5}_{\;\;\;}="b"
                                           \ar@/^2pc/@{=}[dddlll] \\
           && #1 \ar[rr]^{#6}="c" \ar[dr]_{#7} &          & #3   \\
           &&                  & #2 \ar[ur]^{#8} &   \\
           &&                  & 1 \ar[u]^{#9} \ar@/^2pc/@{=}[uuulll]
           \ar@{=} "a";"b"
            \ar@{=} "b";"c"
            \ar@{=} "a";"c"
}}

\newcommand{\bertov}[9]{ 
\xymatrix@C=10pt{
1 \ar[drr]|{#4}="a"
  \ar@/^1pc/@{=}[rrrrrr]
           &&                  &          &
                                           && 1
                                           \ar[dll]|{#5}="b"
                                           \ar@/^2pc/@{=}[dddlll] \\
           && #1 \ar[rr]|{#6}="c" \ar[dr]|{#7} &          & #3   \\
           &&                  & #2 \ar[ur]|{#8} &   \\
           &&                  & 1 \ar[u]|{#9} \ar@/^2pc/@{=}[uuulll]
           \ar@{=} "a";"b"
            \ar@{=} "b";"c"
            \ar@{=} "a";"c"
}}

\newcommand{\abcccc}{ 
        \ar@{=} "a";"b"
        \ar@{=} "b";"c"
        \ar@{=} "a";"c"    }

\newcommand{\cdeeee}{ 
        \ar@{=} "c";"d"
        \ar@{=} "d";"e"
        \ar@{=} "c";"e"    }

\newcommand{\adffff}{ 
        \ar@/_/@{=} "a";"d"
        \ar@/_/@{=} "d";"f"
        \ar@/_/@{=} "a";"f"    }

\newcommand{\beffff}{ 
        \ar@/^/@{=} "b";"e"
        \ar@/^/@{=} "e";"f"
        \ar@/^/@{=} "b";"f"    }

\newcommand{\bertoxx}[9]{ 
\xymatrix@C=10pt{
1 \ar[drr]|{#4}="a"
  \ar@/^1pc/@{=}[rrrrrr]
           &&                  &          &
                                           && 1
                                           \ar[dll]^{#5}_{\;\;\;}="b"
                                           \ar@/^2pc/@{=}[dddlll] \\
           && #1 \ar[rr]|{#6}="c" \ar[dr]_{#7} &          & #3   \\
           &&                  & #2 \ar[ur]_{#8} &   \\
           &&                  & 1 \ar[u]^{9} \ar@/^2pc/@{=}[uuulll]
#9
}}

\newcommand{\bertox}[9]{ 
\xymatrix@C=10pt{
1 \ar[drr]|{#4}="a"
  \ar@/^1pc/@{=}[rrrrrr]
           &&                  &          &
                                           && 1
                                           \ar[dll]|{#5}="b"
                                           \ar@/^2pc/@{=}[dddlll] \\
           && #1 \ar[rr]|{#6}="c" \ar[dr]|{#7}="d" &          & #3   \\
           &&                  & #2 \ar[ur]|{#8}="e" &   \\
           &&                  & 1 \ar[u]|{\x}="f" \ar@/^2pc/@{=}[uuulll]
#9
}}

\newcommand{\pyr}[9]{
\xymatrix@C=30pt@R=13pt{
1 \ar[rr]|{#1} \ar[dr]|{#2} \ar[dddr]|{#3}
                         &                  & #9 \\
                         & #7 \ar[ur]|{#4} \ar[dd]|{#5}    \\ \\
                         & #8 \ar[uuur]|{#6}
}}

\newcommand{\qquadd}{\hspace{.21cm} & \hspace{.21in}}

\newcommand{\alp}{\alpha}
\newcommand{\be}{\beta}
\newcommand{\lemx}{the X~Lemma}
\newcommand{\loone}{\mat{cc} \!\alp\!+\!\be & -\alp\be\! \\ 1&0\tam  }
\newcommand{\ccc}{\bullet}
\newcommand{\trin}[6]{ \xymatrix{& #3\; \ccc\;\; \ar@{-}[dr]^{#4}
\\
#1 \ccc \ar@{-}[rr]_{#6} \ar@{-}[ur]^{#2} && \ccc #5  }}

\newcommand{\functor}{\mathfrak{Func}}  
\newcommand{\ninerule}{9-rule}          
\newcommand{\Perm}{Perm}                
\newcommand{\Gammall}{\Gamma}     
\newcommand{\UUU}{{\mathfrak U}}  
\newcommand{\soutx}[1]{\ignore{#1}} 
\title[On classifying braid representations]{
  Classification of \\ spin-chain braid representations}
\author{Paul Martin}
\address{School of Mathematics\\ University of Leeds \\ Leeds LS2 9JT \\ U.K.}
\email{p.p.martin@leeds.ac.uk}

\author{Eric C. Rowell}
\address{Department of Mathematics\\
    Texas A\&M University \\
    College Station, TX 77843\\
    U.S.A.}
\email{rowell@math.tamu.edu}

\date{\today}

\begin{document}
\maketitle

\module{abstractV1}

\begin{abstract}

A {\em braid representation} is a monoidal functor from the braid category $\Braid$,
for example given by a solution to the constant Yang-Baxter equation.
Given a monoidal category $\mathsf{C}$
{with $ob(\mathsf{C})=\N$},
a rank-$N$ \cc\ representation,
or {\em spin-chain} representation,
is a  strict  monoidal functor
$F$
from
$\mathsf{C}$  
to the category $\Match^N$ of rank-$N$ \cc\ matrices
(see 
below for definition),
{that is {\em natural} in the sense that  
$F(1)=1$}.
In this work we construct all
  spin-chain braid representations,
and classify up to suitable notions of isomorphism.

\end{abstract}


\tableofcontents

\module{tex-eric/intro-eric}

\section{Introduction }
\newcommand{\Aut}{\mathrm{Aut}}
\newcommand{\uN}{\underline{N}}

Classification of braid representations \cite{Birman} is a major open problem ---
except that it is impossibly wild.
For applications,
{as we will discuss below,}
conceding this impossibility is not an acceptable outcome.
So we  
seek a framework for a paradigm change.
{Clues for this can be found, for example, in}
higher representation theory
(works such as
\cite{MazorchukMiemietz14,RumyninWendland,KapVoev,DamianiMartinRowell}
and many others);
and in higher lattice gauge theory (works such as
\cite{MartinsPicken11,Pfeiffer,baez_schreiber}
and 
many others).
Armed with such ideas, 
we are indeed able to
make a paradigm change that yields 
a solution for monoidal functors to a suitable target category.

\newcommand{\nNn}{N}  

One 
formulation
\ppmx{natural formulation of a suitable target category
[I think I am not understanding this wording yet. The YBE does not itself lead to Match
does it? I think it says at most that things are monoidal, but we have already said that.
So I might lean towards reverting this wording to how it was?]
comes from}%
{exemplifying  
such braid representations is given by}
invertible solutions to the
constant ({also known as} quantum)
Yang-Baxter equation (YBE), i.e., $R\in \mathrm{Aut}(V^{\otimes 2})$
satisfying
\[ R_1R_2 R_1=R_2 R_1 R_2
\]
where $R_1=R\otimes I_V$ and $R_2=I_V\otimes R$
and $V=\C^\nNn$
{for some $N$}
\cite{Baxter,Jimbo86}.
\ppmx{[-by the way, if we refer to categories in this para then it might be best to
say what $\otimes$ means? - just in the sense that using the word raises the level
of formality a notch.]}
{The} 
YBE has been a well-spring of deep and beautiful mathematics from its origins in Baxter's approach to finding exactly solvable models in statistical mechanics and Yang's construction of 2-dimensional quantum field theories \cite{Yang67}, to Drinfeld and Jimbo's \cite{drinfeld,Jimbo86} discovery
of quantum groups and the field of quantum topology generalising Jones' celebrated polynomial link invariant \cite{Jones85} in the 1980s.
 Formally the $\nNn$-dimensional YBE is a  system of $\nNn^6$
 cubic polynomial equations in $\nNn^4$ variables--a
formidable problem in symbolic computation for even modest values of
 $\nNn$.
 Indeed, while a classification of the solutions $(R,\mathbb{C}^2)$ for $\nNn=2$ have been known for almost 30 years \cite{Hietarinta:1992} the cases
 $\nNn\geq 3$
 remain wide open. Finding and classifying Yang-Baxter operators (YBOs) 
  remain important classical problems with myriad applications.

%
An alternative approach to attempting  classification in ever-increasing dimension,
is to seek to generalize the forms of the solutions for $\nNn=2$.
%
For example we can point to the isolated 
solution
$\begin{pmatrix} 1 & 0&0&1\\ 0& 1 & 1& 0\\
0& -1 & 1& 0\\ -1 & 0&0&1\end{pmatrix}$
which can be naturally generalised to every dimension as the so-called Gaussian solutions \cite{GoldschmidtJones:1989,RowellWang:2010}.
Two-eigenvalue solutions are relatively well understood (see \cite{Martin92}, \cite{MartinWoodcock98b} and references therein).
And
generalising the
{basic $4\times 4$ flip matrix solution are the permutation-matrix 
solutions
(linearisations of set-theoretical solutions, see e.g. \cite{Etingof99}).}
%
However each of these classes of solutions is
manifestly measure-zero in the space of all solutions,
the second with only two eigenvalues,  and
the third with all eigenvalues of same magnitude. 
As we will see, our construction includes solutions with unboundedly many eigenvalues in general position. 
{
Of course there is  ambiguity in the   
instruction to 
`generalize', and  
a robust method
should employ available symmetries such as local basis changes to reduce the complexity of the problem and the solution description.
}

\newcommand{\mathbeef}{\mathfrak{e}}  

Recently with Damiani \cite{DamianiMartinRowell} we were led to consider
{$N=2$, i.e. }
$4\times 4$ spin-chain representations of the braid group,
{in particular, with only limited basis change symmetry},
in our study of certain finite dimensional quotients of
the loop braid group {algebra}.
The matrices we encountered have the same form as the 
universal $R$-matrix for the quantum supergroup $U_q\mathfrak{sl}_{1|1}$.
A salient feature these YBOs share with the $U_q\mathfrak{sl}_2$ $R$-matrix
{(which does not lift to loop-braid in general)}
is that they are \emph{charge conserving} in the sense that the subspaces spanned by tensor products of standard basis vectors
$\mathbeef_i\otimes \mathbeef_j$ and $\mathbeef_j\otimes \mathbeef_i$ are invariant under $R$ (here for $i,j\in\{1,2\}$),
i.e. of the general form $\begin{pmatrix} a_1& 0&0 & 0\\
0 & a_{12} & b_{12} &0\\
0 & c_{12} & d_{12} &0\\
0&0&0& a_2
\end{pmatrix}$
{as in an XXZ spin chain or 6-vertex model \cite{Baxter}}.
This suggests that an appropriate higher dimensional generalization of these
$4\times 4$ solutions should be
$\nNn^2\times \nNn^2$ charge conserving solutions to the YBE.  
{In \cite{DamianiMartinRowell} we observed that the universal $R$-matrix coming from  $U_q\mathfrak{sl}_{1|1}$ can be supplemented with another \emph{involutive} YBO to lift the braid group representations to \emph{loop} braid group representations, while the $R$-matrix from $U_q\mathfrak{sl}_2$ cannot.  Shedding light on this phenomena was our original motivation for studying charge conserving YBOs.
{Indeed our solution here has an immediate application: as the key step in obtaining a corresponding loop-braid classification \cite{MRT1}}.}

In this article we classify
charge conserving YBOs.
The problem
is first framed in terms of monoidal functors $F$ from the braid category $\Bcat$ to
a category $\Match^N$   of  
charge conserving matrices.  The automorphisms in this target category confer symmetries upon the space of solutions, represented as images of the monoidal generator $\sigma$ of the category $\Bcat$ under $F$.  We then characterize all such functors in terms of polynomial constraints.  Taken in tandem, these symmetries and constraints yield a blueprint for a complete classification.  The description of all solutions is facilitated by a remarkably simple combinatorial encoding of a transversal of the orbits of all solutions under the action of the group of symmetries.  

The tractability of the problem is a consequence of special
{`miracles'}
that do not apply universally to YBOs, principally:
\begin{enumerate}
\item[(a)] Lemma \ref{lem:3 suffice}: An operator $R\in \Aut(V^{\otimes 2})$ is a charge conserving YBO if, and only if, the restriction $R^{ijk}$ of $R$ to the subspace of $V^{\otimes 2}$ spanned by tensor products of $3$ basis vectors $\mathbeef_i,\mathbeef_j,\mathbeef_k$ is also a charge conserving YBO, for all triples $\{i,j,k\}$.
\item[(b)] $X$-symmetry, Lemma
\ref{lem:Xeq}:
Conjugation of a charge conserving YBO $R$ by \emph{any} diagonal matrix $X\in\Aut(V^{\otimes 2})$ again yields a  charge conserving YBO.
\item[(c)] Restricted Symmetry, Lemma 
\ref{lem:chacon1}:
Note that
general local basis changes of the form
$R\mapsto (Q\otimes Q)R(Q\otimes Q)^{-1}$
\emph{do not} preserve charge
conservation. 
Together with $X$-symmetry, the simultaneous basis permutation action of $\Sigma_N$ yields all local basis changes that preserve charge conservation.
\end{enumerate}
The upshot of (a) is that it suffices to first find all solutions with $N=3$ and then
figure out how to
glue them together for $N>3$.  Statements (b) and (c) provide enough
symmetries to organise the cases and ultimately gain purchase on the $N=3$ case.
\ppmx{[-the two parts of this sentence seem to have a very different flavour?
Not sure I understand it yet.]}

Our combinatorial description of solutions begins with an encoding of general charge conserving operator $T\in\Aut(V^{\otimes 2})$ as $N$ scalars $a_i$ associated with the $1$-dimensional spaces $\C\mathbeef_i\otimes\mathbeef_i$ and ${N\choose 2}$ $2\times 2$ matrices $A(i,j)$ associated with the $2$ dimensional spaces spanned by $\mathbeef_i\otimes \mathbeef_j$ and $\mathbeef_j\otimes \mathbeef_i$ for $i<j$.  This can be conveniently visualised as a directed complete graph $K_N$ on $N$ vertices, with the vertices decorated with the $a_i$ and the edges by $A(i,j)$.

Next, we {briefly} describe versions of the combinatorial objects to which we associate solutions.
(In order to
{work with and}
 count the solutions we will use a more refined construction in the main text.)
Let $\underline{N}:=\{1,\ldots, N\}$ be a set of $N$ individuals.
We construct a set
of bi-coloured
{\em composition tableaux}  
from $\underline{N}$ as follows:
\begin{enumerate}
\item First partition $\uN$ into $m$ `nations' $n_1,\ldots,n_m$ respecting the natural order on $\uN$, so that $1,\ldots,|n_1|$ lie in $n_1$;
{$|n_1|+1,\ldots,|n_1|+|n_2|$ lie in $n_2$,} 
etc.
\item Next partition each nation $n_s$ into $\ell_s$ `counties' $c_{s,1},\ldots, c_{s,\ell_s}$, again respecting the natural ordering on individuals. 
\item Next,
\ppmx{[-if the target is tableau then this step is not final?]}%
colour  
a subset of counties blue,
such that the first county in each nation is always
uncoloured.  
\item Finally observe that each nation can be visualised as a
bi-coloured
{composition tableau}
with the counties stacked on top of each other,
and the collection of nations is our set of bicoloured
{composition tableaux}.
\end{enumerate}
  Here is an example of a full configuration for $N=11$:

\begin{centering}

\includegraphics[width=3.7cm]{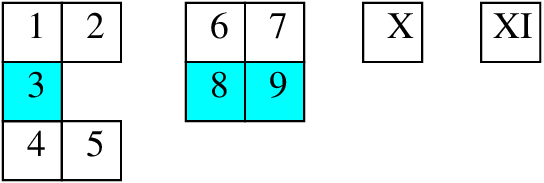}

\end{centering}

For each nation $n_s$ fix two non-zero parameters $\alpha_s,\beta_s$ such that $\alpha_s+\beta_s\neq 0$, and for each pair of nations $n_s,n_t$ with $s<t$ fix a non-zero parameter $\mu_{s,t}$.  To describe a solution it is enough to specify the scalars $a_1,\ldots,a_N$ and matrices $A(i,j)$, which depend on the relationship between the counties/nations that individuals $i<j$ reside in, as well as the colour of their respective counties.
Specifically:

\begin{enumerate}
\item If $i$
\ppmx{[-maybe, if possible, separate out a bit the variable names used for nations and individuals?]}
resides in county $c_{s,z}$ (which is in nation $n_s$) then if the county is 
uncoloured 
$a_i=\alpha_s$ otherwise $a_i=\beta_s$.
\item If $i$ is in nation $n_s$ and $j$ is in nation $n_t$ with $s\neq t$ then $A(i,j)=\begin{pmatrix} 0 & \mu_{s,t}\\ \mu_{s,t} &0\end{pmatrix}$.
\item If $i$ and $j$ are in the same nation $n_s$ but different counties $c_{s,x}$ and $c_{s,y}$ (here we have $x<y$ by construction),
then $A(i,j)=\begin{pmatrix} \alpha_s+\beta_s & \beta_s\\ -\alpha_s & 0\end{pmatrix}$.
\item If $i$ and $j$ are in the same county $c_{s,z}$ then
$A(i,j)=\begin{pmatrix} x_s &0\\0&x_s\end{pmatrix} $
\ppmx{[something not quite right?]}
where $x_s=\alpha_s$ if $c_{s,z}$ is
uncoloured  
and $x_s=\beta_s$
if blue.  
\end{enumerate}
Our main results imply that matrices so constructed are indeed charge conserving YBOs and any charge conserving YBO may be transformed to a solution of this form by means of the symmetries (b) and (c) described above.  In particular we obtain a count of orbits of
varieties of 
solutions by rank as $1, 4, 13, 46, 154,\ldots$ which,
pleasingly,
is  precisely
sequence A104460 in \cite{sloane1}.

\medskip

{A key application of our classification is to answer the question of when
charge-conserving YBOs can be lifted to yield representations of
the loop braid group.
{Indeed using the results in this paper}
we have
{now solved this problem,}
jointly with Fiona Torzewska
{\cite{MRT1}}.  
\\ 
There are also some natural generalisations of charge
conservation  
such as \emph{additive charge conservation} in which the subspace  $V^{i+j}\subset V^{\otimes 2}$ spanned by all the $\mathbeef_i\otimes\mathbeef_j$ with $i+j$ a constant is invariant under $R$.
{So far this problem is solved in smallest non-trivial rank, in work with Hietarinta and Nijhoff \cite{HMNR1}.}}
\ignore{
It is also interesting to ask which general solutions to the YBE are similar to a charge-conserving solution.  Although charge-conserving matrices are sparse, it is conceivable that a significant fraction of general solutions are equivalent to charge conserving solutions via a (potentially non-local) basis change.  
}

Here is a summary of the contents of this paper.
In section \ref{ss:braid} we describe the braid category.
In
section \ref{ss:target} we focus on functors to the target category $\Match^N$
{considered representation theoretically,}
\ppmx{[to contradistinguish with the generalist perspective above]}
and set up our approach to the main theorem, including the $N=2$ case.
The main Theorem \ref{th:mainx} is found in section  \ref{ss:main} which includes a detailed description of the classification in terms of our combinatorial encoding.
Section \ref{section: the proof} contains a sequence of lemmas covering the crucial $N=3$ case, which are then applied in section \ref{ss:proofend} to finish the proof of the main theorem.
We close with a brief discussion of some future directions in the last section.
(Appendices \ref{ss:repthy},\ref{ss:Xeqa}
briefly discuss various intriguing points raised by,
but not internal to, this work.
Appendices \ref{ss:deets},\ref{ss:figs} contain further details on
various aspects treated with brevity in the main text.)


\medskip 

\noindent 
{\bf{Acknowledgements}.}
We thank Lizzy Rowell for the braid laboratory.
We thank Celeste Damiani
for useful conversations.
PM thanks Paula Martin and Nicola Gambino 
for useful conversations. The research of ECR was partially supported by the US NSF grant DMS-1664359.

\medskip

\module{braid-pics}

\section{Braid categories} \label{ss:braid}
\fr{Braids}{
The monoidal category of braids is defined
for example in Mac Lane \cite{MacLane}.
Its  
calculus
is conveniently recalled by the following pictures
of mathematicians.
\begin{figure}[h]

\[
\includegraphics[width=4.85\ccm, trim=25 0 0 35, clip]{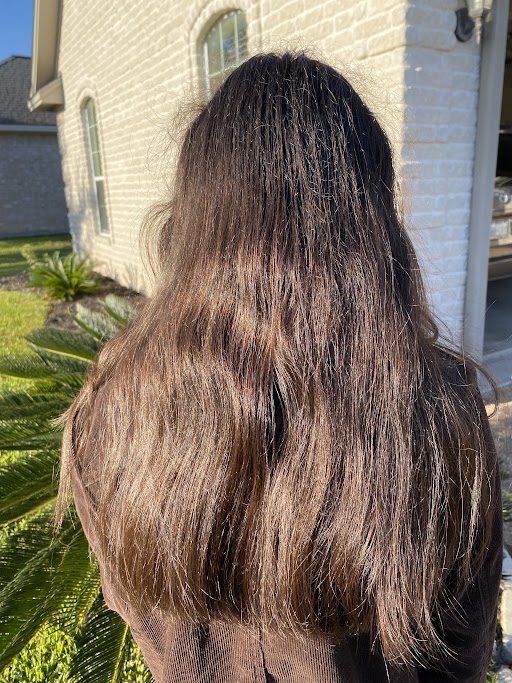}
\hspace{.21in}
\includegraphics[width=4.5\ccm, trim=25 0 25 0, clip]{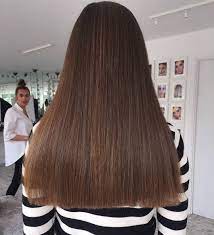}
\]
\caption{Two representatives in the Identity Morphism \label{fig:1}}
\end{figure}
As the name suggests, the braid category is an organisational scheme for
braid patterns - for example in hair. The object class is the skeletal one
for this problem - the set $\N_0$ of possible numbers of hairs on a head.
The key notion is the
appropriate 
notion of equivalence on concrete braids, so that a morphism
is an equivalence class.
The 
picture
in Figure~\ref{fig:1}
illustrates this by showing concrete braids in the same class.
(In case the reader imagines that this construction is somehow frivolous, we
note that organisational schemes essentially identical to this are vital in
making sense of the observable world --- see for example \cite{LeinaasMyrheim77}.)


The category composition comes from thinking of transitions from one braid
class to another as themselves being braid classes.
This is indicated in our schema by the pictures
{in  
Figures \ref{fig:2},\ref{fig:3}},
along with the monoidal composition in Figure~\ref{fig:4}.
For corresponding formalism see for example Birman \cite{Birman} and
Mac~Lane \cite{MacLane}.
}
{Note the category is `diagonal' (endomorphisms only).}

\ignore{{
\befff{}{
\[
\includegraphics[width=7.4\ccm]{images/Braid-Step.jpg}
\]
\[
\includegraphics[width=8.9\ccm]{images/man1.jpeg}
\hspace{-.351in}
\includegraphics[width=4.84\ccm]{images/braid1.jpeg}
\]
}
}}


\befff{}{
\begin{figure}
\[
\includegraphics[width=6.4\ccm]{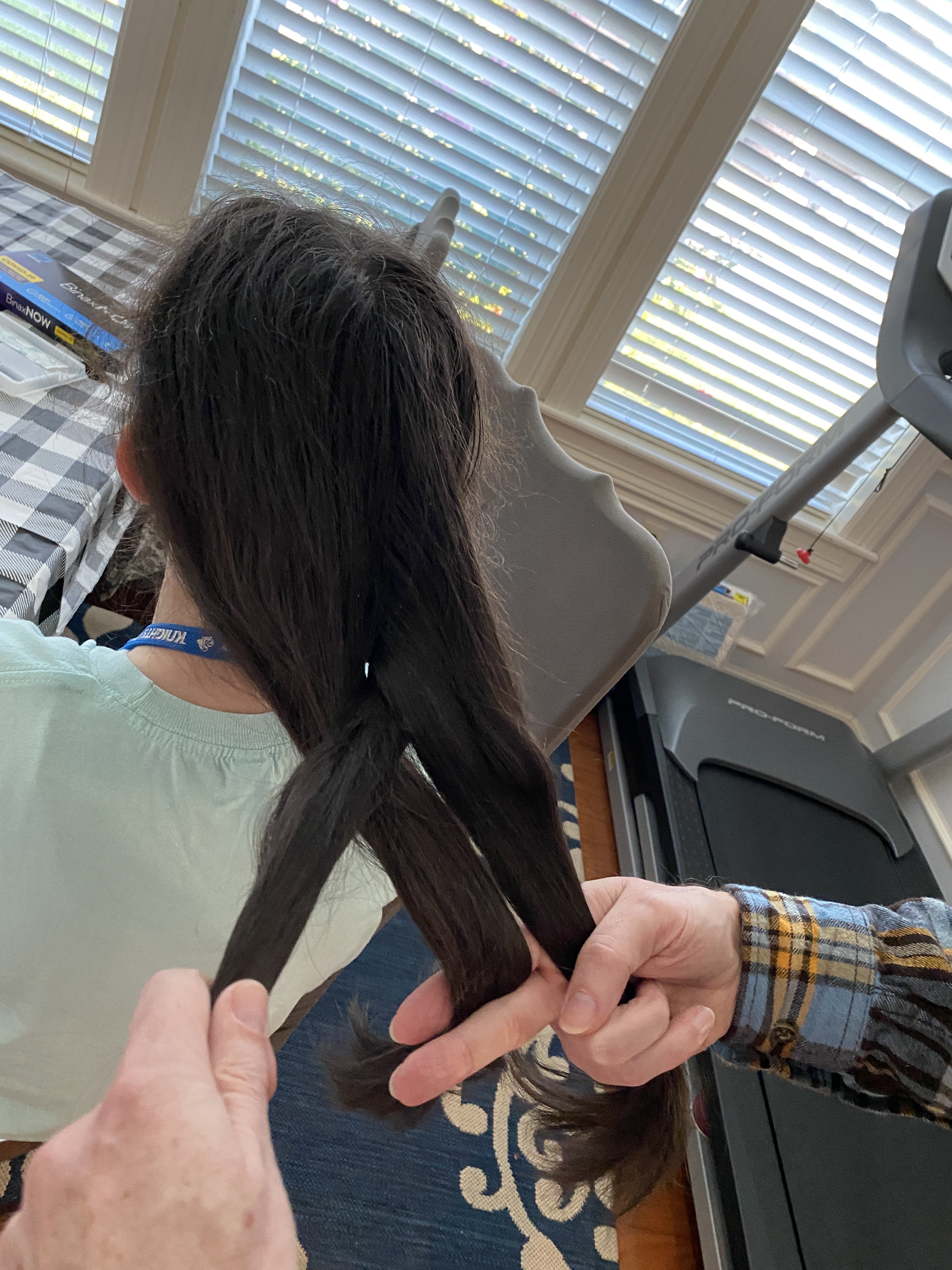}
\]\caption{A composition of morphisms \label{fig:2}}
\end{figure}
\begin{figure}
\[
\includegraphics[width=4.84\ccm]{images/IdentinBbig.jpg}
\hspace{.0351in}
\includegraphics[width=4.84\ccm]{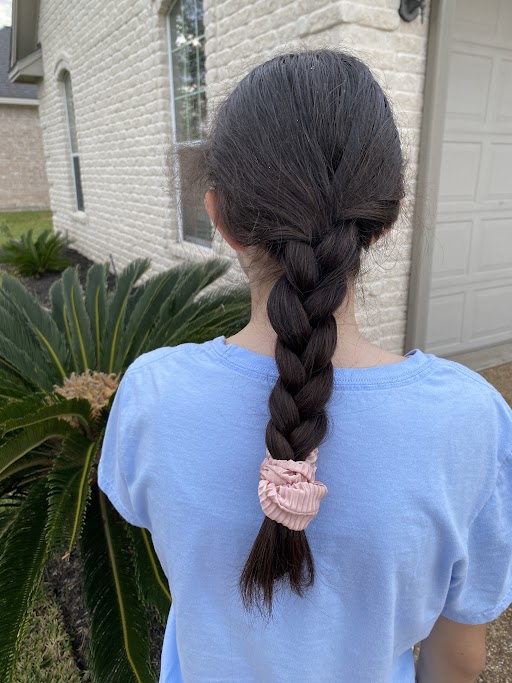}
\]
\caption{Before and after a sequence of compositions \label{fig:3}}\end{figure}
}

\beff{Monoidal composition}{
\begin{figure}
\[
\includegraphics[width=4\ccm]{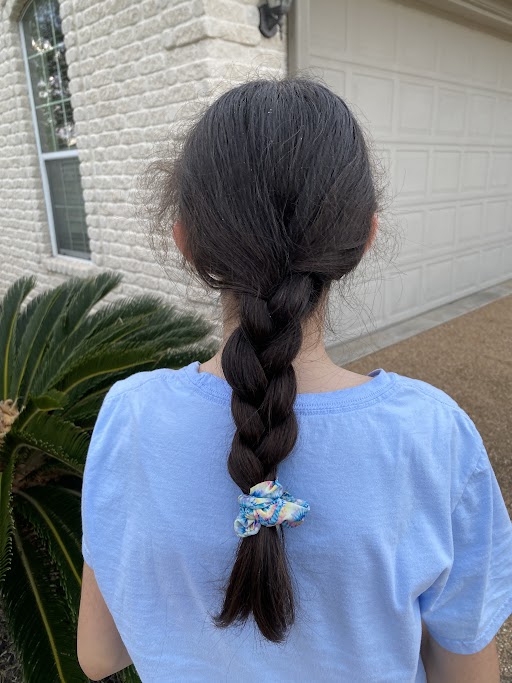}
\raisebox{1.05in}{$\;\;\otimes\;\;$}  
\includegraphics[width=4\ccm]{images/xB.jpg}
\raisebox{1.05in}{\;\; =\;\;\;}
\includegraphics[width=4\ccm]{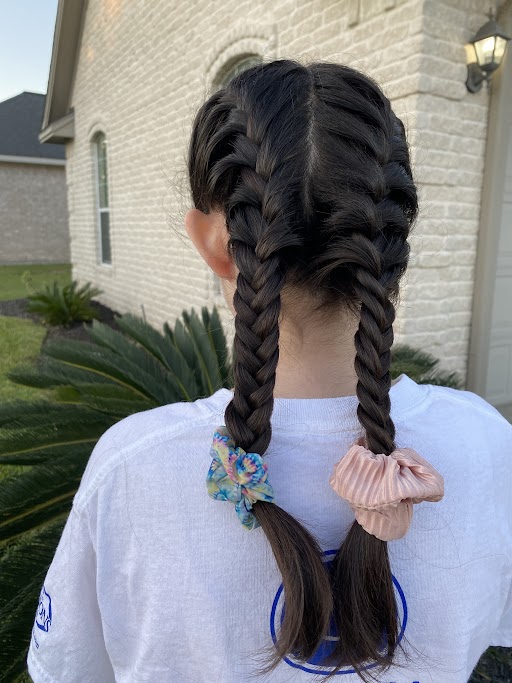}
\]
\caption{The monoidal product on morphisms \label{fig:4}}
\end{figure}
}

\fr{R-matrix}{{

The braid category
$\Bcat$
is monoidally generated by the elementary braids
$\sigma$ and  $\sigma^{-1}$  $ \in \Braid(2,2)$.
Thus a monoidal functor
\[
F : \Braid \rightarrow \Match^N
\]
or indeed from $\Braid$ to any target,
is determined by the image $F(\sigma)$.

}}




\module{notation-braid}





\ignore{{
Here $B_n$ denotes the $n$-string braid group. There are various
alternative definitions defining $B_n$ up to isomorphism.
Informally elements of $B_n$ are braidings of an ordered set of $n$
strings. 

Note that given an ordered pair $(b,b') \in B_n \times B_m$ there is an element 
$b\botimes b' \in B_{n+m}$ obtained by `side-by-side concatenation'.
(Strictly speaking some choices are needed in this construction 
- with details depending on the definition of $B_n$ - 
upon which we shall not dwell here.)
}}

Our monoidal category $\Bcat$
has its geometrical definition
above.
But in order to verify functors it is  
useful to recall a finite presentation.

A  {\em{\integ}} monoidal category is one whose object monoid is $(\N,+)$
  (confer \cite{MacLane2}).

\newcommand{\vsm}{\vs_-}

\begin{definition} \label{de:Bcat0}
\ignore{{
  {$\;$} \\ $\bullet$ 
  The {\em ordinary braid category} $\Bcat_{-}$ is the 
category
with
object set $\Nn$  and
$\Bcat_{-}(n,n) = B_n$, and 
$\Bcat_{-}(n,m)=\emptyset$ otherwise
(see e.g. \cite[XI.4]{MacLane}).
\\
\item $\bullet$
$\Bcat$ denotes  the {\em braid category}:
an \integ\ 
strict monoidal category
with
$\Bcat(n,n) = B_n$,
$\Bcat(n,m)=\emptyset$ otherwise, and monoidal composition $\otimes$
via side-by-side concatenation of suitable braid representatives
(see e.g. \cite[XI.4]{MacLane}).
\\
\item $\bullet$ 
}}
%
  The \integ\ strict monoidal category
$\Bcat'$ is that presented by  generators $\vs$
and $\vsm \in \Bcat'(2,2)$; and relations as follows.
We write
$1_n$ for the identity morphism in $\Bcat'(n,n)$;
$\vsigma_1 = \vsigma\otimes 1_1$;
$\vsigma_2 = 1_1\otimes \vsigma$
in $\Bcat'(3,3)$.
Relations are:
\beq \label{eq:inverse}
\vsigma \; \vsm 
\;=\; 1_2 \; =\;  
{\vsm}  \; \vsigma
\eq
\beq \label{eq:braidx}
\;\;\;
\vsigma_1 \vsigma_2 \vsigma_1 = \vsigma_2 \vsigma_1 \vsigma_2
\hspace{.3in}\mbox{ that is }\hspace{.34in}
\braidrel{\vsigma}{1_1}.
\eq
\end{definition}


\ignore{
Recall that category
  $\Bcat$ is monoidally generated by the positive braid
    $\sigma\in\Bcat(2,2)$ and its inverse (denoted $\sigma^{-1}$).
}
Via Artin's local presentation \cite{Artin:1925} we have the following.

\begin{proposition} \label{pr:Bcat1}
A presentation for $\Bcat$ is given by $\Bcat \cong \Bcat'$
where $\sigma \mapsto \vs$.
\qed

\end{proposition}

\module{mat}

\newcommand{\aaam}{{\mathbf a}} 
\newcommand{\bbbm}{{\mathbf b}} 

\section{Target categories and representation theory generalities} \label{ss:target}

`Charge-conserving' representation theory
as developed here
has several novel features
compared to, for example, classical Artinian representation theory.
Artinian representation theory --- representation theory of
finite-dimensional algebras over an algebraically closed field $k$ ---
can be understood as inheriting features from the appropriate
target category $\Vect_k$
(cf. for example \cite{Pareigis} and also \cite{Bondarenko}).

\newcommand{\VVV}{\mathsf{V}}  

In as much as representation theory of some algebraic structure $A$ say
is the study of the appropriate functor category between $A$ and the target $\VVV$
say, so  representation theory tends to inherit certain properties 
from $\VVV$.
For example if $A$ is a finite group, regarded as a category with one object,
then the functor category
$\functor(A,\Vect_k)$ inherits abelianness from $\Vect_k$,
from which we derive the Jordan--Holder and Krull-Schmidt properties.
These properties tell us that classification of representations should
proceed through classification of isomorphism classes of irreducible representations
and the 
decomposition series of indecomposable projective representations.
In the Artinian setting the theorems appear directly and the categorical
perspective is just nice abstract nonsense.
In more `rigid' settings such as ours, it can be proactively helpful.

Before proceeding the reader will therefore need to consider the definition and
properties of $\Match^N$ that we discuss in
Section~\ref{ss:MatchN}  
{\em et seq}.
We will return to discuss
the corresponding
analogues of Jordan--Holder, irreducibility,
isomorphism etc in
Section~\ref{ss:transversal}.

\subsection{On the target category $\Match^N$} 
\label{ss:MatchN}

\newcommand{\Nm}{N}
\newcommand{\Mato}{\overline{\Mat}}
\newcommand{\RFrip}{\overrightarrow{R}}
\newcommand{\ootimes}{\overline{\otimes}}
\newcommand{\flipp}[1]{\overline{#1}} 

$\;$

\mdef If $A$ and $B$ are matrices, let $A \otimes B$ denote the 
{Ab}-convention Kronecker product
(see e.g. \cite[Ch.3]{Murnaghan});
and $A \ootimes B = B \otimes A$ the
aB-convention
Kronecker product.

Here $\Mat$ is the monoidal category of matrices over a given
commutative ring (and $\Mat_k$ the case over commutative ring $k$),
with
object monoid $(\Nn, \times)$ and
tensor product on morphisms given by  the
aB Kronecker product.
We  write $\Mato$ for the monoidal category with
$\ootimes$.



\mdef
For $N \in \Nn$ 
the monoidal subcategory $\Mat^\Nm$ of $\Mat$ is that generated  by a
single object $\Nm$.  
Then the object monoid  $(\Nm^{\Nn_0}, \times)$ 
is isomorphic to $(\Nn_0,+)$ in the natural way,
so that $\Mat^N(a,b)=\Mat(N^a, N^b)$.

{A monoidal category is `natural' if the object monoid is freely generated
  by a single object, denoted 1 (thus $\Mat^N$ is natural). 
  Note that if a monoidal category $C$ is natural
  then every monoidal functor $F:C \rightarrow \Mat$ factors through
$\Mat^N$ for some $N=F(1)$.}

\newcommand{\mub}{\overline{\mu}}
\newcommand{\cC}{{\mathsf{C}}} 
\newcommand{\GL}[2]{GL_{#1} #2}   

\mdef
Let $\mu$ be an  invertible $N\times N$ matrix.
We write $\mub$ for the inverse here, to neaten superscripts.
Note that the `local base change' given
(on object 1 in the new  $(\Nn_0,+)$ labelling) by
$\mu$ induces a monoidal functor
\[
\cC^\mu : \Mat^N \rightarrow \Mat^N
\]
by taking $T \in \Mat^N(a,b)$ to
$( \mu^{\otimes a} )\; T\; (\mub^{\otimes b} )$.
Thus $\cC^\nu \circ \cC^\mu = \cC^{\nu\mu}$. 
{That is, the object monoid of 
  functor category
  $\functor(\Mat^N,\Mat^N)$
  contains a realisation of $\GL{N}{} = \GL{N}{(k)}$
  (see e.g. \cite{Martin08a}).}
\ignore{{
(Verification: for $U \in  \Mat^N(a',b')$ we have
$T \otimes U \in  \Mat^N(a+a',b+b')$ and 
$\mu^{\otimes a+a'}(T\otimes U)\mub^{\otimes b+b'} =
(\mu^{\otimes a} \;T\; \mub^{\otimes b}) \otimes
(\mu^{\otimes a'} \;U\; \mub^{\otimes b'})$
by Kronecker identities.)
}}

Thus for any monoidal
functor $F:\Bcat \rightarrow \Mat^N$, say, and any
$\mu$  
there is another monoidal functor
$
F^\mu :\Bcat \rightarrow \Mat^N
$
given by applying the `diagonal' action:
$F^\mu = \cC^\mu \circ F$.  

\newcommand{\set}[1]{\ul{#1}} 
\newcommand{\bbb}[1]{\mathbeef_{#1}}  
\newcommand{\bbbb}[1]{\mathbeef_{\underline{#1}}}  
\newcommand{\Setcat}{{\mathsf{Set}}}  

\mdef  \label{de:Set}
For $N \in \N_0$ let  
$\set{N} = \{ 1,2,...,N \}$.
Let $\Setcat$ be the category of sets;
and $\Setcat^i(S,T) \subset \Setcat(S,T)$ the subset of isomorphisms.
Let $\Sym_N$ denote the symmetric group.
Thus $\Sym_N = \Setcat^i(\set{N},\set{N})$. 


For $S$ a set, 
let $S^n$ denote
the set of sequences of length $n$ of elements from $S$.
Equivalently,
$S^n =\Setcat(\set{n},S)$.
If $S$ is a set of symbols 
then $S^n$ can be represented as
the set of words
$s_1 s_2 \cdots s_n$
of length $n$ in $S$.
Let $S^\ast$ denote the set of all words and $S^\star$ all non-empty words:
\[
S^\star = \bigsqcup_{n \in \N} S^n  \hspace{1cm}\mbox{ and }\hspace{1cm}
S^\ast = \bigsqcup_{n \in \N_0} S^n
\]


\mdef  
We may label the row/column index for object $N$ in $\Mat$
(object 1 in $\Mat^N$)
by
$
\bbbb{N} = \{\bbb{1},\bbb{2},...,\bbb{N} \}
$
(or as `charges' $\{ +,- \}$ if $N=2$).
Note that the `standard ordering' is arbitrary.
Thus for $\omega\in \Sym_N$ 
there is a monoidal functor $f^\omega : \Mat^N \rightarrow
\Mat^N$
given by $f^\omega = \cC^{\rho(\omega)}$,
where $\rho(\omega)$ is the perm matrix of $\omega$.
{This is the obvious copy of
$\Sym_N$ in our $\GL{N}{}$.}

Remark: Choosing index sets allows us to use the language of linear algebra.
For a matrix $\alpha\in\Mat(N,M)$ we may write
$
\alpha \bbb{i} = \sum_{j\in \underline{M}} \alpha_{ij} \bbb{j}
$
to give an individual row.
(We call this the `(right) action' of $\alpha$, as if it were a linear
operator.
Think of $\bbb{i} \in \Mat(M,1)$ (an elementary column vector).)

In $\Mat^2$, say, 
the object $2\otimes 2$ has `tensor' index set
$\{ \bbb{1}\otimes \bbb{1},
\bbb{2}\otimes \bbb{1}, \bbb{1}\otimes \bbb{2}, \bbb{2}\otimes \bbb{2} \}$
(which we may abbreviate to
$\{ \ket{11},\ket{21},\ket{12},\ket{22} \}$,
or even just
$\{ 11,21,12,22 \}$)
and so on.
Thus the tensor index set for rows of $\alpha\in \Mat^N(M,M)$, say, is
given by the set $\ul{N}^M$ of words in $\{1,2,...,N\}$ of length
$M$.
Note that $\Sym_M$ acts  on this set by place permutation.


\mdef
Fix $N$. 
We say a matrix $\alpha$ in
$\Mat^N(M,M)$
is {\em charge-conserving} when 
$\alpha_{ij} \neq 0$ implies $j=\omega(i)$ for some $\omega\in \Sym_M$.
That is, $\alpha$ is a (reordered) block matrix with blocks given by
`charge'
(orbits of the $\Sym_M$ action on $\set{N}^M$,
which are, note,
indexed by compositions of $M$).

\begin{lemma} \label{lem:chacon1}
  Fix $N$  and $k$ (by default $k=\C$).
\\
  (I) The charge-conserving matrices form a (diagonal)
monoidal subcategory $\Match^N$ of $\Mat^N$.
\\
(II) Let $w \in \Sym_N$. The functor $f^w$ restricts to a monoidal
functor on $\Match^N$:
\beq \label{eq:wMatc}
f^w : \Match^N \rightarrow \Match^N .
\eq
\ppmx{If $\mu\in GL(N)$ is diagonal (in the diagonal maximal torus) then
  $\cC^\mu$ also restricts to a monoidal functor on $\Match^N$,
  but it is trivial.}

\noindent
(III) For each $M \leq N$ and each injective function
$\psi: \ul{M} \rightarrow \ul{N}$ there is a monoidal functor
$f^\psi : \Match^N \rightarrow \Match^M$ given
on morphisms 
by
$f^\psi (\aaam)_{vw} = \aaam_{\psi v , \psi w}$
(here $\aaam \in \Match^N(m,m)$ say).
\end{lemma}
\proof
(I) The closure of composition follows in the manner of closure for
block diagonal matrices of given block shape,
so it
remains to show the monoidal property.
Consider $\aaam \in \Mat^N(m,m)$ and
$\bbbm \in \Mat^N(m',m')$.
Note that the row index set for $\aaam \otimes \bbbm$ is the concatenation 
$\ul{N}^{m+m'}$
(as noted, it is a convention choice for the overall order).
Thus if $i$ is a row index for $\aaam$ (an $m$-tuple) and $I$ is a row
index for
$\bbbm$)
(an $m'$-tuple) then $iI$ is a row index for $\aaam \otimes \bbbm$.
We have
\[
(\aaam \otimes \bbbm)_{iI,jJ} = \aaam_{i,j} \bbbm_{I,J} 
\]
Now note that if $j=\omega(i)$ and $J=\omega'(I)$ for some perms
$\omega\in \Sym_m$ and $\omega'\in \Sym_{m'}$ then
$
jJ = \omega\otimes\omega'(iI).
$
\\
(II,III) Clear from the construction.
\qed


\mdef
    {Note that if $\mu\in
      \GL{N}{}$ is diagonal (in the diagonal maximal torus) then
  $\cC^\mu$ also restricts to a monoidal functor on $\Match^N$,
  but it is trivial.}
  Note that the functors in (II) do not exhaust  
  the set of isomorphisms
  in the object set of $\functor(\Match^N , \Match^N)$ in general.
  We can classify such functors as `local' (restrictions of $\cC^\mu$);
  `decomposable' (asigning a possibly distinct invertible to each tensor factor
  --- see below);
  and `entangled' (everything else).
  \\
  On the other hand note that (III) does not extend even to a functor on $\Mat^N$
  in general.
  \\
  Note that $\Match^N$ also inherits linearity from $\Mat$.
Here we will use this only for convenience.

\newcommand{\kett}[1]{\ket{\underline{#1}}}
\newcommand{\ssigma}{\varsigma}
\newcommand{\sigaa}{\ssigma\otimes 1}
\newcommand{\sigbb}{1\otimes\ssigma}
\newcommand{\siga}{\ssigma_1}
\newcommand{\sigb}{\ssigma_2}
\newcommand{\ag}{a_{11}}
\newcommand{\af}{a_{22}}
\newcommand{\aaa}{a_{12}}
\newcommand{\aab}{b_{12}}
\newcommand{\aac}{c_{12}}
\newcommand{\aad}{d_{12}}
\newcommand{\alist}{\ag, a,b,c,d,\af}
\newcommand{\varSig}{\varSigma} 



\mdef
  Let $\Ccat$ be an \integ\ monoidal category generated by
  $\varSig \subset \Ccat(2,2)$. 
A level-$N$ monoidal functor
$\FFF : \Ccat \rightarrow \Mat$
is one
taking object 1 (in $\Ccat$) to $N$ (in $\Mat$).
The functor is
{\em charge-conserving}
(the image lies in $\Match^N$)
if every generator $\ssigma\in \varSig$ 
acts so that 
$ \ssigma \kett{ij} $ lies in the span of $\kett{ij}, \kett{ji}$.
%
In other words we have the following
for $\ssigma \in \varSig$
(e.g. the positive braid in $\Bcat(2,2)$):
%
    {{
\beq \label{eq:abcd12}
\ssigma \kett{11}  = a_{11} \kett{11} ,
\hspace{.131in} 
\ssigma \kett{12} = a_{12} \kett{12} + b_{12} \kett{21}
,\hspace{.13in}
\ssigma \kett{21} = c_{12} \kett{12} + d_{12} \kett{21} ,
\hspace{.131in} 
\ssigma \kett{22}  = a_{22} \kett{22} 
\eq
}}
\ignore{{
\newcommand{\ag}{a_{11}}
\newcommand{\af}{a_{22}}
\newcommand{\aaa}{a_{12}}
\newcommand{\aab}{b_{12}}
\newcommand{\aac}{c_{12}}
\newcommand{\aad}{d_{12}}
\newcommand{\alist}{\ag, a,b,c,d,\af}
}}%
\ignore{{
\begin{eqnarray}  \label{eq:Fs}
\ssigma \kett{11} & = & \ag\; \kett{11}
\\
\ssigma \kett{12} &= & \aaa \kett{12} + \aab \kett{21}
\\
\ssigma \kett{21} &= & \aac \kett{12} + \aad \kett{21}
\\ \label{eq:Fs4}
\ssigma \kett{22}  &= & \af\; \kett{22} 
\end{eqnarray}
}}%
%
for some $a_{11},a_{12},b_{12},c_{12},d_{12},a_{22} \in k$;
and similarly replacing
12 with each $ij$, with $i<j$.
\ignore{{
for some
$g=\ag, \; a=a_{12},\; b=b_{12},\; c=c_{12}, \; d=d_{12},\; f=\af$.
}}
\ignore{{
Note (for later use) that the only smooth variations of these
parameters that fix the spectrum of the given action, just of
$\sigma$, are of the form
\beq \label{eq:bcbxcx}
(b,c) \mapsto (bx,cx^{-1}) . 
\eq
}}



\mdef 
Calculus for operators  
in the monoidal category,
such as $\ssigma_1 = \ssigma\otimes 1$,
is exemplified by:  
\[
\siga\kett{112} = (\ssigma\otimes 1) \kett{112} = \ag \kett{112} ,
\hspace{1.015in} 
\ssigma_2 \kett{112} =
(1\otimes\ssigma) \kett{112} = \aaa \kett{112} + \aab \kett{121}
\]
\[  (\ssigma\otimes 1) \kett{121} = \aaa \kett{121} + \aab \kett{211}  
\hspace{1.52in} 
(1\otimes\ssigma) \kett{121} = \aac \kett{112} + \aad \kett{121}   \]
\[
(\ssigma\otimes 1) \kett{211} = \aac \kett{121} + \aad \kett{211}
\hspace{1.52in}
(1\otimes\ssigma) \kett{123} = a_{23} \kett{123} + b_{23} \kett{132}
\]
\ppmx{
\[ 
\siga\sigb\kett{123} = 
\siga( a_{23} \kett{123} + b_{23} \kett{132})
=  a_{23} (a_{12} \kett{123} +b_{12} \kett{213})
 + b_{23} (a_{13} \kett{132} +b_{13} \kett{312})
\]
\[ 
\sigb\siga\kett{123} = 
\sigb( a_{12} \kett{123} + b_{12} \kett{213})
=  a_{12} (a_{23} \kett{123} +b_{23} \kett{132})
 + b_{12} (a_{13} \kett{213} +b_{13} \kett{231})
\]
}%
\[ 
\siga\sigb\siga\kett{123} = 
  a_{12} (a_{23} (a_{12}\kett{123} +b_{12}\kett{213})
        +b_{23} (a_{13}\kett{132} +b_{13}\kett{312}))
\] \beq \label{eq:cond001} \hspace{1in}
 +b_{12} (a_{13} (\aad  \kett{213} +\aac  \kett{123})
        +b_{13} (a_{23} \kett{231} +b_{23} \kett{321})) 
 \eq
\[
 \sigb\siga\sigb\kett{123} =
  a_{23} (a_{12} (a_{23} \kett{123} +b_{23} \kett{132})
        +b_{12} (a_{13} \kett{213} +b_{13} \kett{231}))
\]
\beq \label{eq:cond002} \hspace{1in}
 +b_{23} (a_{13} (d_{23} \kett{132} +c_{23} \kett{123})
        +b_{13} (a_{12} \kett{312} +b_{12} \kett{321}))   
\eq

\mdef
For later reference,
the conditions for the actions in
(\ref{eq:cond001}), (\ref{eq:cond002})
to be equal are thus
\beq \label{eq:N2}
a_{12}^2 a_{23} +a_{13} b_{12} c_{12} \stackrel{123}{=} a_{12} a_{23}^2 + a_{13} b_{23} c_{23} ,
\eq
\beq \label{eq:cond004}
( a_{12} a_{13} - a_{23} a_{12} - a_{13} d_{23} ) b_{23} \stackrel{132}{=}0 ,
\hspace{1in} 
( a_{12} a_{23}  +  a_{13} d_{12} - a_{23}  a_{13} ) b_{12} \stackrel{213}{=}0
\eq
Actions on any $\kett{ijk}$ can be computed similarly, or by applying $\Sym_N$.
For $\kett{112}$ and its orbit:  
\[
 (\ssigma\otimes 1)(1\otimes\ssigma)  \kett{112}
=  (\ssigma\otimes 1) ( \aaa \kett{112} + \aab \kett{121} )
=  ( \ag\aaa \kett{112} + \aab (\aaa \kett{121} + \aab \kett{211} ) )
\]
\ignore{{
\[
(1\otimes\sigma)  (\sigma\otimes 1)(1\otimes\sigma)  \kett{112}
= (1\otimes\sigma)  ( \ag\aaa \kett{112} + \aab (\aaa \kett{121} + \aab \kett{211} ) )
\hspace{1in}
\]}}
\[ 
\sigb\siga\sigb\kett{112}
=  ( \ag\aaa (\aaa\kett{112} + \aab\kett{121})
+ \aab (\aaa (\aad\kett{121}+\aac\kett{112}) + \ag\aab \kett{211} ) )
\]
\ppmx{
\[
 (\ssigma\otimes 1)(1\otimes\ssigma)  \kett{121}
=  (\ssigma\otimes 1) ( \aac \kett{112} + \aad \kett{121} )
=  ( \ag\aac \kett{112} + \aad (\aaa \kett{121} + \aab \kett{211} ) )
\]
}%
\ignore{{
\[
(1\otimes\sigma)  (\sigma\otimes 1)(1\otimes\sigma)  \kett{121}
= (1\otimes\sigma)  ( \ag\aac \kett{112} + \aad (\aaa \kett{121} + \aab \kett{211} ) )
\hspace{1in}
\] }}
\[ 
\sigb\siga\sigb\kett{121}
=  ( \ag\aac (\aaa\kett{112} + \aab\kett{121})
+ \aad (\aaa (\aad\kett{121}+\aac\kett{112}) + \ag\aab \kett{211} ) )
\]
\ppmx{
\[
(1\otimes\ssigma)  (\ssigma\otimes 1) \kett{121}
=  (1\otimes\ssigma) ( \aaa \kett{121} + \aab \kett{211} )
=  ( \aaa(\aac \kett{112} + \aad \kett{121} ) + \ag \aab \kett{211} ) )
\]
}%
\ignore{{
\[
 (\sigma\otimes 1)(1\otimes\sigma)  (\sigma\otimes 1) \kett{121}
= (\sigma\otimes 1)  ( \aaa(\aac \kett{112} + \aad \kett{121} ) + \ag \aab \kett{211} ) )
\hspace{1in} \] }}
\[
\siga\sigb\siga\kett{121}
= \aaa( \ag\aac \kett{112} + \aad (\aaa\kett{121} + \aab\kett{211})
     + \ag \aab (\aad\kett{211} +\aac\kett{121}) )
\]
Thus for example a condition
$
\FFF(\ssigma_1 \ssigma_2 \ssigma_1)  =\FFF(   \ssigma_2 \ssigma_1 \ssigma_2)
$
would imply
(with $\ul{12}\leadsto \ul{ij}$ for all $i<j$ etc.)
\beq \label{eq:N2i}
(\ag^2 \aaa - \ag\aaa^2 - \aab\aac{a_{12}}) \kett{112} =0 ,
\hspace{1.02in} 
\aad\aaa\aab  
\kett{121}=0
\eq
\beq \label{eq:N2ii}
(
\aad\aaa\aac 
) \kett{112}  = 0 ,
\hspace{.31in}  
( 
\aad^2\aaa -\aaa^2\aad
) \kett{121}  = 0 ,
\hspace{.29in} 
\aaa\aad\aab 
\kett{211} = 0
\eq


{Noting Proposition~\ref{pr:Bcat1},
these cubics are the residue in this setting
of the 
intractable
problem described in the
Introduction.
While still daunting, they do now yield  
better to analysis.}




\module{repsN}

%


\ignore{{
\beff{}{{
Our `abcd' labeling conventions for entries in $F(\sigma)$ are as indicated here:

\[
\includegraphics[width=7\ccm]{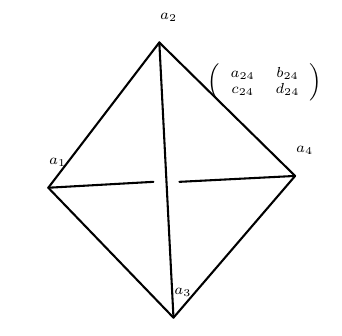}
\]
}}
}}


\beff{}{{

\begin{theorem} \label{th:constraints}
Let $\FF:\Bcat\rightarrow\Match^N$ be a level-$N$ charge-conserving
monoidal functor, thus determined by $F(\sigma)$.
The condition
$ 
(\sigma \otimes 1 )(1 \otimes \sigma)(\sigma \otimes 1 )
=
(1 \otimes \sigma)(\sigma \otimes 1 )(1 \otimes \sigma)
$ 
(from \ref{pr:Bcat1})
{is equivalent to}
the following
  constraints on entries in the $abcd$ form
(equation (\ref{eq:abcd12}))
of $F(\sigma)$:
\beq \label{eq:31xx0}
a_{12} (a_{11}^2 -a_{12} a_{11} -b_{12} c_{12}) =0,
\hspace{1cm}
a_{12} (a_{22}^2 -a_{12} a_{22} -b_{12} c_{12}) =0,
\eq \beq \label{eq:31xxx0}
\hspace{1cm}
a_{12} c_{12} d_{12} = 0 = a_{12} b_{12} d_{12} ,
\hspace{1cm}
a_{12} d_{12} (a_{12}-d_{12})=0,
\eq

\ignore{{
\beq \label{eq:38x}
k(f^2 -fk-ms)=0,
\hspace{1cm}
k(w^2 -kw -ms)=0,
\hspace{1cm}
ksv=kmv=kv(k-v)=0
\eq 
\beq \label{eq:311x}
v(-f^2 + fv + ms)=0,
\hspace{1cm}
v(w^2 -vw -ms)=0  
\eq
}}
%
\beq \label{eq:34xx0}  
c_{12} (d_{13} d_{23} -d_{12} d_{23} -a_{12} d_{13})=0,
\hspace{1cm}
c_{12} (- a_{13} a_{23} +a_{12} a_{23} +d_{12} a_{13})=0  
\eq
\ignore{{ 
\beq \label{eq:36xx0} 
\Or{c_{13}}{b_{13}}  (-a_{23} d_{13} +d_{12} a_{23} -d_{12} a_{13})=0,
\hspace{1cm}
\Or{c_{13}}{b_{13}}  (a_{13} d_{23} -a_{12} d_{23} +a_{12} d_{13})=0
\eq
\beq
\Or{c_{23}}{b_{23}}  (- a_{13} d_{23} -a_{12} a_{23} +a_{12} a_{13}) =0,
\hspace{1cm}
\Or{c_{23}}{b_{23}}  (d_{12} d_{23}+ a_{23} d_{13} -d_{12} d_{13})=0  
\eq
}}
%
\beq \label{eq:35x0} 
-a_{13} d_{23}^2 + a_{13}^2 d_{23} -a_{12} b_{23} c_{23}
               + a_{12} b_{13} c_{13} =0
\eq
\ignore{{
\beq
 -b_{13} a_{23} c_{13}
 +b_{12} c_{12} a_{23} -d_{12} a_{13}^2 +d_{12}^2 a_{13} =0  
\eq
\beq \label{eq:33x0}
-a_{13} b_{23}c_{23}- a_{12} a_{23}^2 +a_{12}^2 a_{23} +b_{12} c_{12} a_{13} =0 , \hspace{1cm} 
b_{13}c_{13}d_{23} -b_{12} c_{12} d_{23} +a_{12} d_{13}^2 -a_{12}^2 d_{13} =0
\eq 
\beq \label{eq:310x0}
d_{12}b_{23}c_{23}-a_{23}d_{13}^2 +a_{23}^2 d_{13} -d_{12}b_{13}c_{13} =0, \hspace{1cm}
d_{12}d_{23}^2 - d_{12}^2 d_{23} +b_{23}d_{13}c_{23} -b_{12}c_{12}d_{13} =0 
\eq
}}
  together with the images of these constraints which are their orbits
  under the action
{
of $\Sigma_N$ by $\omega a_{ij} = a_{\omega(i) \omega(j)}$ and so on
(where $a_{ji}=d_{ij}$ and $b_{ji}=c_{ij}$ is understood).}
\end{theorem}
}}

\beff{}{{
\noindent{\bf Proof.}
Observe that the YBE
can be verified in $\Mat^N(3,3)$
(since it can be formulated with three tensor factors; and the higher
versions simply contain copies of the same entries, by the tensor
construction).
This yields
equations 
across three indices,
 computed as in (\ref{eq:N2})-(\ref{eq:N2ii}).
\ppmx{
for example with $N=4$ the row of the matrix that must vanish
$F(\vsigma_1 \vsigma_2 \vsigma_1) -F( \vsigma_2 \vsigma_1 \vsigma_2)$
{(index notation as in Definition \ref{de:Bcat0})}.
with label $1\otimes2\otimes3$ (hereafter simply 123) is
\begin{align*} &
[0, 0, 0, 0, 
0, 0, a_{13} b_{12} c_{12} - a_{13} b_{23} c_{23} - a_{12} a_{23}^2  + a_{12}^2  a_{23}, 0, 
\\ & \hspace{2.82cm} 
0, a_{12} a_{13} b_{23} - a_{13} b_{23} d_{23} - a_{12} a_{23} b_{23} , 0, 0, 
0, 0, 0, 0,
\\ &
0, 0, a_{13} b_{12} d_{12} - a_{13} a_{23} b_{12} + a_{12} a_{23} b_{12}, 0, 
0, 0, 0, 0, 
0, 0, 0, 0, 
0, 0, 0, 0, 
\\ &
0, 0, 0, 0, 
0, 0, 0, 0, 
0, 0, 0, 0, 
0, 0, 0, 0, 
\\ &
0, 0, 0, 0, 
0, 0, 0, 0, 
0, 0, 0, 0, 
0, 0, 0, 0]
\end{align*}
--- column order is the usual lex order (so
111 112 113 114 121 122 123 124 then 131 132 ...).
Observe that the non-zero entries are at matrix entries given by 
$123\rightarrow 123$,
$123\rightarrow 132$,
$123\rightarrow 213$.
The point we need here is (from charge conservation) that for a
non-zero matrix entry the column index is at most a perm of the row
index.
Since the row index can contain at most three different colours (it
contains three colours), the column index contains the same three.
It follows that the level 3 case
produces representatives of all possible orbits of constraint.
}
\qed
}}



The following is an immediate consequence, but is a key fact so we record it as:
\begin{lemma}\label{lem:3 suffice} 
A level-$N$ charge-conserving operator $F(\sigma)$ defines a monoidal functor $F:\Bcat\rightarrow\Match^N$
if, and only if all of the ${N \choose 3}$ restrictions of $F$ to $\Match^3$ satisfy equations (\ref{eq:31xx0})-(\ref{eq:35x0}) together with the images of these constraints which are their orbits
  under the action of
  $\Sigma_3 \ppmx{\times \Zz_2?}$.
\qed
\end{lemma}





\module{Xlemma}

\ppmx{done: promote X-Lemma here. (but proof must remain after nec Theorem)}

\newcommand{\Xequivalence}{X-equivalence} 
\newcommand{\Xequivalent}{X-equivalent}  

\ppmx{[maybe promote following lemma to immediately after
(\ref{eq:cond001}-\ref{eq:N2ii})?]}

\mdef {\bf{Lemma}}. \label{lem:Xeq} (X Lemma)
{\em
Let $\FF:\Bcat\rightarrow\Match^N$ be a braid representation, hence given by
$\FF(\sigma)$.
(I)
Let $X$ be a diagonal invertible matrix in $\Match^N(2,2)$.
Then $\FF'(\sigma)= X \FF(\sigma) X^{-1}$ also gives a representation.

\newcommand{\simm}{\stackrel{-}{\sim}}

\noindent
(II)
%
The relation: $\FF\simm\FF'$ if $\FF'(\sigma)= X \FF(\sigma) X^{-1}$ for some $X$,
is an equivalence.
\\ The relation $\simm$ is here called 
{\em \Xequivalence},
{and written $\stackx$}.
}

%

\smallskip

\noindent
{\em Proof}.
{Noting presentation \ref{pr:Bcat1}, this }
follows from how $b_{ij},c_{ij}$ appear in
(\ref{eq:N2}-
\ref{eq:N2ii}).
{See also (\ref{proof:X}).}
\qed

The analogue of
\ref{lem:Xeq} 
in classical representation theory is
stronger, and is elementary
and subject independent. 
Here in the monoidal setting
we see that
it is not elementary
or general.






\module{mat}


\mdef
We have 
a monoidal functor between
monoidal categories
$
\RFrip: \Mat^N \rightarrow \Mato^N
$
given by the identity on objects
and on morphisms by
$\RFrip(A)_{v,w}  = A_{\flipp{v},\flipp{w}}$
where $\flipp{w}$ is the reverse word.
\ignore{{
, and
such that for each pair $A,B$   
$\RFrip (A \otimes B) = A \ootimes B = B \otimes A$.
}}
\ignore{{
\ppm{-well defined!? - it is not quite right to say `given' here.
  we have to check that no issues arise from the
  fact that different pairs $A',B'$ might give
  $A\otimes B=A'\otimes B'$...
  We claim, however, that $A\otimes B = A'\otimes B'$ implies
$B\otimes A = B'\otimes A'$, so we are good. ...ok?}
(they are identical as categories,
but $A \ootimes B = B \otimes A$).
}}


\newcommand{\cF}{{\mathsf F}}
\newcommand{\Flip}{\overleftrightarrow{\cF}} 
\newcommand{\Frip}{\overrightarrow{\cF}} 
\newcommand{\PPP}{{\mathtt P}}  


\mdef \label{de:flip1}
The `lateral flip' inner automorphism on an individual braid group:
\beq
\Flip: B_n \rightarrow B_n
\eq
(see e.g. \cite[\S5.7.2]{Martin91} or \cite{Hietarinta:1992})
corresponds, on the 
braid generators, to
$\sigma_i \mapsto \sigma_{n-i}$.
Note that this extends to an automorphism of
$\Bcat$, except that the monoidal product convention is reversed,
so there we denote it $\Frip : \Bcat \rightarrow \overline{\Bcat}$.
Similarly we have  $\Frip : \overline{\Bcat} \rightarrow \Bcat$.
Combining with $\RFrip$ we get an involutive action on the functor category
$\functor(\Bcat,\Mat)$,
with generator that we will denote $\PPP$.
$\;$
Note that this fixes the target subcategory $\Match^N$.

\module{braid-mat}

{{

\mdef\label{restriction-and-symmetries}
In ordinary representation theory organisation is provided by
intertwiners/module morphisms, which are the morphisms in the representation
category. In our case this view lifts so that intertwiners between
representation-functors  
are natural transformations.
The various key principles of representation theory have lifts to our setting,
but these lifts have novel features,
as in Lemma~\ref{lem:chacon1}.
And then there are
features that depend on specific properties of the 
`subject-algebra'   
(in our case $\Braid$)
{such as (\ref{lem:Xeq}), and the following}.

\ignore{{
The following observations will play an important role in what follows:
\begin{itemize}
\item If $F : \Braid \rightarrow \Match^N$
is a monoidal functor  
then by restricting to any subset $J\subset \ul{N}$ we obtain another monoidal functor $F_J:\Braid\rightarrow\Match^{|J|}$.
\item  The aB Kronecker product convention choice is arbitrary, so that there is $\Zz_2$ action on the monoidal functors $F:\Braid \rightarrow \Match^N$ given by choosing the Ab convention.   Later we denote this symmetry by $P$ as it is implemented by conjugating by the usual `flip' operator, often denoted $P$.
\item The symmetric group $\Sigma_N$ acts on the category of monoidal functors $F : \Braid \rightarrow \Match^N$ by permuting the basis  $\ul{N}$.
\item
\ppm{The $\Zz_2$ action from \ref{de:flip1} and
$\Sigma_N$ action from \ref{lem:chacon1}(II)
commute, giving an action of $\Zz_2 \times \Sigma_N$.}
\end{itemize}
Notice that while $GL(N)$ acts on the category of monoidal functors $F : \Braid \rightarrow \Mat^N$ by (local) change of bases, these do not respect the charge-conserving property in general.  Thus the category
of monoidal functors to $\Match^N$ seems to have less symmetry.  On the other hand we shall see later that there are continuous symmetries that the space of monoidal functors to $\Match^N$ have that are not present for monoidal functors
$F : \Braid \rightarrow \Mat^N$ in general.

}}

{The $\Zz_2$ action from \ref{de:flip1}
fixes the target $\Match^N$.  
It
and the
$\Sigma_N$ action induced from \ref{lem:chacon1}(II)
commute, giving an action of $\Zz_2 \times \Sigma_N$
on the category of monoidal functors $F : \Braid \rightarrow \Match^N$.}

}}

\subsection{
{Geometrical setup for}
functors to Match categories} $\;$
\ignore{{

\fr{Target category $\Match^N$}{{

Fix a commutative ring $k$ (we take $k = \C$).

\mdef \label{de:Mat}
Here $\Mat = \Mat_k$ denotes  
the monoidal category of matrices
    and aB-convention Kronecker product.
Label rows of object $N$ by $\ul{N} = \{1,2,...,N\}$.

\mdef
    $\Mat^N$ is the monoidal subcategory generated by object $N$
    (which is renamed as object 1 in $\Mat^N$).
    $\Mat^N(m,n) = \Mat(N^m, N^n)$

\mdef
    Consider $R \in \Mat^N(m,n)$ with entries
    \\
    $\langle w | R | v \rangle \in k$ 
    for $w \in \ul{N}^m$ and $v \in \ul{N}^n$
    \\
    $R \in \Mat^N(m,m)$ is \cc\ if $\langle w | R | v \rangle \neq 0$ 
    implies $w$ a permutation of $v$.
    \\
    These matrices form a subcategory, $\Match^N$. 
}}

}}
\fr{R-matrix}{{

The focus of this paper is to classify charge-conserving solutions to the YBE,
i.e., monoidal functors
$F:\Braid \rightarrow \Match^N$;
and provide  
the
algorithm for producing them.

The braid category is monoidally generated by the elementary braids
$\sigma,\sigma^{-1} \in \Braid(2,2)$.
Thus a monoidal functor
$
F : \Braid \rightarrow \Match^N
$
is determined by the image $F(\sigma) \in \Match^N(2,2)$.
\footnote{
{ The braid category is also
equivalent to 
the free braided monoidal category on
a single object \cite{JoyalStreet93}, but the cost of this extra level of 
wiring 
makes it less of a boon than it is a distraction
in our immediate setting.
}}

}}
{{

\mdef \label{pa:sparse}
Recall that the rows (and columns) of a matrix $M\in\Mat^N(2,2)$
may be indexed by
ordered pairs $(i,j)$ with $i,j \in \ul{N} = \{ 1,2,...,N \}$.
Thus matrix entries are $M_{ij,kl}$.
A matrix $M\in\Match^N(2,2)\subset\Mat^N(2,2)$ is sparse.
The non-zero blocks are $1\times 1$ and $2 \times 2$,
naturally in correspondence with the vertices $i$ and edges $ij$ respectively
of the complete graph $K_N$.
\ignore{
  The image
$F(\sigma)$ 
  is a sparse matrix whose rows (and columns) may be indexed by
ordered pairs $(i,j)$ with $i,j \in \ul{N}
 = \{ 1,2,...,N \}$.
The non-zero blocks are $1\times 1$ and
$2 \times 2$, naturally in correspondence with the vertices and edges
respectively
of the complete graph $K_N$. }

}}


\frx{Complete graph visualisation}{{
    Solutions can thus be visualised using $K_N$
    \[
    \includegraphics[width=7.5\ccm]{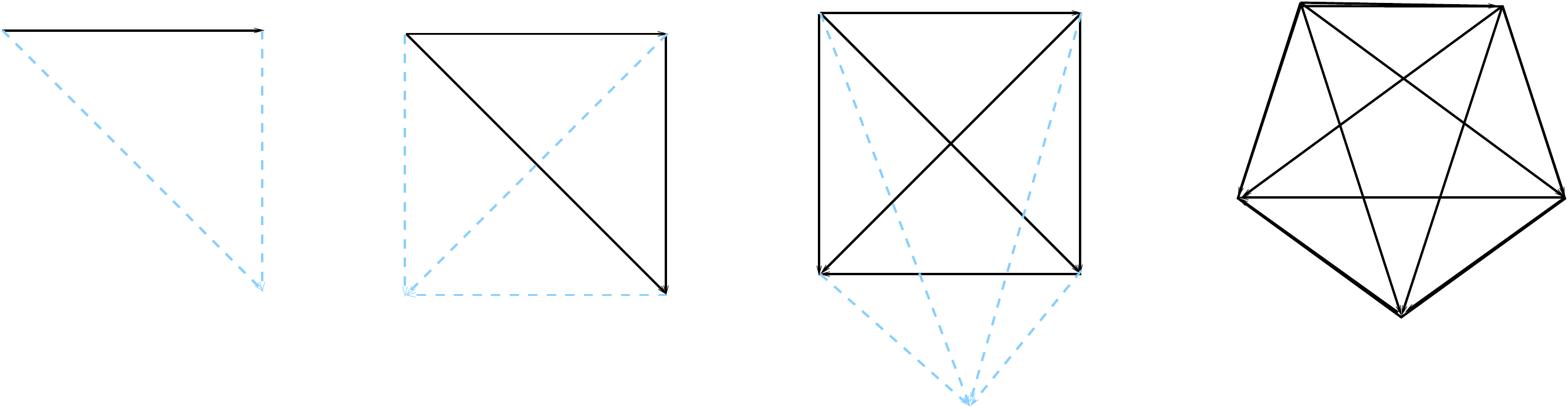}
\hspace{.51\ccm} ...\hspace{.51\ccm}
    \includegraphics[width=3.5\ccm]{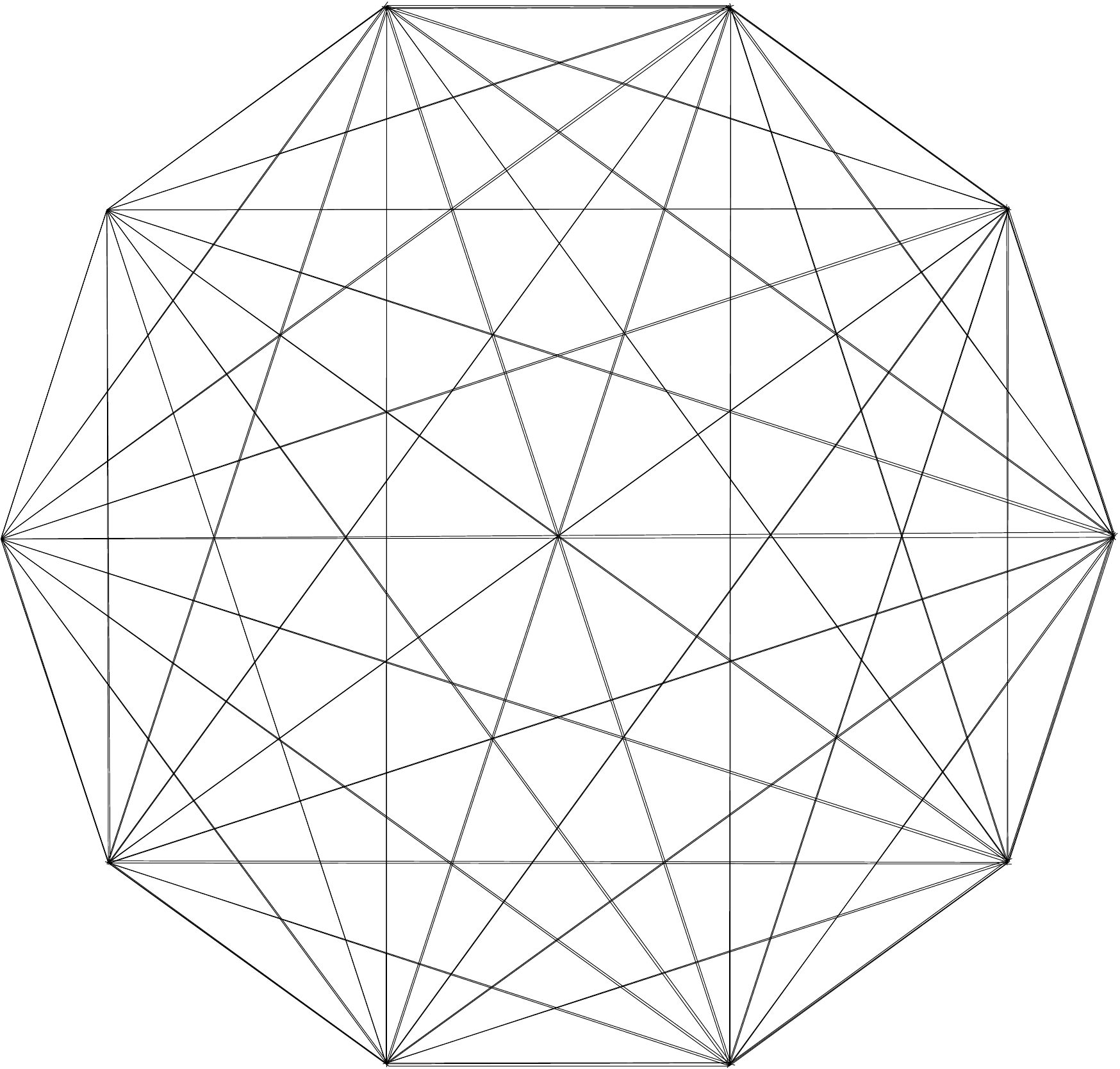}
    \]
%
    
}}

\frx{}{{
    ...So that

    \includegraphics[width=3.5\ccm]{xfig/reg10gon7.eps}
    becomes \hspace{-1\ccm}
    \raisebox{-1\ccm}{
      \includegraphics[width=8.5\ccm]{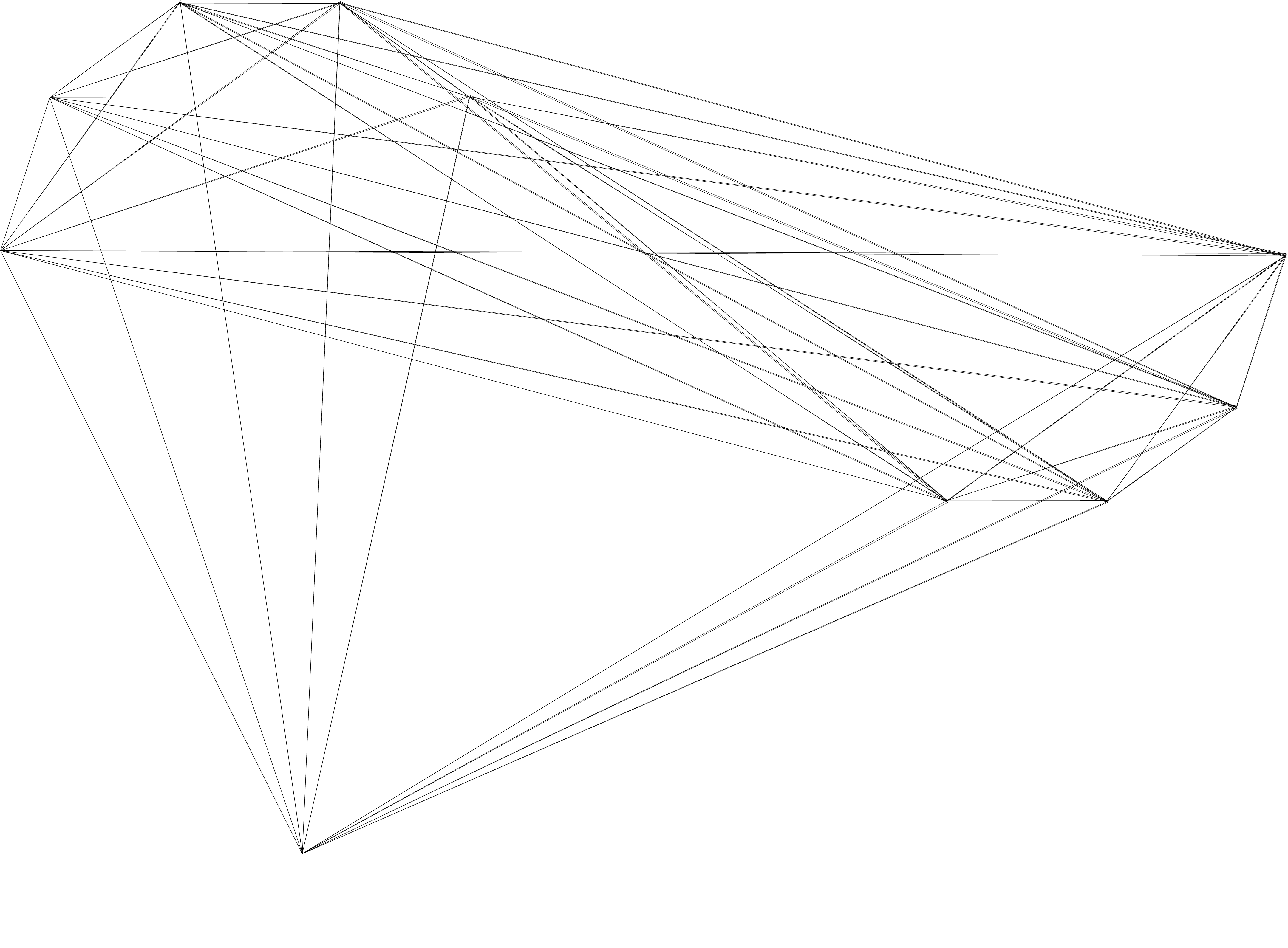}

} 
    
}}


\begin{figure}
\[
\includegraphics[width=6cm]{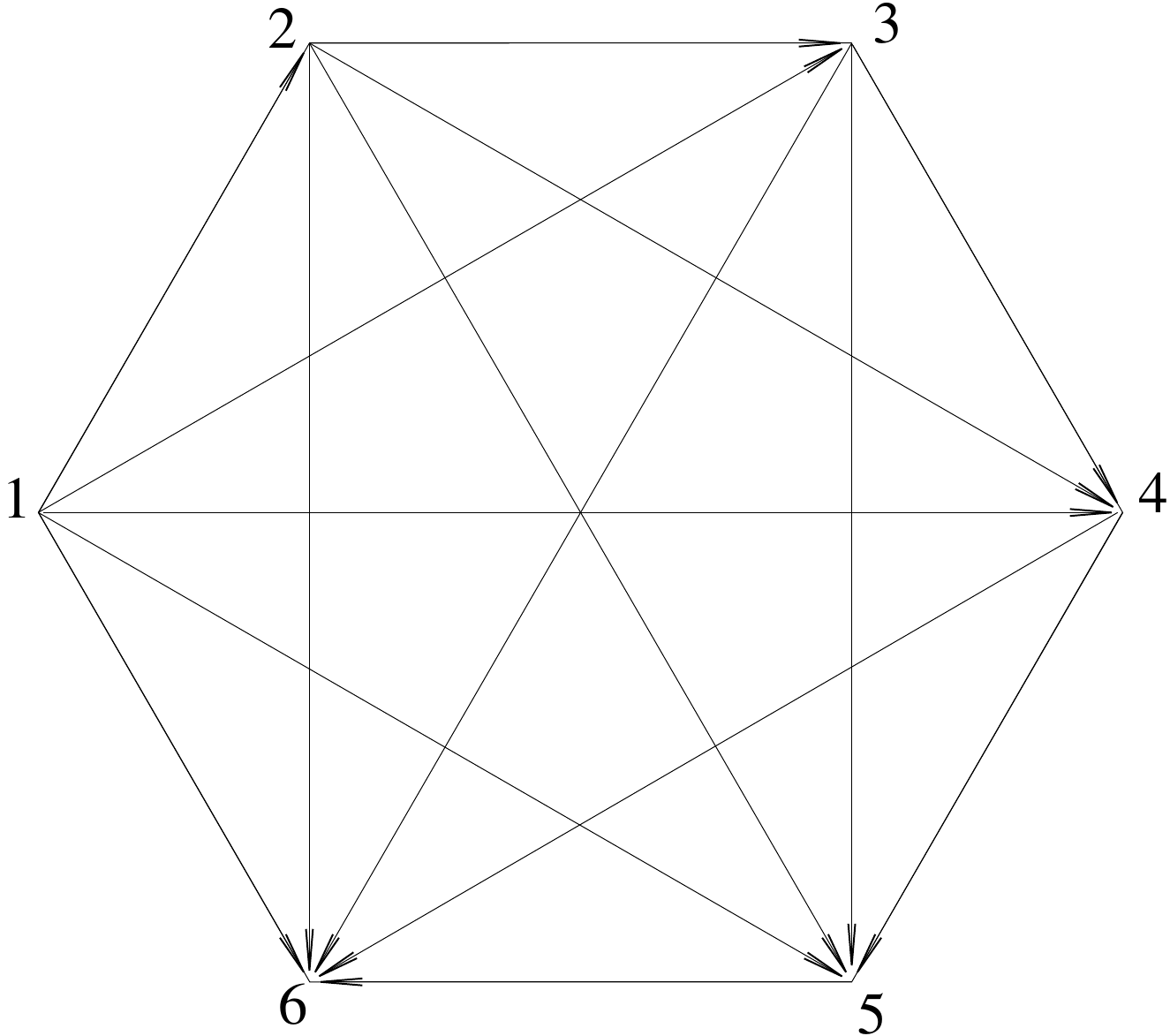}
\hspace{2.1cm} 
\includegraphics[width=6\ccm]{tex-eric/XCvXCcrop.pdf} 
\]
\caption{(a) The 
graph $K_6$; and
(b) conventions for a Higgs configuration on $K_4$.  \label{fig:K6}}
\end{figure}

\fr{}{{

\mdef \label{de:KN}
    For $N \in \N$, $K_N$ denotes the directed graph with vertex set
$\ul{N} = \{1,2,...,N \}$ and an edge $(i,j)$ whenever $i<j$.
If $S$ is an ordered set then we write $K_S$ for the
corresponding complete graph on vertex set $S$.
\hspace{.1in}
For example, with $N=6$:
see Fig.\ref{fig:K6}(a).
}}


\newcommand{\alphain}{{\underline{\alpha}}} 
\newcommand{\alphaout}{{\overline{\alpha}}} 

\fr{Gauge configuration}{{

\mdef \label{pa:config}
    A 
    configuration
$\alphaa$
on $K_N$ is an assignment of a variable to
each vertex and a $2\times 2$ matrix of variables to each directed
edge.
By (\ref{pa:sparse}) this encodes an element of $\Match^N(2,2)$.
%
See Fig.\ref{fig:K6}(b),
where we name matrix entries as in (\ref{eq:abcd12}), with $a_i =a_{ii}$.


We should think of these assignments in this geometrical form.
But sometimes it is convenient to give them in-line.
Then the order we shall take is
with vertex assignments first in the natural order, and then edge assignments
in the order 12,13,23,14,24,34 and so on.
That is, for $M \in \Match^N(2,2)$ we may write it,
{first generally, and then using (\ref{eq:abcd12})
with $a_i = a_{ii}$}:
\beq \label{eq:alconv1}
\alphain(M) = (M_{11,11}, M_{22,22}, \ldots, M_{NN,NN},
\mat{cc} M_{12,12} & M_{12,21} \\ M_{21,12} & M_{21,21} \tam,
\mat{cc} M_{13,13} & M_{13,31} \\ M_{31,13} & M_{31,31} \tam,
\ldots)
\eq
\[
=
(a_1, a_2, \ldots, a_N,
\mat{cc}\! a_{12} & b_{12} \!\\ \! c_{12} & d_{12} \!\tam\!,
\mat{cc}\! a_{13} & b_{13} \\ \! c_{13} & d_{13} \tam,
\ldots)
\hspace{1.6in}
\]
\[ \hspace{.562in}
=
(a_1, a_2, \ldots, a_N, A(1,2), A(1,3),\ldots, A(N\!-\!1,N)),
\hspace{.2in}
\mbox{ where }
A(i,j)  := \mat{cc}\!\!\! a_{ij}&b_{ij} \!\!\\ \!\!\! c_{ij}&d_{ij} \!\!\tam
\] 


}}





\begin{figure}
  \newcommand{\ddir}{tex-p/}
\newcommand{\ecm}{-18.76cm}
\newcommand{\eecm}{-18.99cm}
\newcommand{\furk}{
\includegraphics{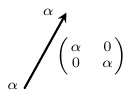} \hspace{\ecm} }
\newcommand{\ferk}{
  \includegraphics[trim=5cm 22cm 14cm 3.5cm]{\ddir/labelscheme2x3a.pdf} }
\newcommand{\ferko}{
  \includegraphics[trim=5cm 22cm 14cm 3.5cm]{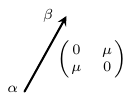} }
\newcommand{\ferkoo}{
  \includegraphics[trim=5cm 22cm 14cm 3.5cm]{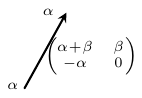} }
\newcommand{\ferkooo}{
  \includegraphics[trim=5cm 22cm 14cm 3.5cm]{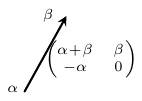} }
\newcommand{\ferkoooo}{
  \includegraphics[trim=5cm 22cm 14cm 3.5cm]{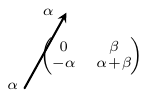} }
\newcommand{\ferkooooo}{
  \includegraphics[trim=5cm 22cm 14cm 3.5cm]{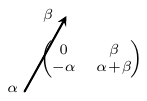} }
\newcommand{\feerk}{
  \includegraphics[trim=5cm 22cm 14.23cm 3.5cm]{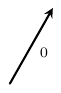} }
\newcommand{\feerko}{
  \includegraphics[trim=5cm 22cm 14.23cm 3.5cm]{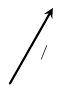} }
\newcommand{\feerkoo}{
  \includegraphics[trim=5cm 22cm 14.23cm 3.5cm]{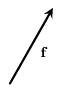} }
\newcommand{\feerkooo}{
  \includegraphics[trim=5cm 22cm 14.23cm 3.5cm]{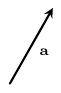} }
\newcommand{\feerkoooo}{
  \includegraphics[trim=5cm 22cm 14.23cm 3.5cm]{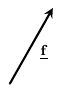} }
\newcommand{\feerkooooo}{
  \includegraphics[trim=5cm 22cm 14.23cm 3.5cm]{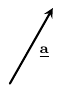} }
\centering
\[
{ }\hspace{-.6cm}
\ferk
\ferko
\ferkoo \ferkooo \ferkoooo \ferkooooo
\hspace{-15.5cm}
\raisebox{-1.75cm}{
\feerk \feerko
\feerkoo \feerkooo\feerkoooo \feerkooooo
}
\hspace{-15cm}
\raisebox{-2.19cm}{(1) \hspace{1.65cm} (2)\hspace{1.875cm} (3a)\hspace{2.15cm} (3b)
  \hspace{2.15cm}(4a)\hspace{1.85cm}(4b)}
\]
\caption{Schematic for possible edge submatrices up to X-equivalence. \label{fig:N2}}
\end{figure}

\module{reps-N2}


\newcommand{\AAA}{{A}} 
\newcommand{\aaf}{a_2}
\newcommand{\aag}{a_1}  


\begin{proposition} \label{pr:L2class2}
  For $N=2$ the following gives a complete classification
  of charge conserving functors $\FF$ from $\Bcat$ up to
  X-equivalence.
  We use the coefficient names
  given by
$\FF(\sigma) =
 \mat{cccc} \!\!\aag \!\\ &a_{12}&b_{12} \\ &c_{12}&d_{12} \\ &&&\!\aaf \!\tam$.
That is, $\alphain(\FF(\sigma))=(\aag,\aaf,\AAA(1,2))$.
We have:
\begin{enumerate}
\item
The $0$ case: $a_1=a_{12}=d_{12}=a_2$, and $b_{12}=c_{12}=0$,
i.e., $\FF(\sigma)=\alpha\cdot Id$, with $\alpha\neq 0$.
\item
The $/$ case: $a_{12}=d_{12}=0$, and $a_1=\alpha$, $a_2=\beta$ and $b_{12}=c_{12}=\mu_{}$,
say, with $\alpha,\beta,\mu\neq 0$.
\item
The $+$ cases: $d_{12}=0$, $a_{12}=\alpha+\beta$, $b_{12}=1$, $c_{12}=-\alpha\beta$ where $\alpha+\beta\neq 0$ and either 
\begin{enumerate}
\item $\pai$ case: $a_1=\alpha\neq a_2=\beta$ or
\item $\pfii$ case: $a_1=a_2=\alpha$
\end{enumerate}
\item The $-$ cases: $a_{12}=0$, $d_{12}=\alpha+\beta$, $b_{12}=1$, $c_{12}=-\alpha\beta$ where $\alpha+\beta\neq 0$ and either 
\begin{enumerate}
\item $\underline{\pai}$ case: $a_1=\alpha\neq a_2=\beta$ or
\item $\underline{\pfii}$ case: $a_1=a_2=\alpha$.
\end{enumerate}

\end{enumerate}

\end{proposition}
%
\noindent
{NB, the $/$ and $\pfii,\pai$ varieties `touch', but the given partition is
good for treating higher ranks, as we shall see.}
%
%
%
\label{pa:schema2}%
    Possible edge submatrices are
    as shown in Figure~\ref{fig:N2}.
    (The labels $\pfii$ and $\pai$ are
    taken from 
    \emph{ferromagnetic} and \emph{anti-ferromagnetic},
indicating
if vertex labels are the same or not.) 

\ignore{{
\begin{figure}
\[\includegraphics[width=14.5\ccm]{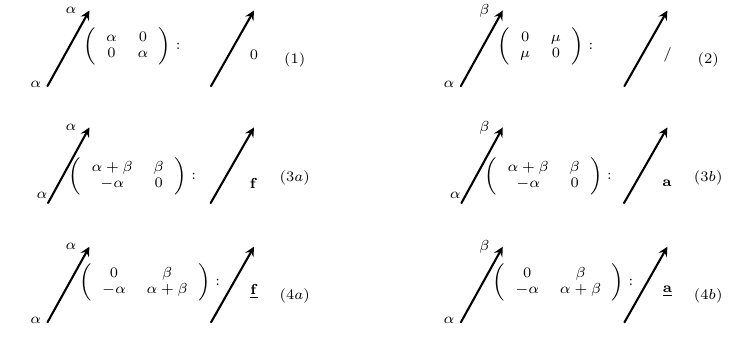}
\]
\caption{\label{fig:N2}
Schematic for possible edge submatrices
{up to X-equivalence}.}
\end{figure}
}}



\medskip
\module{reps-N2proof2}


\beff{}{{
\mdef\label{proof of N2}
\noindent
{\bf Proof.}   
%
Equations (\ref{eq:31xx0}-(\ref{eq:31xxx0})
(and their images under the action of $\Sigma_2$) are relevant for $N=2$.
\ppmx{[surely we don't need to write them again here??-]}
\ignore{{
\begin{eqnarray*}
&&a_{12} (a_{1}^2 -a_{12} a_{1} -b_{12} c_{12}) =0,\hspace{1cm}
a_{12} (a_{2}^2 -a_{12} a_{2} -b_{12} c_{12}) =0,\\
&&a_{12} c_{12} d_{12} = 0 = a_{12} b_{12} d_{12} ,\hspace{1cm}
a_{12} d_{12} (a_{12}-d_{12})=0.
\end{eqnarray*}
}}
Under the transposition $(1\/2) \in \Sym_2$ we obtain 
similar
equations 
but with $a_{12}$ and $d_{12}$ interchanged.
Observe that $a_1,a_2$ and $a_{12}d_{12}-b_{12}c_{12}$ are non-zero by invertibility.
We take it
case by case 
to find necessary and sufficient conditions for these constraints.

If $b_{12}c_{12}=0$ then $a_{12}d_{12}\neq 0$, so the equations
imply $a_{12}=d_{12}=a_1=a_2$, which is case $0$.
It is routine to check that this  indeed gives a solution for any
nonzero value.

If $b_{12}c_{12}\neq 0$ then there are 3 cases to consider:
$a_{12}=d_{12}=0$ or exactly one of $a_{12}$ or $d_{12}$ is non-zero.
In the first case all equations are satisfied, so that $a_1,a_2,b_{12},c_{12}$ may be chosen arbitrarily.
The X transformation then allows $b_{12}=c_{12}$,
{so this gives the $/$ case}.

Next
suppose that $b_{12}c_{12}a_{12}\neq 0$, and $d_{12}=0$.  This implies that $b_{12}c_{12}=a_2(a_{12}-a_2)=a_1(a_{12}-a_1)$.
Now the characteristic polynomial of the $2\times 2$ matrix
(call it $A(1,2)$)
is $x^2-a_{12}x-b_{12}c_{12}$ so that $a_1$ and $a_2$ are eigenvalues of $A(1,2)$.  Let $\alpha,\beta$ be these eigenvalues.
Since $Tr(A(1,2))=a_{12}$ we have that $a_{12}=\alpha+\beta$, and
{similarly} 
$Det(A(1,2))=-b_{12}c_{12}=\alpha\beta$.
Using the X transformation we may assume that
$b_{12}=-\alpha\beta$ and $c_{12}=1$ or vice versa.
Notice that the only restriction on $\alpha,\beta$ is that $\alpha+\beta\neq 0$.
Now there are two cases, up to the obvious labeling ambiguity of $\alpha$ and $\beta$: either $a_1=a_2=\alpha$
(case $\pfii$)  
or $a_1=\alpha$ and $a_2=\beta$
(case $\pai$).  

The case $a_{12}=0$ is completely analogous:
we find that $d_{12}=\alpha+\beta$, and 
$b_{12}=-\alpha\beta$, $c_{12}=1$
{(or vice versa)}
and $a_1=a_2=\alpha$
(case $\mfii$)   
or $a_1=\alpha$ and $a_2=\beta$
(case $\mai$).   

In
any of the last 4 cases it is possible that $\alpha=\beta$,
whereupon the 
$\pfii$
case prevails.
%
%
\qed

}}



\ignore{{
Note that in either
`two-eigenvalue' 
case $b,c$ are only constrained through $bc$.
Thus we may take $b=c$ or for example $c=1$; and we may continuously
vary between these. This is a variant of the classical
semi-normal/normal isomorphism and for $\Bcat$ itself these
representations are generally classed together.
However we are interested in extension to $\LBcat$,
and it is not true in general that conjugations in
images of $\Bcat$ commute
with  extension. Thus we should treat them separately. 
}}





\medskip
\module{braid-theorem0}

\def\barroman#1{\sbox0{#1}\dimen0=\dimexpr\wd0+.1pt\relax
  \makebox[\dimen0]{\rlap{\vrule width\dimen0 height 0.06ex depth 0.06ex}%
    \rlap{\vrule width\dimen0 height\dimexpr\ht0+0.03ex\relax 
            depth\dimexpr-\ht0+0.09ex\relax}%
    \kern.5pt#1\kern.5pt}}

\section{Constructions for the main Theorem}\label{ss:main}




  We write $\Lambda_N$ for the set of integer partitions of $N$.
  We may write $\lambda \vdash N$ for $\lambda\in \Lambda_N$. 


\beff{On classification theorems}{{
    
\subsection{ 
      How to read  the classification  
      Theorem~\ref{th:mainx}:
      } $\;$ 
      
An example of a classification theorem
is of course Young's classification of
irreducible representations of the symmetric group over $\C$.
Here one says that {\em irreducible} representations
of $\Sym_n$ may be classified
up to {\em isomorphism} by the set
$\Lambda_n$.
Each of the three aspects of this formulation flags  challenges that we
must address in our problem.

Firstly, 
  analogously
to the integer partitions
we will need to introduce  
some 
  combinatorial structures.
  (Our Theorem also gives a {\em construction}.
  We will introduce notation for this too.)

Secondly,
group/algebra representations form an additive category
in a natural way.
The category property yields a `good' notion of isomorphism.
The additive property
simultaneously implies that there are unboundedly many
classes of representations,
but also that these can be effectively
classified by classifying their additive components.
%
There is no canonical lift of these properties to the monoidal category setting
(confer for example \cite{MazorchukMiemietz16,MazorchukMiemietz14,RumyninWendland,KapVoev}),
so our lift  
sheds new light on this aspect of higher representation theory.

And finally, 
  treated individually symmetric group algebras over algebraically closed
  fields are Artinian, so the Jordan--Holder property tells us that every
  representation has a decomposition series with irreducible factors.
  Indeed in the complex case the decomposition multiplicities characterise a
  representation up to isomorphism.
  There is no canonical lift of this property to monoidal categories
either
(roughly since decomposition series are additive decompositions,
  while the monoidal structure is multiplicative).
  This is one of the places where the rigidity of the $\Match^N$
  target categories
  saves the day --- we will be able to classify {\em all} representations directly.


}}


\medskip


\befff{Classification: Notation}{{

    Next
    we construct two sets for each $N\in \N$.
    One is the set $\SSS_N$ to which we may apply an
    algorithm (given in (\ref{pa:rec}))
    to construct all varieties of solutions
    at level-$N$.
    The other gives a transversal of this set under the $\Sigma_N$ action
    from (\ref{restriction-and-symmetries})
    (and frames the effect of the $\Zz_2$ action),
    thus addressing the classification up to isomorphism.
}}

\subsection{Braid representations
            from multisets of row-2-coloured Young diagrams}\label{ss:comb0} $\;$ 
\newcommand{\jjj}{{\mathsf j}}  

Just as
multisets of integers -
represented as Young diagrams 
-
index $\Sym_N$ representations,
we will see that
for braid representations we can use
multisets of row-2-coloured 
composition
diagrams.
In fact these are a natural combinatorial generalisation of $\Lambda_N$.
First we explain this combinatoric;
then we show how each diagram leads to
a class of charge-conserving braid representations;
then we prove that this construction
exhausts the set of all such braid representations. 


\newcommand{\staro}{*}  

\mdef \label{41}
Given a set $\Omega$ and a map $\psi\in \Setcat(\Omega,\N)$ define
the
sets of multisets
of total degree $N$:
\[
J_N(\Omega) = J_N(\Omega,\psi) =
\{ f \in \Setcat(\Omega,\N_0) \; | \; \sum_{w\in\Omega} f(w) \psi(w) =N \}
\]
Thus for example
if $\Omega=\N$ and $\psi(n)=n$ then $J_N(\N) = \Lambda_N$.

\mdef
For $S$ a set, define $\psi: S^\staro \rightarrow \N$
by $\psi(w) = length(w) +1$.
(Cf. (\ref{de:Set}).)
Observe 
that
up to isomorphism
$J_N(S^\staro)$ depends only on the cardinality of $S$.
Consider $J_N(S^\staro)$ for $S=\{1\}$, $\{1,2\},\{1,2,3\}$
and so on.

Of course $\psi : \{1\}^\staro \stackrel{\sim}{\rightarrow} \N$.
So
$J_N(\{1\}^\staro ) $  is a realisation of $\Lambda_N$.
For example
$f \in \Setcat(\{1\}^\star,\N_0)$
given by
$f(1)=1$; $f(11)=2$; $f(w)=0$ otherwise,
becomes
$(3,3,2)\vdash 8$.

In this $J_N(S^\staro)$ notation,
$\{1,2\}^N  $
encodes the set $\Gamma_N$ of 
{\em compositions} of 
$N$ (via the $\jjj$-function given just below);
and $\{1,2\}^\staro  $
encodes the set $\Gamma$ of all compositions;
so $J_N(\{1,2\}^\staro ) $
encodes the set  of multisets of
compositions of total degree $N$.

The set we need ---
of multisets of  
2-coloured compositions
---
is $J_N(\{1,2,3\}^\staro ) $. 
The construction works as follows.  
Firstly $\Gamma^2_N$ denotes the set of 2-coloured compositions of $N$,
whose elements are compositions of $N$
together with a two-part partition
of the components.
We may draw $\lambda\in\Gamma^2_N$
as a stack of rows of boxes, where the top row is unshaded
but other rows can be unshaded or shaded
(in the first or second part respectively).
For example:
\raisebox{-.1in}{\includegraphics[width=.783cm]{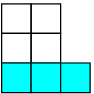}}
$\in\Gamma^2_7$,
and see Fig.\ref{fig:N1234}.
{(This box formulation  
is useful, but
note
{\em columns} have 
no significance,
cf. ordinary Young diagrams%
.)}
Next we show that $\{1,2,3\}^*$ encodes the set
$\Gamma^2$ of all 2-coloured compositions.



\mdef \label{de:j}
Given a 2-coloured composition
or equivalently a row-2-coloured Young/composition diagram
$\lambda\in\Gamma_N^2$
then diagram
$\jjj_1(\lambda)\in\Gamma_{N+1}^2$ is obtained by adding a box to the last row;
$\jjj_2(\lambda)\in\Gamma_{N+1}^2$ is obtained by adding a box in a new unshaded row;
and
$\jjj_3(\lambda)\in\Gamma_{N+1}^2$ is obtained by adding a box in a new shaded row.
Starting with the one-box (unshaded) diagram,
a word  
$w \in \{1,2,3\}^*   $
then defines a row-2-coloured diagram by applying
the $\jjj_i$ maps in the sequence given by $w$.
That is
$\jjj : \{1,2,3\}^n \rightarrow \Gamma_{n+1}^2$
is given by
$w_1 w_2 ... w_n \mapsto
 \jjj_{w_n} \circ ... \circ \jjj_{w_2} \circ \jjj_{w_1} (\square)$.
Thus $\jjj(121311) $ is given by
$\square
\stackrel{\jjj_1}{\mapsto}
\includegraphics[width=.53cm]{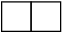}
\stackrel{\jjj_2}{\mapsto} 
\includegraphics[width=.53cm]{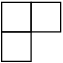}
\stackrel{\jjj_1}{\mapsto} 
\includegraphics[width=.53cm]{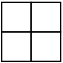}
\stackrel{\jjj_3}{\mapsto} 
\includegraphics[width=.53cm]{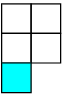}
\stackrel{\jjj_1}{\mapsto} 
\includegraphics[width=.53cm]{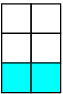}
\stackrel{\jjj_1}{\mapsto} 
\includegraphics[width=.783cm]{xfig/y223.eps}
$
and so on.   
Note from the construction
that $\jjj$ is
reversible and hence
a bijection in each rank.


\ignore{{
\mdef
\ppm{[time for bed, but in the morning the idea is that
part of this para will come
after the following one (!) just to give a definite order...,
and the rest deleted.]}
It will be convenient to have a total order on $\Gamma^2$.
We can pull one through from a total order on $\{1,2,3\}^*$
using the isomorphism $\jjj$.
For $\{1,2,3\}^*$
let us say firstly that longer sequences come before shorter ones:
$111<11<1$ and so on;
and for equal-length sequences we use dictionary order,
that is $111<112<113<121<...<333$.

Armed with this order, we have a way of
writing out a sequence of diagrams corresponding to a multiset
$f \in J_N(\{1,2,3\}^*)$. We use the above order, and the multiplities from $f$.
We can now assign vertices of $K_N$ to the $N$ boxes
in a multiset of 2-coloured compositions by filling the boxes with the
vertex numbers.
We fill the first diagram starting with  
1 in the top left corner then
increasing by 1 along rows and then down columns.
Then continue with the next diagram; and so on.

Thus for example:
\[
\includegraphics[width=3.494cm]{xfig/tree012ex10bbb1XI.eps}
\]

}}


\mdef \label{44}
Note that 
$f \in J_N(\{1,2,3\}^\staro )$ assigns a multiplicity $f(w)$ to each
word $w$ and hence
diagram $\jjj(w)$,
such that the total degree is $N$.
We will show how
each such $f$ corresponds to an equivalence class of varieties of
configurations on $K_N$, and hence
matrices
that 
give the functors $\FF:\Bcat\rightarrow\Match^N$.

Firstly we
arrange 
2-coloured composition diagrams in multiplicities as given by $f$,
so that we have $N$ boxes altogether.
For example
with $f$ given by $f(1321)=1$, $f(131)=1$, $f(\emptyset)=2$ and
$f(w)=0$ otherwise we obtain an 11-box composite diagram such as:
\beq \label{eq:54110}
\jjj(f) \; =\;\;
\raisebox{-.21in}{
\includegraphics[width=2.94cm]{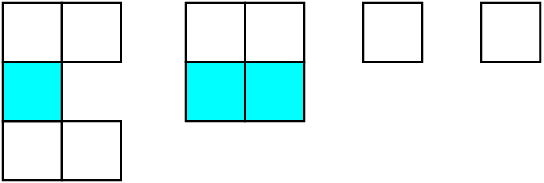}
}
\eq
Secondly we
assign 
the $N$ vertices of $K_N$ to these $N$ boxes.
For example:
\beq \label{eq:5411}
\jjj(f) \; \mapsto\;\;
\raisebox{-.21in}{
\includegraphics[width=3.494cm]{xfig/tree012ex10bbb1XI.eps}
}
\eq

Write
$\SSS_N$
for the set of {\em all} assignments
{\em up to the order within each row}.
Thus
\beq \label{eq:SSS2}
\SSS_2 = \{\;
\includegraphics[width=.8cm]{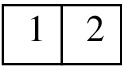}\;,\;\;
\includegraphics[width=.41cm]{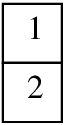}\;,\;\;
\includegraphics[width=.41cm]{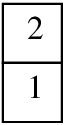}\;,\;\;
\includegraphics[width=.41cm]{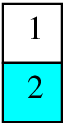}\;,\;\;
\includegraphics[width=.41cm]{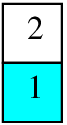}\;,\;\;
\includegraphics[width=1.03cm]{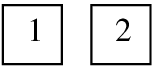}\;
\}.
\eq
We can choose any such assignment for each $f$.
But 
it will be convenient to have a prefered assignment,
of vertices in `book' order as in (\ref{eq:5411}).
This requires the diagrams in a composite to be written in a specific
sequence. We can do this using
a declared  total order on
$\Gamma^2$ (or equivalently on $\{1,2,3\}^\staro $)
 --- noting that identical diagrams are indistinguishable.

\newcommand{\bbbbb}{{\mathsf b}}  

For
a total order on
$\{1,2,3\}^*$
let us say firstly that longer sequences come before shorter ones:
$111<11<1$ and so on;
and for equal-length sequences we use dictionary order,
that is $111<112<113<121<...<333$.
Then
we get our example
$\jjj(f)$ 
in (\ref{eq:54110}) above,
and hence (\ref{eq:5411}), on the nose.
Write
$\bbbbb : \jjj(J_N(\{1,2,3\}^*)) \rightarrow \SSS_N $
for the book order assignment.
Define $\TTT_N = \bbbbb( \jjj(J_N(\{1,2,3\}^*)) )$.
For example
$\TTT_2 =
\{\;
\includegraphics[width=.8cm]{xfig/tree012ex2aab.eps}\;,\;\;
\includegraphics[width=.41cm]{xfig/tree012ex2ab1b.eps}\;,\;\;
\includegraphics[width=.41cm]{xfig/tree012ex2abb.eps}\;,\;\;
\includegraphics[width=1.03cm]{xfig/tree012ex2acb.eps}\;
\}
$.
\ignore{{
for example, with $f$ given by $f(1321)=1$, $f(131)=1$, $f(\emptyset)=2$ and
$f(w)=0$ otherwise we obtain an 11-box diagram:
\[
\includegraphics[width=2.94cm]{xfig/tree012ex10bbb1XI0.eps}
\]
}}
\ignore{{
We can now assign vertices of $K_N$ to the $N$ boxes
`in' $f$ 
by filling the boxes with the
vertex numbers.
We fill the first diagram starting with  
1 in the top left corner then
increasing by 1 along rows and then down columns.
Then continue with the next diagram; and so on:
\[
\includegraphics[width=3.494cm]{xfig/tree012ex10bbb1XI.eps}
\]
}}


\medskip
{There is a correspondence between $\SSS_2$ and our classification
for $N=2$ in Prop.\ref{pr:L2class2}.}
The rest of our exposition 
  amounts to showing $\SSS_N$ yields an index scheme for all
varieties of 
representations.
Noting
Lemma~{\ref{lem:chacon1}(III)}, the
$N=2$ solutions
also give the possible edge configurations 
for all $N$.


\module{notation-partition}

\ignore{{

\ppm{[...thus these next bits would now be demoted:]}
To obtain an element of the index set $\SSS_N$
that we shall use
one simply
draws $f(\lambda)$ copies of the diagram for each diagram $\lambda$
(so there are $N$ boxes in total),
then 
fills the boxes with the
numbers $1,2,...N$, such that each row is increasing
(noting that identical diagrams are indistinguishable,
so unordered amongst themselves).
For example with $f$ given by $f(1321)=1$, $f(131)=1$, $f(\emptyset)=1$ and
$f(w)=0$ otherwise we obtain a 10-box diagram which we may fill:
\[
\includegraphics[width=2.94cm]{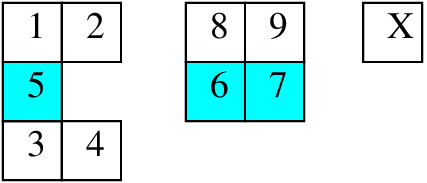}
\]
Altogether for $N=2$ we have $\;$ 
 $
\SSS_2 = \{\;
\includegraphics[width=.8cm]{xfig/tree012ex2aab.eps}\;,\;\;
\includegraphics[width=.41cm]{xfig/tree012ex2ab1b.eps}\;,\;\;
\includegraphics[width=.41cm]{xfig/tree012ex2ba1b.eps}\;,\;\;
\includegraphics[width=.41cm]{xfig/tree012ex2abb.eps}\;,\;\;
\includegraphics[width=.41cm]{xfig/tree012ex2bab.eps}\;,\;\;
\includegraphics[width=1.03cm]{xfig/tree012ex2acb.eps}\;
\}.
$
\\
Next we show how each element
of each $\SSS_N$ 
gives a class of braid representations.
That is, each element gives a configuration on the complete graph $K_N$
as in (\ref{pa:config}) that gives an $\FF(\sigma)$.
\\
\medskip
\ppm{
[maybe we
can delete or demote most of the rest of this section.]}


\mdef
We write $\Gamma_N$ for the set of compositions of $N$.
  For a composition (or partitition)
  $\lambda=(\lambda_1, \lambda_2, ..., \lambda_L)$
  one may 
  simply
  write
  $\lambda=\lambda_1 \lambda_2 \cdots\lambda_L$
where no ambiguity arises.
Otherwise write 
{$\lambda = \lambda_1 \ovar\lambda_2 \ovar\cdots \ova \lambda_L$}
(or
$\lambda=\lambda_1 +\lambda_2 +\cdots+ \lambda_L$).
\\
Note that there is a lexicographic order on $\Gamma_N$.
For example
$\Gamma_2 = \{ 11, 2 \}$ and
$\Gamma_3 = \{ 111, 12, 21, 3 \}$ are written in the lex order.

}}


\beff{Set partitions}{{

\mdef \label{pa:child1st}
  Given a set $S$ we will write $\PP(S)$ for the set of partitions;
  and $\PP_2(S)$ for the set of partitions into at most 2 parts.
  If $N \in \N$ then $\PP(N) = \PP(\ul{N})$.

``Children-first order'':
  Given a partition $p$ of $\{ 1,2,...,N \}$
  (or any subset thereof),
we will order the parts according to their lowest numbered element;
and write $p_i$ for the $i$th such part.

\mdef
A refinement of set partition $p$ is a further partition of the parts of
$p$
(`nations', say)
into possibly smaller parts (`counties', say).
We will write $p<q$ (or $q>p$) if $q$ is a refinement of $p$.
\\
We may write $q|p_i$ for the restriction of $q$ to the partition of $p_i$.

}}


\ignore{{  

\newcommand{\PPu}{\PP^u}
\newcommand{\PPm}{\PP^{-}}

\beff{Partition combinatorics}{{ 

\mdef
Note that there  are natural bijections between the following sets:
\\
$\{ (p,q) \; | \; p,q \in \PP(S), p<q \}$,  $\;\;$
$\{ (q,p) \; | \; p,q \in \PP(S), p<q \}$,  $\;\;$
\hspace{.01in}
$ \{ (q,p') \; | \; q \in \PP(S), \; p' \in \PP(q) \}$
$\;$
and
$$
\PPu(S) =
 \bigcup_{q \in \PP(S)} \PP(q) = 
 \{ p' \; |  \; p' \in \PP(q), \; \mbox{for some } q \in \PP(S) \}
$$
 \\
 The map from first to second is simply to flip the pair.
 The map from second to third uses that $p<q$ has elements that are unions of
 parts from $q$, so we obtain a $p'$ by partitioning $q$ according to this data.
 Finally
 given $p' \in \PP(q)$ we can recover $q$ by taking the union over parts,
 so the $q$ is redundant.
 That is we can map $(q,p') \mapsto p'$.
 
 Note that $\PP$ has a kind of inverse $\PPm$ that takes an element of
any $\PP(S)$ and forms the union, thus restoring $S$
(caveat: this construction does not play well with the convention that
the result of a function applied to a subset of the domain is the set of images,
since here it may not be clear if an argument is a subset of the domain or an
element of it).
\\
Example:
$$
\PPu(\{1,2\})=
\PP(\{ \{1,2\}\}) \cup \PP(\{ \{1\},\{2\}\})
  = \{ \{\{\{1,2\}\}\}, \{\{ \{1\},\{2\}\}\}, \{ \{\{1\}\},\{\{2\}\}\} \}
$$

\ignore{{
Sequences of form $(p,q,r,...)$ with $p<q<r<...$
(any length of sequence, $l$)
are sometimes called ramified partitions.
}}

}}


\beff{Permutations}{{

\mdef
  Given a  
  finite set $S$ then 
  $\Perm(S)$ 
  is the set of total orders of the  
  elements, i.e.
  $\Perm(S) = Hom^{iso}(S,\ul{|S|})$,
  the set of isomorphisms
  (although we express elements simply by giving the order).
\\
If $p<q$ as above
then $\Perm_p(q)$
denotes the set of total orderings of  the elements of $q$ in each $p_i$
(counties in each nation)
without imposing an overall total order:
$$
\Perm_p(q) = 
            \prod_{p_i \in p}  \Perm(q|p_i)
$$ 
Example: \hspace{.1in}
$
\Perm_{ \{ \{1\},\{2,3\}\}}( \{ \{1\},\{2\},\{3\}\}) =
\{  \{ (\{1\}),(\{2\},\{3\}) \} , \{ (\{1\}),(\{3\},\{2\}) \} \}
$ 
\\
Finally,
$\PP_2(q/p) = \prod_{p_i \in p} \PP_2(q|p_i) \;$
--- the set of partitionings of counties of each nation into at most two parts.
Example:
$ \PP_2(\{\{1\},\{2\},\{3\},\{4\}\}/\{\{1,2\},\{3,4\}\})
  = \PP_2(\{\{1\},\{2\}\}) \prod \PP_2(\{\{3\},\{4\}\})$.  
}}

\subsection{Constructions for solutions} $\;$ 

\sout{
Beside the geometrical 
encoding
we have a combinatorial encoding
of $(p,q,\rho,s)$
in terms of collections of Young tableaux
--- one tableau for each nation; one row for each county,
rows stacked  
according to $\rho$,
with each nation coloured so that counties  
that are in the same part in $s$ have the same colour.
The {\em same} element (but now  
complete with $s$ data)
becomes:}

Each of our diagrams yields a decoration of the vertices and edges of $K_N$,
and hence a candidate for a functor as follows.
Firstly we visualise the clustering of vertices into counties and nations
as in Fig.\ref{fig:puppyears}.
\begin{figure}
\[
\raisebox{.251in}{
\includegraphics[width=3.84cm]{xfig/tree012ex10bbb.eps}
}
\hspace{.2in}
\raisebox{.53in}{
$\leadsto$}
\hspace{-.431in}
\includegraphics[width=10.84cm]{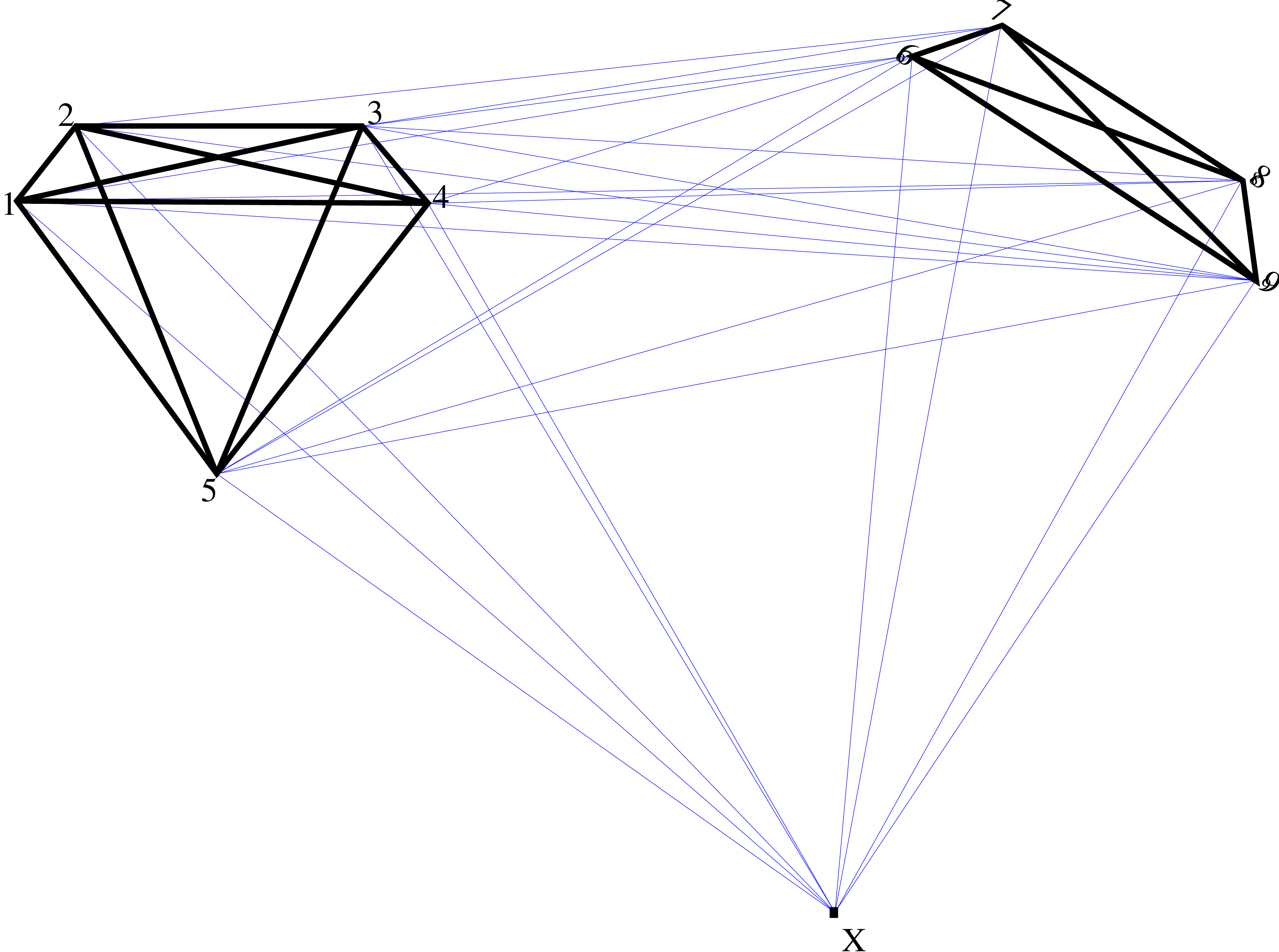}
\]
\caption{Geometric complete graph from 
a multiset of 2-coloured compositions. \label{fig:puppyears}}
\end{figure}

\beff{Classification: first schematic}{{

\mdef \label{de:SSSN}
For $N \in \N$:
\[
\SSS_N =
\{ (p,q,\rho, s) \; | \; p<q \in \PP(N);\; \rho\in \Perm_p(q); \; s \in \PP_2(q/p)   \}
\]
%
We call the partition $p$ a partition into `nations';
and refinement $q$ a partition into `counties'.
Schematically
(neglecting only $s$ for now)
elements
of $\SSS_N$ 
can be depicted:
    \[
    \includegraphics[width=2.52\ccm]{xfig/reg10gon7.eps}
\hspace{.71\ccm}   \raisebox{1cm}{ $\leadsto$}
     \hspace{.41\ccm}
     \raisebox{-1\ccm}{   \includegraphics[width=7.9\ccm]{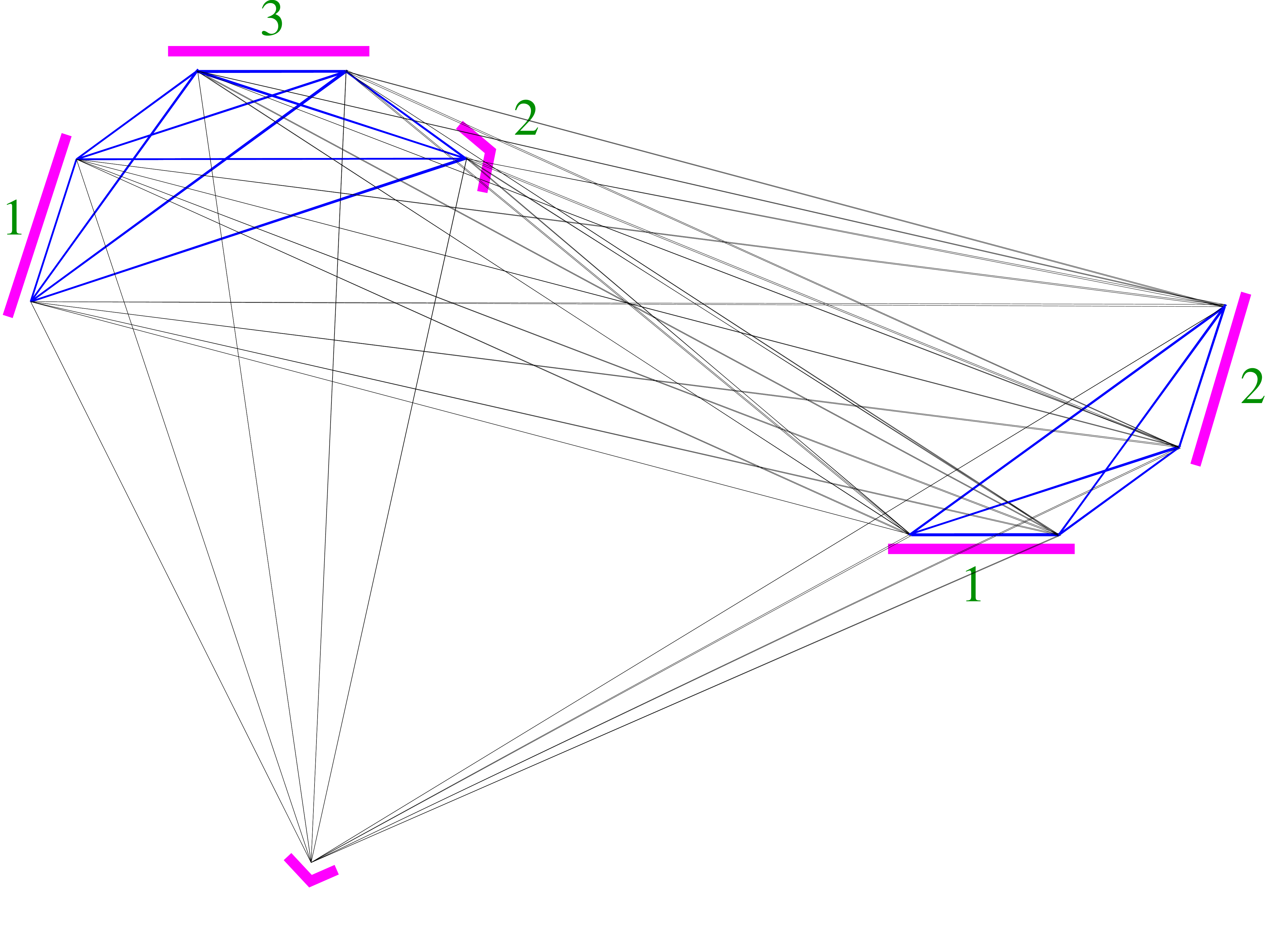} }
\]
This example has $N=10$
(here we will use Roman $\barroman{X}$ as connected symbol for vertex 10).
In the schematic the separation of vertices indicates $p$, thus
$\;$ 
$p= \{ \{\{ 1,2,3,4,5\},\{6,7,8,9\}, \{\barroman{X}\}\}$.
The coloured bars indicate $q$, thus
$\;$ 
$q=\{\{ 1,2 \},\{3,4\},\{5\}, \{6,7\}, \{8,9\}, \{\barroman{X}\}\}$.
And from the numbering
$\;$ 
$\Perm_p(q/p) = \{ (\{ 1,2 \},\{5\},\{3,4\}),( \{8,9\},\{6,7\}),( \{\barroman{X}\})
\}$.

}}


\newcommand{\ovb}{\overbrace}

\mdef
Example. For $N=2$ we have $\;$ 
 $
\SSS_2 = \{\;
\includegraphics[width=1cm]{xfig/tree012ex2aab.eps}\;,\;\;
\includegraphics[width=.51cm]{xfig/tree012ex2ab1b.eps}\;,\;\;
\includegraphics[width=.51cm]{xfig/tree012ex2ba1b.eps}\;,\;\;
\includegraphics[width=.51cm]{xfig/tree012ex2abb.eps}\;,\;\;
\includegraphics[width=.51cm]{xfig/tree012ex2bab.eps}\;,\;\;
\includegraphics[width=1.3cm]{xfig/tree012ex2acb.eps}\;
\}.
\ignore{{
\\ { } \hspace{.18in} =
\{ \;
(\; \ovb{\{\{1,2\}\}}^p, \; \ovb{\{\{ 1,2 \}\}}^q, \;\;\; \ovb{(\{1,2\})}^\rho ,
 \;\;\; \ovb{\{\{1,2\}\}}^s ),
$\\ $ { } \hspace{1.3\ccm}
(\; \{\{1,2\}\}, \{\{ 1\},\{2 \}\}, (\{1\},\{2\}), \{\{1,2\}\}), \;
$ \\ $ { } \hspace{1.3\ccm}
(\; \{\{1,2\}\}, \{\{ 1\},\{2 \}\}, (\{2\},\{1\}), \{\{1,2\}\}), \;
$ \\ $ { } \hspace{1.3\ccm}
(\; \{\{1,2\}\}, \{\{ 1\},\{2 \}\}, (\{1\},\{2\}), \{\{1\},\{2\}\}), \;
$ \\ $ { } \hspace{1.3\ccm}
(\; \{\{1,2\}\}, \{\{ 1\},\{2 \}\}, (\{2\},\{1\}), \{\{1\},\{2\}\}),
$ \\ $ { } \hspace{1.3\ccm}
(\{\{1\},\{2\}\}, \{\{ 1\},\{2 \}\}, (\{1\},\{2\}), \{\{1\},\{2\}\})
\}
}}
$
\\
To enumerate
we first vary over nation partitions
$\PP(2) = \{ \{ \{1,2\}\}, \{\{1\},\{2\}\} \}$. 
For the single-nation $p =  \{ \{1,2\}\} $
cases we  
have two possible refinements into counties,
again given by $\PP(2)$.
Here each county partition has two total orders; and there are two possible $s$
partitions.
Finally we have
an element with 
two nations, each a singleton so the other components are forced.

}}






\module{braid-theorem1}

\ignore{{  
\newcommand{\vi}{v}
\newcommand{\vj}{w}
\newcommand{\offmat}{\mat{cc} 0&1\\1&0\tam}
\newcommand{\omat}{\mat{cc} 1&0\\0&1\tam}
\newcommand{\Rec}{{\mathsf R}}     
\newcommand{\Cpq}{\C^{(p,q)}}
\newcommand{\germ}{\SSS^\C}
\newcommand{\grm}{\TTT^\C}
\newcommand{\Cx}{\C^{\times}} 
}}

\befs{20}{Recipe for solutions from $\SSS_N$}{{

\ignore{{

\mdef
    For $N \in \N$:
\[
\SSS_N =
\{ (p,q,\rho, s) \; | \; p<q \in \PP(N);\; \rho\in Perm(q); \; s \in \PP_2(q)   \}
\]
}}

\mdef
Next we need a recipe for constructing
a variety of
elements of $\Match^N(2,2)$
that give solutions $\FF(\sigma)$
from an element of $\SSS_N$. 
An element of $\SSS_N$ can be viewed in various 
ways.
Firstly it encodes a partition $p$ of the vertices into the individual
diagrams --- `nations';
and a refinement $q$ of that partition into rows --- `counties'.
It also gives a total order $\rho$ on the rows in each nation,
from top to bottom;
and a partition $s$ of them into two parts, according to shaded or not.

This formulation $(p,q,\rho,s)$ will control the varieties
upon which our representations will lie.



\mdef
The following example is in $\SSS_{12}$
(but not $\TTT_{12}$):
\beq \label{eq:onegnat1}
\lambda = \;\;
\raisebox{-.2in}{
\includegraphics[width=3.15cm]{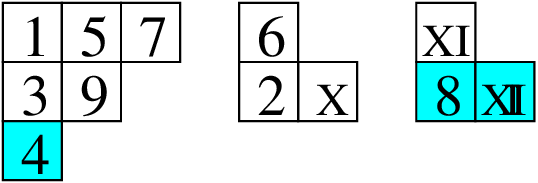}}
\eq
The $(p,q,\rho,s)$ form of  $\lambda$
as in (\ref{eq:onegnat1})
has
$ \;  p =\{ \{ 1,3,4,5,7,9 \}, \{2,6,
               X\},\{ 8,X\!\!I,X\!\!I\!\!I\}
           \} \;$
and
$\; q=\{\{ 1,5,7 \}, \{ 3,9 \}, \{ 4\},
        \{6\},\{2,X\},\{  X\!\!I\},\{8,X\!\!I\!\!I\}
\}$; 
and
the  order $\rho$ of the counties in the first nation
happens to be the nominal children-first order;
and, just focussing on the first nation,
$s$ given by lumping the first two counties together.
%


\mdef \label{de:typespace}
Consider $\lambda \in \SSS_N$ in the form $(p,q,\rho,s) $. 
Order/number
parts of $p$
as in (\ref{pa:child1st}).
Consider the complete graph $K_p$ on the ordered set $p$.
We associate
a variable  $\mu_{ij}$ to each edge, i.e. to each pair  
$p_i, p_j \in p$ with $i<j$;
and a  variable $\alpha_i$ to each vertex $p_i$;
and a  variable $ \beta_i $ to each vertex
$p_i \in p$
containing more than one part of $q$.
(NB the given {\em names} of variables
serve to distinguish them, but
have no extrinsic significance.)

The {\em type-space}
$\Cpq$
of  $(p,q,\rho,s) \in \SSS_N$ is the space whose elements
give an evaluation of each variable to an element of $\Cx$,
a nonzero complex number, such that $\beta_i \neq -\alpha_i$.
\ignore{{
thus each $\mu_{ij} \map
consist of a nonzero
complex number $\mu_{ij}$ for each distinct pair  
$p_i, p_j \in p$;
and a  nonzero complex number $\alpha_i$ for each $p_i$,
and a  nonzero complex number $ \beta_i \neq -\alpha_i$ for each
$p_i \in p$
containing more than one part of $q$.
}}
That is the complex space of variables excluding only the variety
$\prod_{ij} \mu_{ij} \prod_i \alpha_i \beta_i (\alpha_i +\beta_i)$. 
\\
A {\em germ} is a pair consisting of an element   $(p,q,\rho,s) \in \SSS_N$
and an element from $\Cpq$.
Write $\germ_N$ for the set of germs.
Write  $\grm_N$ for the subset of germs derived from
$(p,q,\rho,s) \in \TTT_N$.

\newcommand{\xx}{(\mu_- , \alpha_- , \beta_- )} 

\mdef \label{pa:rec}
The recipe
\[
\Rec: \germ_N \rightarrow \Match^N(2,2)
\]
for constructing  
elements of $\Match^N(2,2)$
(that we shall prove later are all solutions to the braid conditions)
from $(p,q,\rho,s) \in \SSS_N$
and $\xx \in \Cpq$
is as follows.
\\
\ignore{{
(Strictly speaking our target here is  $\Match^N(2,2)$ over a ring with 
indeterminates as indicated below; and we pass to $\C$ by a free choice of non-zero
evaluation of each.)
\\
}}
(i)
Each edge
{of $K_N$}
that is between nation $i$ and nation $j$ is decorated with
$\off_{\mu_{ij}} = \mu_{ij}\offmat$.

\noindent (ii)
Each vertex
in nation $i$
is 
in one of the two parts of $s_i\in\PP_2(q|p_i)$ obtained from $s$
(partitioning the counties of each nation).
Simply for definiteness in naming conventions, we order the two parts (using $\rho$).
For vertices in the first 
part in $s_i$ 
we decorate with $\alpha_i$; 
and otherwise with $\beta_i$.
(Note that the decoration  $x\in\{ \alpha_i, \beta_i \}$, say,
is the same on all vertices 
in the same county.)
\\ (iii)
Each edge between vertices in the same county
is decorated with 
$ x\omat$ where $x\in\{ \alpha_i, \beta_i \}$
is the common vertex decoration.
(We 
abbreviate to  
symbol 0 for this edge decoration.) 

\noindent (iv)
The remaining edges ---
between vertices in different counties in the same nation are decorated as
follows.
If the order on counties given by $\rho$
agrees with
the natural order
(e.g. if $\rho$ places vertex $\vi$ before $\vj$ and also $\vi < \vj$)
then the $\vi\vj$
edge is +; otherwise it is $-$.
Finally for each signed edge, if it is $+$ (resp. -)
and the end vertices are 
both $\alpha_i$
then
decorate with $\fii_{\beta_i}$ (resp. $\mfii_{\beta_i}$);
or both $\beta_i$
then
decorate with $\fii_{\alpha_i}$ (resp. $\mfii_{\alpha_i}$);
else if end vertices are different then
decorate with $\ai$ (resp. $\mai$).
See 
Fig.\ref{fig:N2}
for the matrix implications of these decorations.


}}

\begin{figure}         
\[
\includegraphics[width=6cm]{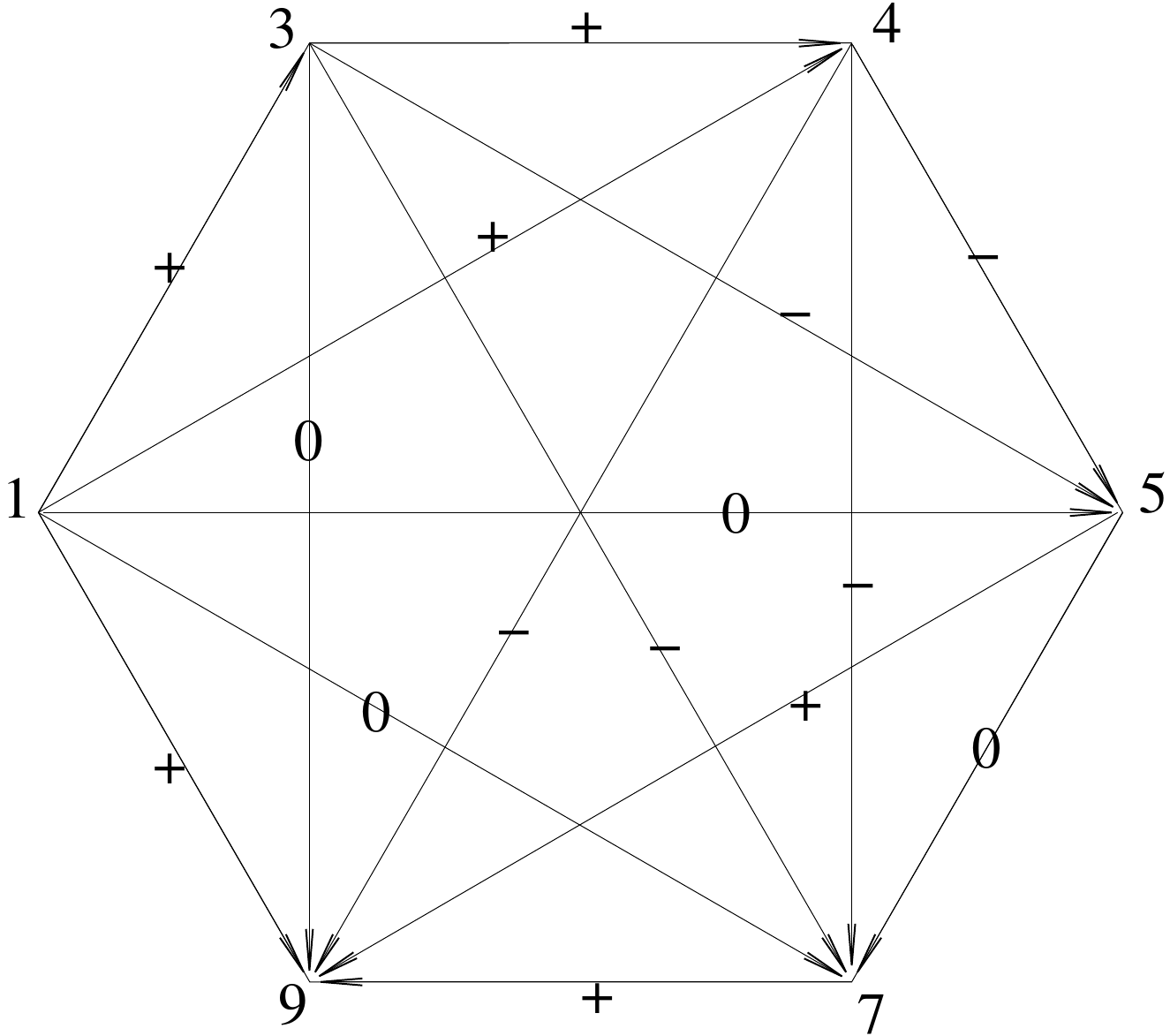}
\hspace{1.7591\ccm}
\includegraphics[width=6cm]{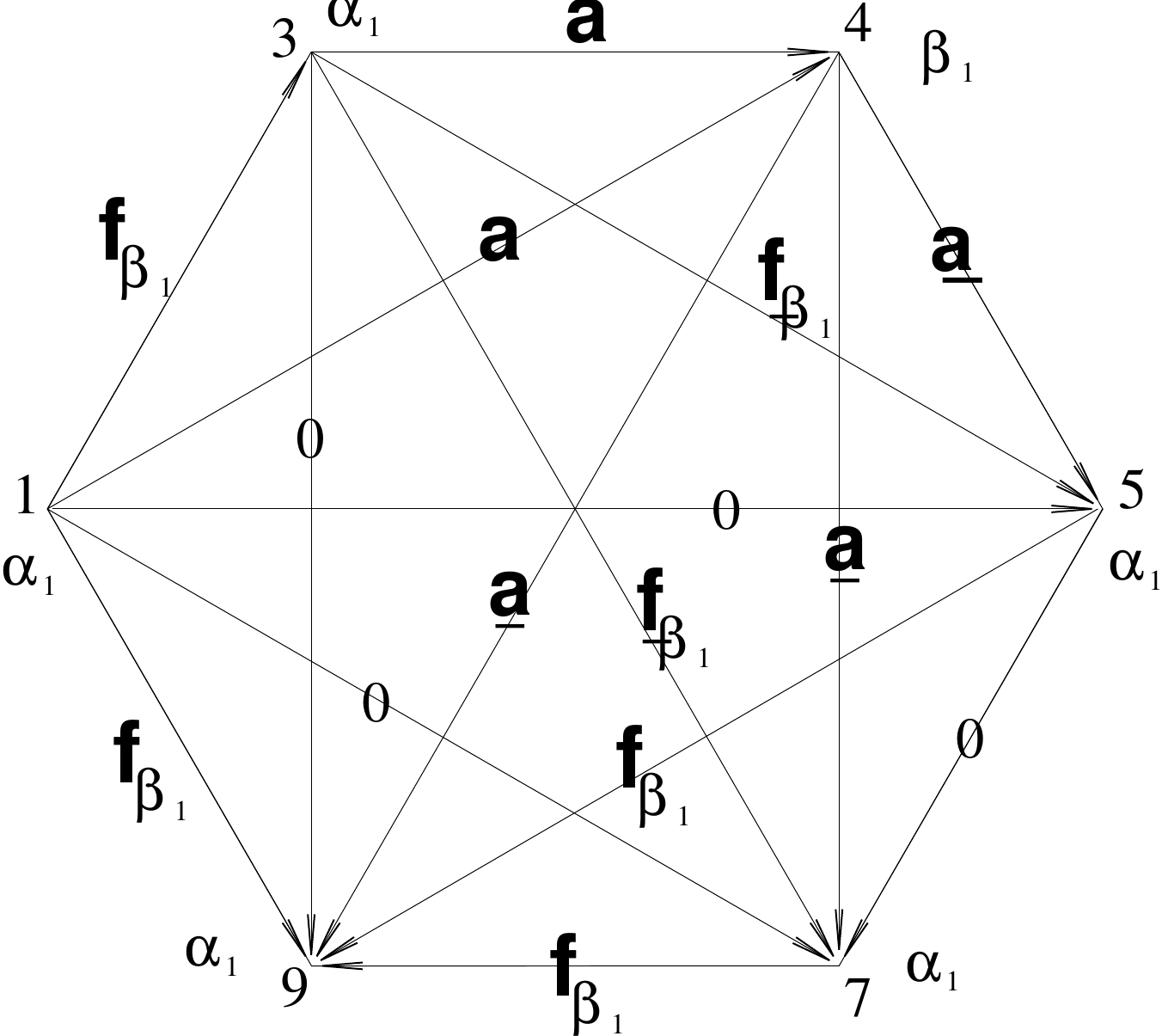}
\]
\caption{Constructing a 
part-solution
from the first nation in (\ref{eq:onegnat1}).
\label{fig:hexa321a}}
\end{figure}

\mdef
Example.
Considering just the first nation in $\lambda$ in (\ref{eq:onegnat1}),
then we have the 
edge signs
as on the left in Fig.\ref{fig:hexa321a}. 
Thus for {example} the 35 edge is $-$
because $3<5$ 
but 5
comes first in the chosen county order.
\ignore{
Finally for $s$ 
we lumped
the first two counties together:
$s = \{\{1,5,7,3,9\},\{4 \}\}$ (restricting here only to the first nation). 
Then t}
The vertex
and f/a
variables
are as shown on the right 
in Fig.\ref{fig:hexa321a}.


\beff{Theorem}{{ }}

\begin{theorem} \label{th:mainx}

For each $N \in \N$ the
set
$\Rec(\germ_N)$ gives a complete set 
of charge-conserving braid representations
up to  
\Xequivalence. 

\ppmx{-need to say what this is - from X-Lemma}
\ignore{{\\
\ppm{Task: I have demoted it, but I think
Eric wants to include here the construction of
a transversal of the symmetry actions. We can do this.}

\\
A1. Pick a partition
$\lambda = (\lambda_1, \lambda_2, ..., \lambda_m )$ of $N$ and hence a partition of
$\ul{N} = \{1,2,...,N \}$ into
`nations',
$$
\ppp(\lambda) = \{
\{ 1,2,...,\lambda_1 \},
\{ \lambda_1 +1, \lambda_1 +2, ..., \lambda_1 + \lambda_2 \},
\hspace{1.65in}$$ $$\hspace{1.65in}
\{ \lambda_1  + \lambda_2+1, \lambda_1  + \lambda_2+2, ..., \lambda_1  + \lambda_2+ \lambda_3 \},
..., \{ ..., N \}\}
$$
\\ 2. Pick a partition of each nation into counties.
Thus we have a refinement $(\ppp(\lambda),q)$.
(Note that $q$ includes a partition, $q(i)$ say, of each nation.)
\\ 3. Pick a total order $\orr$ on counties $q(i)_j$ in each nation.
(For definiteness the
arbitrary labelling
$j$ order can be given by a convention, and the total order then
given in terms of this.)
\\ 4. Pick a partition $r$ of counties in each nation into two parts.
\\ B1. Form a $K_N$-configuration and hence a matrix
from $(\ppp(\lambda),q,\orr, r)$
as follows.
\\
...}}
\ppmx{paul improve this. recall/cite I set and define a controlling map
from configs in MatchN22 to solutions; and use to state theorem cleanly.}

\end{theorem}



\noindent
Outline of proof.  
(The proof
is in Sections~\ref{section: the proof}-\ref{ss:proofend}.)
Theorem~\ref{th:constraints}
formulates sufficient conditions for a solution
{- noting Lemma~\ref{lem:3 suffice}}.
We solve these conditions in rank 2 in Prop.\ref{pr:L2class2}
and in rank 3 in Prop.\ref{pr:N3}.
Then
for general rank
we will need  
Lemmas exploiting some `magic' in these solutions
- see \S\ref{ss:proofend}.
\soutx{The  
subsequent proof of \ref{th:JJN} in 
(\ref{pa:proofJJN}),
(\ref{pa:proofJJN0})
is then relatively straightforward.}
$\;$ The following is a corollary.

{{ }}   


\begin{theorem} \label{th:JJN}

\newcommand{\JJ}{{\mathsf J}} 

\ppmx{[drafting here for comparison. to be moved later.]}
\soutx{{ The set
$\TTT_N$
{is the injective image of $\JJ_N(\{1,2,3\}^*)$ in $\SSS_N$ and }
gives a transversal of the set of
orbits of the $\Sym_N$ action on $\SSS_N$. 
}}
For each $N\in\N$
the set $\TTT_N$ gives a transversal 
of the set of varieties
of  
charge-conserving braid representations
up to $\Sym_N$
(action induced from (\ref{lem:chacon1})(II))
and  
\Xequivalence.
That is $\Rec(\grm_N)$ gives a complete set of representations up to
$\Sym_N$ and \Xequivalence.

\end{theorem}


Note that the size $|\TTT_N|$ is given by the Euler transform
(see e.g. \cite{SloanePlouffe95}) of sequence $3^{N-1}$, by
(\ref{41}) and (\ref{44}).

\begin{remark}
As noted in Section~\ref{ss:target},
if $\VVV$ is the target category then
isomorphisms in the functor category $\functor(\VVV,\VVV)$ induce isomorphisms
in our representation category $\functor(\Bcat,\VVV)$.
Thus in our case
the action of the symmetric group
$\Sym_N \hookrightarrow \functor(\VVV,\VVV)$
leads,
via $\Rec(\germ_N)$,
to a decomposition of $\germ_N$ into orbits.
In this way we may access all solutions from a  
smaller
transversal
set, $\TTT_N$,  
via the $\Sym_N$ action.
(We can also use the flip symmetry $\PPP$
to reduce
slightly
further - see \S\ref{ss:flyp}.)
\end{remark}




\module{higgs}


\newcommand{\KK}{\ul{K}}
\newcommand{\GGG}{\ul{G}}
\newcommand{\III}{\ul{\II}}
\newcommand{\homm}{\ul{\hom}}


\mdef
\ppmx{[strictly speaking this could be later, not needed before Thm.?]}
A {\em (Higgs) \configuration} of $K_n$ is an assignment,
cf. (\ref{pa:config}),
of a scalar
from $GL_1(\C)$ to each vertex;
and matrix from $GL_2(\C)$ to each edge,
so that each $K_2$ subgraph is assigned
from the set of six types
$\II = \{ \off, 0, \uppfii, \mfii, \uppai,\mai \} $
as in
Fig.\ref{fig:N2}.

We observe  
from {Proposition}~\ref{pr:L2class2}
and \ref{lem:chacon1}(III)
that
a necessary condition for an element of $\Match^N(2,2)$
to give a braid representation is that it gives an \configuration\ of $K_N$.

Write $\homm(\KK_N,\GGG)$ for the set of \configuration s of $K_N$
($\KK_N$ is $K_N$ as a complex and $\GGG=(GL_1,GL_2)$).
Note that the $K_2$ subgraphs of $K_N$ are in bijection with  the edges.
Define
$$
\III : \homm(\KK_N,\GGG) \rightarrow \Setcat(K_N^1 , \II )
$$
by mapping each $K_2$ configuration to its type
in the obvious way.
Observe then that there is a (possibly empty) fibre of solutions over each element
$\alphaa$ of
$\Setcat(K_N^1 , \II )$.
More generally the fibre $\III^{-1}(\alphaa)$ consists of the possible choices of
the parameters in each matrix.
\\
Remark: note that there is no {\it ab initio} guarantee here that there are
{\em any} solutions in ranks above two,
even before imposing higher-rank constraints,
since the $K_2$-subgraphs of $K_N$  
overlap with each other.





\section{Proof of Theorem~\ref{th:mainx}: first steps}\label{section: the proof}
In the same spirit as Fig.\ref{fig:N2}  
we may express solutions for $N=3$ schematically on triangles. 

\module{summary03}


\medskip

\newcommand{\signa}[1]{signature: \; #1}



\beff{}{{
\renewcommand{\ooo}[1]{}   
\mdef \label{pr:N3} 
{\bf Proposition}.
(`\ninerule') {\it
For $N=3$
the following
types of 
triangle
configurations
are allowed
for a charge-conserving functor $\FF$
(showing one per $\SymXZ_3$ orbit
{as defined  in (\ref{restriction-and-symmetries})}):
}
\[
\ber{2}{3}{\;\;\;/_{\mu} }{\;\; /_{\nu} }{\;/_{\lambda} \;}
\hspace{.1in}
\]
\[
\bermss{1}{2}{1}{\;/_{\mu}\;}{\;/_{\mu}\;}{\uppfii_\beta} \hspace{.1in}
\hspace{.1241in} 
\bermss{1}{2}{3}{\;/_{\mu}\;}{\;/_{\mu}\;}{\uppai}
\hspace{.21in}
\bermss{1}{2}{1}{\;/_{\mu}\;}{\;/_{\mu}\;}{0} 
\]
\vspace{.05in}
\[ \hspace{.014in}
\bermpp{1}{1}{1}{\uppfii_\beta \;}{\;\uppfii_\beta \;}{\;\uppfii_\beta \;}
\hspace{.24in}
\bermpp{1}{2}{1}{\uppai}{\uppai}{\;\uppfii_{a_2} \;}
\hspace{.24in}
\bermpp{1}{2}{2}{\uppai}{\uppfii_{a_1}}{\; \uppai\;}  
\]
\vspace{.041in}
\[
\bermxx{1}{2}{2}{\;\pai\;}{0}{\pai}{\abcac}
\hspace{.21in}
\bermxx{1}{1}{1}{\pfii_\beta \;}{0}{\pfii_\beta}{\abcac}
\hspace{.21in}
\ber{1}{1}{0}{0}{0}   
\]
}}


\befff{}{{
\prooff { } 
From the $N=2$ result in Prop.\ref{pr:L2class2},
a solution must have an 
edge decoration type from the set
$\{ 0,\off,\pfii,\pai,\mfii,\mai \}$ on each edge of the triangle.
Thus there are $6^3$ such combinations to try
with the constraints from Theorem~\ref{th:constraints}. 
As we show
in Lemmas~\ref{lem:slash}-\ref{lem:sundry0s},
various
little
`miracles' occur,
so that within each of these finitely many cases the analysis becomes tractable,
yielding the claimed configurations.
\\
(This proof will conclude with the proof of Lemma~\ref{lem:sundry0s} below.)
\qed
}}

\ignore{{
\beff{}{{
\noindent {\bf Proposition}.
For $N=3$
the following triangles are allowed (showing one per $S_3$ orbit?):
\[
\ber{2}{3}{\;\;\;\;/_{\pm\mu} }{\;\;\;\; /_{\pm\nu} }{\;/_{\pm\lambda} \;}
\hspace{.1in}
\hspace{-.02471in} signature: \; (1^9),6
\]
\ignore{{
- here there are (up to) 9 distinct eigenvalues
($a_1, a_2, a_3, \pm\mu, \pm\nu,\pm\lambda$),
each of multiplicity 1;
given by 6 free parameters ($a_1, a_2, a_3, \mu, \nu,\lambda$), as shown.
We record this as the `signature'.}}
\[
\bermss{1}{2}{1}{\;/_{\pm\mu}\;}{\;/_{\pm\mu}\;}{\uppfii_\beta} \hspace{.1in}
\hspace{-.41in} signature: \; (3^1 2^2 1^2),4
\hspace{.41in} 
\bermss{1}{2}{3}{\;/_{\pm\mu}\;}{\;/_{\pm\mu}\;}{\uppai} \hspace{.1in}
\hspace{-.41in} signature: \; (2^4 1^1),4  
\]
\[
\bermss{1}{2}{1}{\;/_{\pm\mu}\;}{\;/_{\pm\mu}\;}{0} 
\hspace{-.41in} \signa{(4^1 2^2 1^1),3}
\]
\[
\bermpp{1}{1}{1}{\uppfii_\beta \;}{\;\uppfii_\beta \;}{\;\uppfii_\beta \;}
\hspace{-.41in}
\signa{(6^1 3^1),2} \hspace{.4in}
\bermpp{1}{2}{1}{\uppai}{\uppai}{\;\uppfii_{a_2} \;}
\hspace{-.41in} \signa{(5^1  4^1 ),2}
\]
\ppm{-note we force $\beta=a_2$.}
\[
\bermxx{1}{2}{2}{\;\pai\;}{0}{\pai}{\abcac}
\hspace{-.41in} \signa{(6^1 3^1),2} \hspace{.4in}
\bermxx{1}{1}{1}{\pfii_\beta \;}{0}{\pfii_\beta}{\abcac}
\hspace{-.41in} \signa{(7^1 2^1),2}
\]
\[
\ber{1}{1}{0}{0}{0}   
\hspace{-.41in} \signa{(9^1),1}   
\]
}}
}}

\ignore{{

\whisper{
Note 
that we have a formal degeneracy of aa0 with fff above - however
in the fff case above there are three non-trivial blocks if we do not
have an inverse to $a_1 - \beta$, whereas for aa0 there are two such blocks.
---Also need to put in the other type from this case.
}

}}



\ignore{{
Note in particular that single / does not appear...
%
%
\ppm{The six rows of types here correspond to a (very interesting) partial classification used
by Eric to assist counting.}
}}

\ignore{{

\beff{}{{
We should note what the symmetry action does to these.
Applying (12) say, we `flip' the triangle at the 12 edge
so this edge is reverses and the others are swapped.
Examples:
\[
\bermxx{1}{2}{1}{\mai}{\pai}{0}{} \hspace{.21in}
\bermxx{1}{1}{1}{\mfii}{\pfii}{0}{} 
\]
thus the = annotation must be treated with caution --- it becomes a
kind of twisted-=.

}}
}}

\ignore{{
{Cf. a full `equilateral'
triangle symmetry we might have assumed gets taken here:}
 \[
\bermss{1}{2}{1}{\pai}{\pai}{0} \hspace{.21in}
\bermss{1}{1}{1}{\pfii}{\pfii}{0}
\hspace{.51cm}\mbox{(--- NB  INvalid move)}
\]
That is, we do not have full equilateral symmetry in general, and need
to keep track of either long/short edge data or possibly arrow directions.
}}

\ignore{{
\ppm{Time now perhaps to take these thoughts over to try them on
Eric's tetrahedron results?}

A possible notation is to ignore the edge arrows for / and 0; and to
draw the arrow reversed with $\pfii$ resp $\pai$ instead of writing out $\mfii$ and $\mai$.
}}




\module{slashlemma}

\begin{lemma}\label{lem:slash} 
	Let $\FF:\Bcat \rightarrow\Match^N$.
	If there is one $/$ among the edge labels for a triple $i<j<k$
        then there are at least 
        two such edges.
\end{lemma}
\begin{proof}
  Without loss of generality we may assume
  the labels are $1,2$ and $3$, with the edge between $1$ and $2$ being a /.
  Thus $d_{12}=a_{12}=0$ and $c_{12}b_{12}\neq 0$.
  Suppose (for a contradiction) that 
  the $13$ and $23$ edges are not
  /.
  Equations (\ref{eq:34xx0}) are: $c_{12} (d_{13} d_{23} -d_{12} d_{23} -a_{12} d_{13})=0$ and 
  $c_{12} (- a_{13} a_{23} +a_{12} a_{23} +d_{12} a_{13})=0$ which imply $d_{13} d_{23}=a_{13} a_{23} =0$.
  By assumption, this leaves two choices:
  $d_{13}=a_{23}=0$
  with
  $a_{13}d_{23}\neq 0$; or vice versa, since otherwise we have another / edge.
  In either case we have $b_{13}c_{13}\neq 0$ and $b_{23}c_{23}\neq 0$.
  Applying the permutation $(132)$ to the second equation above we find:
$-b_{13}(d_{23}a_{12}-d_{13}a_{12}-a_{13}d_{23})=0$.  But since $b_{13}\neq 0$ and $a_{12}=0$ we see that $a_{13}d_{23}=0$, so that the first choice is impossible.  Applying $(132)$ to the first equation above yields: $-b_{13}(a_{23}d_{12}-a_{13}d_{12}-d_{13}a_{23})=0$ which eliminates the second choice.
\end{proof}

\module{reps-03110}



\subsection{Cases for $N=3$} $\;$ 
\begin{table}
\[
\begin{array}{|l|lllllll|} \hline
\off\off\off &&&&&&&
\\ \hline
\off\off\doo \qquadd \off\off\upp \qquadd \off\doo\off \qquadd
\off\upp\off \qquadd \upp\off\off \qquadd \doo\off\off  &&
\\ \hline
\off\doo\doo \qquadd \off\upp\upp \qquadd \doo\off\upp \qquadd
\upp\off\doo \qquadd \upp\upp\off \qquadd \doo\doo\off &&
\\ \hline
\off\doo\upp \qquadd \off\upp\doo \qquadd \doo\off\doo \qquadd
\upp\off\upp \qquadd \upp\doo\off \qquadd \doo\upp\off &&
\\ \hline
\doo\doo\doo \qquadd \upp\upp\upp \qquadd \upp\doo\doo \qquadd
\doo\upp\upp \qquadd \doo\doo\upp \qquadd \upp\upp\doo &&
\\ \hline
\doo\upp\doo \qquadd \upp\doo\upp &&&& &&
\\ \hline
0\off\off \qquadd \off 0\off \qquadd \off\off 0 &&& &&
\\ \hline
0\off\upp \qquadd 0\off\doo \qquadd 0\upp\off \qquadd
0\doo\off \qquadd \off 0\doo \qquadd \off 0\upp \qquadd
\doo 0\off \qquadd ...
\\ \hline
0\doo\upp \qquadd 0\upp\doo \qquadd \doo 0\doo \qquadd
\upp 0\upp \qquadd \upp\doo 0 \qquadd \doo\upp 0 &&
\\ \hline
0\doo\doo \qquadd 0\upp\upp \qquadd \upp 0\doo \qquadd
\doo 0\upp \qquadd \upp\upp 0 \qquadd \doo\doo 0  &&
\\ \hline
00\off \qquadd 0\off 0 \qquadd \off 00 &&& &&
\\ \hline
00\upp \qquadd 00\doo \qquadd 0\upp 0 \qquadd 0\doo 0 \qquadd \upp
00 \qquadd \doo 00 &&
\\ \hline
000 &&&&& &&
\\ \hline
\end{array}
\]
\caption{Level $\NL=3$. The $4^3$ $\HH$-types arranged
so that each row is a $\GG_3$-Orbit.
(The first column thus gives a transversal.)
\label{tab:1}}
\end{table}
\ignore{{  
\begin{table}[h]
\[
\begin{array}{|llllllll|} \hline
\off\off\off &&&&&&&
\\ \hline
\off\off\doo \qquadd \off\off\upp \qquadd \off\doo\off \qquadd
\off\upp\off \qquadd \upp\off\off \qquadd \doo\off\off  &&
\\ \hline
\doo\doo\doo \qquadd \upp\upp\upp \qquadd \upp\doo\doo \qquadd
\doo\upp\upp \qquadd \doo\doo\upp \qquadd \upp\upp\doo &&
\\ \hline
0\off\off \qquadd \off 0\off \qquadd \off\off 0 &&& &&
\\ \hline
0\doo\doo \qquadd 0\upp\upp \qquadd \upp 0\doo \qquadd
\doo 0\upp \qquadd \upp\upp 0 \qquadd \doo\doo 0  &&
\\ \hline
000 &&&&& &&
\\ \hline
\end{array}
\]
\caption{
Level $\NL=3$.
The
6 classes containing the 23 
admissible $\HH$-types arranged by $\GG_3$-Orbit. 
\label{tab:1ad}}
\end{table}
}}

Recall that 
if $F(\sigma)\in\Match^N(2,2)$
gives a solution then so does its image under the action of the group
$\GG_\NL = \Zz_2 \times \Sigma_N$ from (\ref{restriction-and-symmetries}).
Observe that the assignment within the set  $\pfii,\pai$ to an edge commutes with
the $\GG_\NL = \Zz_2 \times \Sigma_N$ action,
up to $\pfii \leftrightarrow \mfii$.
Accordingly we may  
first group together the edge labels $\pfii,\pai$ as type-+;
and $\mfii,\mai$ as $-$.
So now the set of label types is $\HH = \{ 0 ,\off,+,- \}$.
With $N=3$,
all $4^3$ possible triangle $\HH$-configurations are then shown in  Table~\ref{tab:1}.

\newcommand{\fHH}{\ul{\HH}} 
\newcommand{\fII}{\ul{\II}} 

\ppmx{[probably need to say a bit more about the $\HH$ and $\II$ function notations that we use
repeatedly below! We used to have a commuting diagram defining them, but not in this
document (yet). Since $\HH$ and $\II$ are the label sets, perhaps underline them to notate
the functions? (as in (5.19)?)...]}

For $N=3$
we define $\fHH$ as the map from solutions to $\HH^3$ --- triples of edge labels from $\HH$;
and $\fII$ is the corresponding map to $\II^3$.
(And similarly for general $N$.)

%
Since  solutions belong to $\GG_\NL$-orbits
and $\fHH$-fibres  
we can
characterise all solutions 
by giving 
those
for the
fibres of a transversal. 
Applying  
to $N=3$,
Table~\ref{tab:1} is thus organised by orbits.




\begin{lemma}
 Table~\ref{tab:1} gives the complete list of $\HH$-types
for $N=3$
arranged by
$\GG_3$-orbit.
\end{lemma}
\proof
Both distinctness and saturation of the $4^3$ bound are clear from the arrangement.
\qed
\ignore{{
But it might be worth spelling out in prepartion for higher ranks.
The number of 0s and number of /s is fixed in an orbit.
We organise firstly by the number of 0s --- in this case from none to
three.
Within each such grouping we organise by number of /s --- in this case
decreasing. The remaining entries are $\pm$s and  we
have not used any particular order, but there are not many cases.
The orbits follow
from applying \ppm{...the actions of generators as in .... .}
}}

\newcommand{\bonet}[1]{\lambda_{#1}^2}
\newcommand{\bl}{\lambda}

\begin{lemma} \label{le:///}
The fibre of solutions $\fHH^{-1}(\off,\off,\off)$ is non-empty.
The solutions $F(\sigma)$  
may be written here,
up to $X$ equivalence (\ref{lem:Xeq})
in the form
\[
\alpha(F(\sigma)) = (a_1, a_2, a_3,
  \mat{cc} 0 & \bonet{12}\\ 1 & 0 \tam ,  \mat{cc} 0 & \bonet{13}\\ 1 & 0 \tam
  ,
  \mat{cc} 0 & \bonet{23}\\ 1 & 0 \tam )
\]
\end{lemma}
\proof
Here $a_{ij} =0$ and $d_{ij}=0$ so all the constraints from
\eqref{eq:31xx0}-\eqref{eq:35x0} 
are satisfied trivially.
For invertibility, $b_{ij} c_{ij} \neq 0$.
{Thus $b_{ij}, c_{ij}$ are non-zero.
Using the X~Lemma there is no loss
of generality in setting the $c_{ij}$s to 1.}
\qed

\ignore{{
\noindent Here there are (up to) 9 distinct eigenvalues
($a_1, a_2, a_3, \pm\mu, \pm\nu,\pm\lambda$),
each of multiplicity 1;
given by 6 free parameters ($a_1, a_2, a_3, \mu, \nu,\lambda$), as shown.
We record this as the `signature'.
}}

The existence of this 6-dimensional variety is betokened in
Prop.\ref{pr:N3} by the following figure:   
\[
\ber{2}{3}{\;\;\;\;/_{\mu} }{\;\; /_{\nu} }{\;/_{\lambda} \;}
\]



\medskip

\ppmx{For cases involving + Eric's successful approach uses $\HH$, and gives rise
to some very interesting results. Paul's approach uses $\II$ and so is
the route towards complete classification. For these reasons we discuss both.}

\ignore{{
\begin{lemma} \label{lem:offoffupp}
  The fibre of solutions $\HH^{-1}(\off,\off,\upp)$ is non-empty.
And hence so are the fibres in the orbit, as given in Table~\ref{tab:1}.
If $\HH(F(\sigma)) = (\off,\off,\upp) $ then $F(\sigma)$ is
  given, up to the basis scaling as in \ref{X},
by
\ignore{{
\[
\alpha(F(\sigma)) =
  (a_1,a_2,-b_{23},
  \mat{cc} 0 & b_{12}\\ 1 & 0 \tam ,
  \mat{cc} 0 & b_{12}\\ 1 & 0 \tam ,
  \mat{cc} c_{23}-b_{23} & b_{23}\\ c_{23} & 0 \tam
  )
\]
where either $a_2 = a_3 = -b_{23}$ or $a_2 = c_{23}$. 
\ppm{-but now what happens if we rescale?
should be able to replace $c_{23}\leadsto 1$ and
$b_{23}=\leadsto b_{23} c_{23}$. but this seems odd..}
\\
Paul version:
}}
\[
\alpha(F(\sigma)) =
  (a_1,a_2,a_{3},
  \mat{cc} 0 & \bonet{12} \\ 1 & 0 \tam ,
  \mat{cc} 0 & \bonet{12} \\ 1 & 0 \tam ,
  \mat{cc} a_3+\bl & -a_3 \bl \\ 1 & 0 \tam
  )
\]
where either $a_2 = a_3 $ or $a_2 = \bl$;
and 
(since the lefthand-side is $a_{23} \neq 0$ by assumption, which we
will use in the proof), $a_3 +\bl \neq 0$. 
\end{lemma}
}}



\ppmx{Let's try to state it again in the $\II$ form
(keep in mind the convention (1,2,3,12,13,23)):}

\begin{lemma} \label{le://+}
Solutions of type $\off\off+$.
The fibre of solutions $\fII^{-1}(\off,\off,\pai)$,
  respectively $\fII^{-1}(\off,\off,\pfii)$, is non-empty.
And hence so are the fibres in the orbit, as given in Table~\ref{tab:1}.
Specifically,
if $\fII(F(\sigma)) = (\off,\off,\pfii) $
then $F(\sigma)$ is given, up to X equivalence (\ref{lem:Xeq}), by
\ignore{{
\[
\alpha(F(\sigma)) =
  (a_1,a_3,a_{3},
  \mat{cc} 0 & \bonet{12} \\ 1 & 0 \tam ,
  \mat{cc} 0 & \bonet{12} \\ 1 & 0 \tam ,
  \mat{cc} a_3+\bl & -a_3 \bl \\ 1 & 0 \tam
  )
\]
or maybe organisationally better convention
}}
\[
\alpha(F(\sigma)) =
  (a_1,a_2,a_{2},
  \mat{cc} 0 & \bonet{12} \\ 1 & 0 \tam ,
  \mat{cc} 0 & \bonet{12} \\ 1 & 0 \tam ,
  \mat{cc} a_2+\bl & -a_2 \bl \\ 1 & 0 \tam
  )
\]   
where
$a_1, a_2, \lambda_{12}, \lambda \neq 0$ and 
{$a_2 +\bl \neq 0$}. 
$\;$
If $\fII(F(\sigma)) = (\off,\off,\pai)$ then 
\newcommand{\bll}{a_2}
\[
\alpha(F(\sigma)) =
  (a_1,\bll,a_{3},
  \mat{cc} 0 & \bonet{12} \\ 1 & 0 \tam ,
  \mat{cc} 0 & \bonet{12} \\ 1 & 0 \tam ,
  \mat{cc} a_3+\bll & -a_3 \bll \\ 1 & 0 \tam
  )
\]
where 
$a_1, a_2, a_3, \lambda_{12} \neq 0$ and 
$a_3 +a_2 \neq 0$. 
\end{lemma}


\proof
For 
$\off\off+$
we have:   $\;\;$  
$d_{ij} = 0$ and
$a_{12} = a_{13} = 0$,
\hspace{.3cm}
$b_{ij} c_{ij} \neq 0$   
and
$a_{23} \neq 0$.
\\
We must now check
\eqref{eq:31xx0}-\eqref{eq:35x0} and orbits.
\\
\eqref{eq:31xx0}-\eqref{eq:31xxx0}:
the `12' versions  are satisfied immediately;
and similarly for the `13' versions.
For the 23 versions we have firstly
(this is just retracing the relevant case from the level-2 Proposition)
$a_2^2  -a_{23} a_{2} =b_{23} c_{23}$
and
$a_3^2 - a_{23} a_3 = b_{23} c_{23}  $.
This gives $ a_2 (a_2 -a_{23}) = a_3 (a_3 - a_{23}) = b_{23} c_{23}$
and so (at most) the indicated two solutions for $(a_2 , a_3, \mate_{23})$.

\noindent
\eqref{eq:34xx0}:
All are satisfied immediately here.

\noindent
\eqref{eq:35x0}:
All are immediate, except that 
the image under the (13) perm
now implies
$b_{12} c_{12} = b_{13} c_{13} $ (since $a_{23} \neq 0$).
It follows, via the X Lemma \ref{lem:Xeq},
that 
$\mate_{12}$ and $\mate_{13}$
can be given the form as above.
%
%
\qed


\medskip

\mdef
As for pictures 
for the $\off\off+$ class we have the following.
These pictures are not exactly those for the given solutions, but note
that here the $\pai/\pfii$ can be moved anywhere using the symmetry orbit.
That is:
\ppmx{
(say what the eigenvalue story is here!)}
\[
\bermss{1}{2}{1}{/_\mu \;}{\;/_\mu }{\pfii_\lambda} \hspace{.631in}
\bermss{1}{2}{3}{/_\mu}{\;/_\mu}{\pai} \hspace{.1in}
\]
\ppmx{--- the $\mu$ and $\lambda$ are just there to give names to all the
eigenvalues ---}
respectively.
Note the `identification' of the two edge matrices.
Note the directed edges and the longer 13 edge.
\ppmx{Are these things important or not?
Yes so this is not quite right yet!!! Except that here it is, because $\off$ edges
are effectively undirected.}
\\
\ppmx{Note the new notation ferro $=\pfii$ for $+_{ii}$;
$\mfii$ for $-_{ii}$
(ferro in the sense of ferromagnetic ground state, where the state of
the vertices at either end of the edge is forced to be the same);
and antiferro $=\pai$ for $+_i$ and $\mai$ for $-_i$
(antiferro in the sense of not-ferro; really it is a misnomer and
something like `hot' --- i.e. the states at the two ends are indendent
--- would be more appropriate)
(there are macros for all this).}

\ppmx{(Task: Now need to propagate these changes below!...)}


\begin{lemma}
The fibres of solutions
$\fHH^{-1}(\off,\upp,\upp),\fHH^{-1} (0,\off,\upp ),\fHH^{-1} (0,0,\off )$
and  
$\fHH^{-1}(\off,\upp,\doo)$ are  empty.
\end{lemma}
\proof

This is an immediate consequence of Lemma \ref{lem:slash}.
\qed

\begin{lemma} \label{le:+++}
Solutions of type +++.
If $\fHH(F(\sigma)) = (+++)$
(NB expressing it this way uses up most of the symmetry, so we will 
list the three aaf cases below explicitly!)
then either
$\fII(F(\sigma)) = (\pfii,\pfii,\pfii) $   
which implies{,  up to X equivalence (\ref{lem:Xeq}),}
\[
\alpha(F(\sigma)) = (\alpha,\alpha,\alpha,\loone,\loone,\loone)
\]
or $\fII(F(\sigma))$ is one of
$ ( \pai , \pai , \pfii)$,
$ ( \pai , \pfii , \pai)$,
$ (  \pfii, \pai , \pai )$,
which imply (respectively)
\[  \alpha(F(\sigma)) = (\alpha,\beta,\beta,\loone,\loone,\loone) \]
\[  \alpha(F(\sigma)) = (\alpha,\beta,\alpha,\loone,\loone,\loone) \]
\[  \alpha(F(\sigma)) = (\alpha,\alpha,\beta,\loone,\loone,\loone) \]
(NB $\alpha+\beta\neq 0$.)
(See the pictures in \eqref{pa:+++}.) 
\end{lemma}

\ignore{{
\newcommand{\alp}{\alpha}
\newcommand{\be}{\beta}
\newcommand{\lemx}{the X~Lemma}
\newcommand{\loone}{\mat{cc} \!\alp\!+\!\be & -\alp\be\! \\ 1&0\tam  }
\newcommand{\ccc}{\bullet}
\newcommand{\trin}[6]{ \xymatrix{& #3\; \ccc\;\; \ar@{-}[dr]^{#4}
\\
#1 \ccc \ar@{-}[rr]_{#6} \ar@{-}[ur]^{#2} && \ccc #5  }}
}}


\proof
\ppmx{DOUBLE CHECK signs/symmetries now!!!}
Here $d_{ij} =0$, while $a_{ij},b_{ij},c_{ij} \neq 0$.
Considering  
\eqref{eq:31xx0}
we find
$a_1 (a_1 -a_{12})=b_{12} c_{12} = a_2 (a_2-a_{12})$
and similarly for the 13 and 23 versions.
Meanwhile all versions of
\eqref{eq:31xxx0}  
are satisfied identically.
\\
From \eqref{eq:34xx0} 
most are satisfied identically,
however we have $a_{13}=a_{12}$ from
{the P image of} \eqref{eq:34xx0}
and $a_{23}=a_{13}$ from \eqref{eq:34xx0}(i).
\\
From the \eqref{eq:35x0} 
orbit
we have $b_{23} c_{23} = b_{13} c_{13} =
b_{12} c_{12}$.
Given the identities already established, \eqref{eq:35x0}(i)  and its images
are satisfied.
And all the remaining equations are satisfied identically.
\\
It is convenient now to write $\alp,\be$ for the eigenvalues of
$A(1,2)$; and hence (by the above identities) of all
$A(i,j)$
{$=\mat{cc} a_{ij} & b_{ij} \\ c_{ij} & d_{ij}\tam$}.
Then our conditions become $(a_i -\alp)(a_i -\be)=0$.
Finally we put $A(1,2)$
(and hence $A(i,j)$ --- all may be taken equal here)
in the `lower-1' form
$\mat{cc} \alp+\be & -\alp\be \\ 1&0\tam$ using \lemx.
\ppmx{-note Eric seems to have a few small errors in this one.}
{\qed}

\medskip


\ignore{{
\begin{lemma}
The fibre of solutions $\HH^{-1}(\upp,\upp,\upp)$ is non-empty.
The solutions $F(\sigma)$ are given up to $X$ equivalence (Lemma \ref{lem:Xeq}) by
\[
\alpha(F(\sigma)) = (a_1, a_2, a_3,
\loone, \loone, \loone
)
\]
where $a_i \in \{\alp,\be \}$.
(NB $\alp+\be\neq 0$; and if $\alp=\be$ the Jordan form is {\em not} diagonal.)
\end{lemma}
}}

\mdef \label{pa:+++}
The  
classes of cases
of type +++
may be illustrated as follows:
\ignore{{
\[
\trin{a_1}{\uppi}{a_2}{\uppi}{a_1}{\uppii}
\hspace{2.1cm}
\trin{a_1}{\uppii}{a_1}{\uppii}{a_1}{\uppii}
\ignore{{
\xymatrix{& a_2 \circ \ar@{-}[dr]^{\uppi}
\\
a_1 \circ \ar@{-}[rr]_{\uppii} \ar@{-}[ur]^{\uppi} && \circ a_1  }
\hspace{1in}
\xymatrix{& a_1 \circ \ar@{-}[dr]^{\uppii}
\\
a_1 \circ \ar@{-}[rr]_{\uppii} \ar@{-}[ur]^{\uppii} && \circ a_1  }
}}
\]
Or in the $\pfii\pai\mfii\mai$ form:
}}
\[
\bermpp{1}{2}{1}{{\pai}}{{\pai}}{{\pfii}}
\hspace{1in}
\bermpp{1}{1}{1}{{\pfii}}{{\pfii}}{{\pfii}}
\]
\[
\bermpp{1}{1}{2}{{\pfii}}{{\pai}}{{\pai}}
\hspace{.21in}
\raisebox{-.31in}{$\stackrel{\PPP(13)}{\leftrightarrow}$}
\hspace{.2in}
\bermpp{1}{2}{2}{{\pai}}{{\pfii}}{{\pai}}
\]
(note that the symmetry `breaks' our indeterminate labelling convention,
in an unimportant way).
As one will see in the proof, in the
$\pai\pai\pfii$  
cases the
$\pfii$
can be placed in
any position (unlike some configuration aspects, where the long
edge - or the orientation - induces different behaviour).
See fig.\ref{fig:sga1} for the symmetry group actions in one case.
\ignore{{
; or some here:
\[
\bermpp{1}{2}{1}{{\pai}}{{\pai}}{{\pfii}}
\]
\[
\bermpp{1}{2}{1}{{\mai}}{{\pfii}}{{\pai}}
\hspace{.3in}
\bermpp{1}{2}{1}{{\pfii}}{{\mai}}{{\pai}}
\]
}}


\module{figSymGpAct2}

\begin{figure}
\[
\xy        
\xygraph{  
!M{\mbermupp{\alpha}{\beta}{\alpha}{{\pai}}{{\pai}}{{\pfii}}{}{}{} }
[dl]  
:@{-->} [dl]^{(12)}
[urrrrr] :@{-->} [dr]^{(23)}
 [dlllllll]
 !M{\mbermupp{\beta}{\alpha}{\alpha}{{\mai}}{{\pfii}}{{\pai}}{\sim}{\sim}{} }
[rrrrrr]
!M{\mbermupp{\alpha}{\alpha}{\beta}{{\pfii}}{{\mai}}{{\pai}}{\sim}{}{\sim} }
 }
 \xygraph{
 [ddlllllll]
 :@{-->} [d]^{(23)} [urrrrrrrrr] :@{-->} [d]^{(12)}
[dllllllllll]
!M{\mbermupp{\beta}{\alpha}{\alpha}{{\pai}}{{\mfii}}{{\mai}}{\sim}{\sim}{} }
[rrrrrr]  
!M{\mbermupp{\alpha}{\alpha}{\beta}{{\mfii}}{{\pai}}{{\mai}}{\sim}{}{\sim} }
 }
 \xygraph{
 [ddlllll]
 :@{-->} [dr]^{(12)} [urrrr] :@{-->} [dl]^{(23)}
[lll]  
!M{\mbermupp{\alpha}{\beta}{\alpha}{{\mai}}{{\mai}}{{\mfii}}{}{}{} }
}
 \endxy
\]
\caption{
Action of 
symmetry group $\Sym_3$ on an $\pai\pfii\pai$ solution.
Note that some edge matrices are `twisted equal' (equal up to bar) as shown.
\ignore{{
Note that while the  $\pai\pfii$ labelling has full `equilateral'
symmetry (we can start with the $\pfii$ anywhere),
the twists do not: it is not possible to have both `twisted equalities'
to the long edge.
}}
\label{fig:sga1}}
\end{figure}
\ignore{{
and
\[
\xy
\xygraph{
!M{\mbermpp{1}{2}{1}{{\pai}}{{\pai}}{{\pfii}} }
[dd] :@/^/ [dl]^{\sigma_1}
[urrr]  :@/^/ [dr]_{\sigma_2}
[dlllll]
!M{ \mbermpp{1}{2}{1}{{\mai}}{{\pfii}}{{\pai}} }  
 [rrrr]
 !M{     \mbermpp{1}{2}{1}{{\pai}}{{\pai}}{{\pfii}} }  
[dllll] [u]  :@ [d]^{s}
\\
[dlllllll] 
!M{ \mbermpp{1}{2}{1}{{\mai}}{{\pfii}}{{\pai}} } & 
 [rrr]   !M{     \mbermpp{1}{2}{1}{{\pai}}{{\pai}}{{\pfii}} }
 \\
[dlllllll]
!M{     \mbermpp{1}{2}{1}{{\pai}}{{\pai}}{{\pfii}} }
  }
\endxy
\]
}}
\ppmx{NEXT TASK: Write a version of this in the $\II$ case.}

\begin{lemma}
The fibre of solutions
$\fHH^{-1} (\upp\doo\upp )$ is empty. 
\end{lemma}
\proof
Here $d_{12}=d_{23}=a_{13}=0$ and
$a_{12},a_{23},d_{13},b_{ij},c_{ij} \neq 0$.
But then 
\eqref{eq:34xx0}
has no solution. \qed


Given a matrix in $\Match^N(2,2)$,
let
$\AAAA{i,j}=\mat{cccc} a_i \\ &a_{ij}&b_{ij} \\ &c_{ij}&d_{ij} \\ &&&a_{j}\tam$,
the full $ij$-submatrix. 
So  for us $\fHH(\AAAA{i,j}) \in \HH$.

\begin{lemma} \label{lem:0XX}
Let $F(\sigma)$ be a solution. If $\fHH(\AAAA{1,2}) = 0$ then
$\AAAA{1,3} \stackx \AAAA{2,3}$.
\end{lemma}
\proof
By assumption $\AAAA{1,2}
= \lambda\oplus\lambda\oplus\lambda\oplus\lambda$,
i.e. $b_{12} = c_{12} = 0$ and $a_{12} = d_{12} = a_2 = a_1$.
Plugging $b_{12}=0$ into \eqref{eq:35x0} 
and its orbit
we have that
\[
a_{13} d_{23} (a_{13} -d_{23}) = a_{23} d_{13} (d_{13}-a_{23}) =
    a_1  (b_{23} c_{23} - b_{13} c_{13})
\]
thus
either $a_{13} d_{13} a_{23} d_{23} \neq 0$ --- i.e.
$\fHH(F(\sigma)) = (0,0,0)$ and we are done
--- or $b_{23} c_{23} = b_{13} c_{13} \neq 0$.
In the latter case by \eqref{eq:34xx0}
and orbit
we have $d_{13} = d_{23}$ and $a_{13} = a_{23}$ as claimed.
\ppmx{...unpack calcs a bit.}
\qed

\begin{lemma} \label{lem:sundry0s}
(I) The fibre of solutions
$\fHH^{-1} (0\off\off )$ is non-empty. The solution is given up to
{X-equivalence}
by
\[
\alpha(F(\sigma)) =
  (a_1,a_1,a_3,\mat{cc} a_1 &0\\ 0&a_1\tam, \mat{cc} 0 &b\\ 1&0\tam,
  \mat{cc} 0 &b\\ 1&0\tam)
  \]
  where $a_1,a_3,b$ take any non-zero values. \\
(II) $\fHH^{-1} (0\doo\upp )$ is empty. \\
(III.x)   $\fII^{-1} (0\pai\pfii )$ and $\fII^{-1} (0\pfii\pai )$
\ppmx{[-replace with] $\fII^{-1} (0\pfii\pai )$ ??} are empty.
\\
(III.i)   $\fII^{-1} (0\pai\pai )$  is non-empty. The solution is given
(in lower-1 form) by
\newcommand{\aaone}{\alpha}
\newcommand{\aathree}{\beta}
\[
\alpha(F(\sigma)) =
  (\aaone,\aaone,\aathree,\mat{cc} \aaone &0\\ 0&\aaone \tam,
  \loone,\loone
  )
  \]
(NB $\aaone,\aathree\neq 0$ for invertibility; $\aaone+\aathree\neq 0$
for case $\uppi$).
\\
(III.ii)  $\fII^{-1} (0\pfii\pfii )$ is non-empty. The solution is
given (in lower-1 form) by
\newcommand{\aaaone}{\beta}
\newcommand{\aaathree}{\beta}
\[
\alpha(F(\sigma)) =
  (\aaaone,\aaaone,\aaathree,\mat{cc} \aaaone &0\\ 0&\aaaone \tam,
  \loone,\loone
  )
  \]
(IV)  $\fHH^{-1} (00\doo )$ is empty. \\
(V)  $\fHH^{-1} (000 )$ contains only the `constant' solutions. 
\end{lemma}
\proof
(I) Here $b_{12}=c_{12}=a_{13}=a_{23}=d_{13}=d_{23}=0$, while
$a_{12},d_{12}\neq 0$, $b_{13}c_{13}\neq 0$ and $b_{23} c_{23} \neq 0$. 
Thus from \eqref{eq:31xx0}(i,ii) we have $a_1 = a_{12} = a_2$.
And from \eqref{eq:31xx0} (or \eqref{eq:31xxx0}(ii)) then also $d_{12} = a_1$. 
And \eqref{eq:31xxx0}(i) is directly satisfied.
\\
The 13 and 23 versions of \eqref{eq:31xx0}-\eqref{eq:31xxx0} are all directly satisfied.
\\
Also \eqref{eq:34xx0} and its images are directly satisfied.
\\
From \eqref{eq:35x0}(i) we have $b_{23} c_{23} = b_{13} c_{13}$.
From this \eqref{eq:35x0}(i) is then satisfied.
Meanwhile \eqref{eq:35x0}(ii) and its images are directly satisfied.
Note that the eigenvalues and type of $A(1,3)$ and $A(2,3)$ are the
same.
By \lemx\ we can put them in lower-1 form and claim (I) is established.\\
(II) Here  $b_{12}=c_{12}=a_{13}=d_{23}=0$, while
$a_{12},d_{12},d_{13},a_{23}\neq 0$, $b_{13}c_{13}\neq 0$ and
$b_{23} c_{23} \neq 0$.
Thus the corresponding image of \eqref{eq:34xx0}(ii) has no solution.
\\
(III)
{The $x$ cases are covered by Lemma~\ref{lem:0XX}.}
{For $i$ and $ii$ consider the following.}
Here  $b_{12}=c_{12}=d_{13}=d_{23}=0$, while
$a_{12},d_{12},a_{13},a_{23}\neq 0$, $b_{13}c_{13}\neq 0$ and
$b_{23} c_{23} \neq 0$.
Thus from \eqref{eq:31xx0}(i,ii) we have $a_1 = a_{12} = a_2$.
And from \eqref{eq:31xx0} (or \eqref{eq:31xxx0}(ii)) then also $d_{12} = a_1$. 
And \eqref{eq:31xxx0}(i) is directly satisfied.
\\
The 13 and 23 images of \eqref{eq:31xx0} are all then satisfied.
\\
Next \eqref{eq:34xx0} is directly satisfied.
Meanwhile an image of \eqref{eq:34xx0}(i) gives $a_{23} = a_{13}$;
and then the remaining images of  \eqref{eq:34xx0} are directly satisfied.
\\
From \eqref{eq:35x0} 
we get that $b_{23} c_{23} = b_{13} c_{13}$.
From
one of the $\Sym_N$ images of 
\eqref{eq:35x0} 
we get $-b_{13}c_{13}
-d_{12}a_{13}+d_{12}^2=0$.
The remaining images are all directly satisfied.\\
(IV,V) These are  corollaries of Lemma~\ref{lem:0XX}.
\qed

\ignore{{
\begin{lemma}
The fibre of solutions $\HH^{-1} (000 )$ contains only the `constant' solutions. 
\end{lemma}
\proof ...
}}

\medskip \vspace{.5cm}
This concludes the proof of Proposition~\ref{pr:N3}. \qed


\vspace{1cm}



\module{repsNproof} 

\section{Proof of 
Theorem~\ref{th:mainx}: conclusion} \label{ss:proofend}
\subsection{A key combinatorial lemma} $\;$ 

\ignore{{  
\newcommand{\chippy}{parted}
\newcommand{\choppy}{choppy}
\newcommand{\roo}{Rule-of-1}
\newcommand{\configuration}{$\II$-configuration}
\newcommand{\SET}{{\mathsf{ Set}}}  
\newcommand{\Hom}{Hom}
\newcommand{\Power}{{\mathcal P}}   
\newcommand{\pig}{\overline{\pi}}
\newcommand{\bb}{{\mathsf b}} 
\newcommand{\rrrr}{{\mathsf r}}
}}


\newcommand{\Homs}{\Setcat}

\mdef
Recall
from (\ref{de:KN})
that 
$K_n$ denotes the complete graph on vertices $\{ 1,2,...,n \}$.
Note
{again}
the natural orientation of edges $i \rightarrow j$ when $i <j$.
We will have in mind also the corresponding 2-graph, the complex over
$K_n$ that includes the triangular faces.
We write $K_n^0$, $K_n^1$, $K_n^2$ for the sets of vertices, edges and faces respectively.
Thus $K_n^0 = \{ 1,2,...,n \}$;
%
and
$\Homs(K_n^1, \{x,y\})$ is the set of colourings of the edges of $K_n$
from the colour set $\{x,y\}$ (here $x,y$ are arbitrary symbols/`colours').

\mdef
For a set $A$, we write $\Power(A)$ for the power set. 
\ignore{{
Observe that  $\Hom(K_n^1, \{x,y\})$ is in bijection with $\Power(K_n^1)$.
(A specific bijection
$$
\bb_y : \Hom(K_n^1, \{x,y\}) \stackrel{\sim}{\rightarrow} \Power(K_n^1)
$$
arises by choosing $x=0$, $y=1$ say.)
}}
Recall that the map
$$
\bb_y : \Homs(A, \{x,y\}) \stackrel{\sim}{\rightarrow} \Power(A)
$$
$\bb_y(f) = f^{-1}(y)$,  
is a bijection.

In turn the set $\Power(K_n^1)$ can be considered as the set of full subgraphs
(i.e. subgraphs retaining all vertices) of $K_n$.
Furthermore each subgraph $\gamma$ defines a partition
$\pig(\gamma) \in \PP(K_n^0)$ 
of the vertex set into
connected components.
That is, we have a map
$$
\pig: \Power(K_n^1) \rightarrow \PP(K_n^0)
$$

Note that an element $c$ of $\Homs(K_n^1, \{x,y\})$ gives an element
$\rrrr_y(c)$
of the set $\Power(K_n^0 \times K_n^0)$ of relations on $K_n^0$ by 
$(i,j) \in \rrrr_y(c)$ if $c(i,j)=y$.
The partition
$\pig_x(c) := \pig(\bb_y(c))$ is simply the
partition corresponding to the reflexive-symmetric-transitive closure
of $\rrrr_y(c)$.


\newcommand{\yy}{y} 
\newcommand{\xx}{x} 
\newcommand{\Homx}{\Hom_x}

\begin{lemma} \label{le:chip} 
(\roo\ Lemma)
An edge 2-colouring of $K_n$,
with colour set $\{ \xx,\yy \}$,
i.e. an element of $\Homs(K_n^1, \{ \xx,\yy \} )$,
is called \xx-{\em \chippy}
if every triangle that has two \yy's has three; i.e.
if no triangle has a single \xx.
We write $\Homx(K_n^1, \{ \xx,\yy \} )$ for the set of  \xx-\chippy\ colourings.
\\
(I) The set of \xx-\chippy\ colourings is in bijection with the set
$\PP(K_n^0)$ of partitions of the vertices.
In particular
$\pig_x = \pig \circ \bb_y$ restricts to a bijection
$\pig_x : \Homx(K_n^1, \{ \xx,\yy \} ) \stackrel{\sim}{\rightarrow} \PP(K_n^0)$.
\\
(II) The subset of colourings where no triangle has three \xx's is in
bijection with the subset of partitions into two parts.
\end{lemma}

\proof
(I) This  
amounts to 
a rearrangement of the definition of equivalence relation,
where \yy\ on  edge $\{ v,v' \}$ means that $v \sim v'$.
\\
(II) This gives the subset where no triple has three vertices in
different parts.
\qed

An \xx-\chippy\ configuration can be
realised  
geometrically: think of $\yy$s as being asigned to edges much shorter that $\xx$s.
If $ij$ and $jk$ are short then $ik$ is short by the triangle inequality.
Conversely if $ij$ is long then at least one of $jk$, $ik$  
is long.
See Fig.\ref{fig:puppyears} {at (\ref{de:SSSN}).}
\\ (Caveat: we will use the 
lemma in several different ways.  E.g. with $x=\off$ and then $y=0$.)

\ignore{{

\mdef Remark.
Let $\lambda$ be the integer partition corresponding to a given set
partition of a set of order $n$.
Then note that the number of \xx's in the corresponding equivalence relation,
encoded as an \xx-\yy-colouring of $K_n$ as in Lemma~\ref{le:chip},
is
\beq
h(\lambda) = \frac{n^2 - \sum_i \lambda_i^2 }{2} 
\eq
For example $h(31) = (16-9-1)/2 =3$ while $h(22)=(16-4-4)/2=4$.

}}

\subsection{Admissible configurations} $\;$ 

\mdef \label{de:Higgs}
Recall that
a {\em (Higgs) \configuration} of $K_n$ is an assignment of a scalar to
each vertex;
and matrix to each edge,
so that each $K_2$ subgraph is assigned
from the set of six types
$\II = \{ \off, 0, \uppfii, \mfii, \uppai,\mai \} $
as in \eqref{pa:schema2}.
(Note that there are variables in all these components.
And that we can work up to X-equivalence, as in Lemma~\ref{lem:Xeq}.)

By Proposition~\ref{pr:L2class2}
and Theorem~\ref{th:constraints}
a necessary condition for an element of $\Match^N(2,2)$
to give a braid representation is that it gives an \configuration\ of $K_N$.

\newcommand{\groundstate}{groundstate}

\mdef \label{de:admissible}
A `\groundstate'
or `admissible'
configuration of $K_n$ is an \configuration,
as in (\ref{de:Higgs}),
such that every
triangle configuration is
in the orbit of 
one of those 
ten
forms given in Proposition~\ref{pr:N3}.
(Note that variables are constrained by these conditions, as well as types.)

By Proposition~\ref{pr:N3}
and the form of the constraints in Theorem~\ref{th:constraints}
(from which we observe that the total constraint set is the union of constraint
sets across all $K_3$ subgraphs)
a necessary and sufficient condition for an element of $\Match^N(2,2)$
to give a braid representation is that it gives an admissible configuration of $K_N$.

\subsection{The proof for general $N$} $\;$ 



\newcommand{\bI}{{\bf I}}
\newcommand{\bII}{{\bf II}}

We require to 
show that $\Rec(\SSS_N^\C)$ gives the set of 
braid representations;
and by \eqref{de:admissible} this is the same as the set of \groundstate s.
That is, 
we require to show
\\
(\bI) that every solution is in  $\Rec(\SSS_N^\C)$ {up to X-equivalence},
and
\\
(\bII) that every element of  $\Rec(\SSS_N^\C)$ gives a solution.

For (\bI)
We proceed as follows. Suppose $\alphaa$ is a solution, i.e. an
admissible configuration.
We require to show that there is a $(p,q,\rho,s) \in \SSS_N$
and a point in $\C^{(p,q)}$ 
such that applying $\Rec$
to this germ
gives $ \alphaa$. 
First let us consider the data of the $\off$ edges.
For this (and in various guises hereafter) we will make use of
Lemma~\ref{le:chip}.




\newcommand{\partp}{\pi_{\mathsf{p}}}
\newcommand{\partq}{\pi_{\mathsf{q}}} 
\newcommand{\partr}{\pi_{\rho}} 
\newcommand{\parts}{\pi_s} 

\begin{lemma} \label{le:off1}
Let $\alphaa$ be  an admissible configuration of $K_N$.
\\
(I)
The $\off$ data of
$\alphaa$ 
induces a partition
$\partp(\alphaa)$ 
of the vertices (into `nations');
and in particular there is a $\off$ on an edge if and only if this edge is between
vertices in different parts.
\\ (II)
Furthermore, every $\off$ between the same two nations carries the same parameter
{$\mu_{ij} = \sqrt{b_{ij} c_{ij}}$}.
\\ (III)
Thus a solution $\alphaa$ is partially characterised by
the partition $\partp(\alphaa)$ and the parameters it associates to the
edges of $K_{\partp(\alphaa)}$
(each vertex of  $K_{\partp(\alphaa)}$ represents a nation; and 
each edge of  $K_{\partp(\alphaa)}$ represents the collection of edges of
$K_N$ that are between the same two nations).
\\ (IV) The $\off$ data of $\alphaa$ agrees with what we obtain by applying
$\Rec$ to a germ with $p=\partp(\alphaa)$.
\end{lemma}

\proof
(I)
Observe from the list of allowed triangle configurations
in Prop.\ref{pr:N3}
that this
data is \chippy\ with $\off$ in the role of $\xx$ 
and all other decorations together in the role of $\yy$.
Now apply the \roo\ Lemma, \ref{le:chip}(I).
\\ (II)
The parameter constraint comes from observing the configurations in
Prop.\ref{pr:N3} with exactly two $\off$s.
We see firstly that if two edges between nations meet then they carry the same
parameter. So now suppose
$\{v,v' \} $ are 
distinct vertices in one nation, and
$\{w,w'\}$ in another.
Thus $vw$ and $v'w'$ are edges between nations, hence $\off$ edges.
But of course so is $vw'$, and this has the same parameter as both.
{Note that $\mu_{ij}=\sqrt{b_{ij} c_{ij}}$ is X-invariant, thus the X-orbit of $\alphaa$ contains
a rep with these $b_{ij} =c_{ij}$ as in $\Rec$.}
\\
(III,IV) follow immediately. Note that there are 
points in $\C^{(p,q)}$ with  $\mu$ parameters as in $\alphaa$. 
\qed

\begin{lemma} \label{le:0s}
Let $\alphaa$ be  an admissible configuration of $K_N$.
\\
(I)
There is a partition $\partq(\alphaa)$ of the vertices (into `counties')
such that
\ignore{{
Consider a part of the partition $\partp(\alphaa)$ induced by $\off$. 
Edge 0's induce an equivalence 
on vertices
(i.e. a partition, of this nation into `counties') ---  in particular
}}
two vertices are in the same part if and only if their common edge has a 0.
The partition $\partq(\alphaa)$ refines $\partp(\alphaa)$.
\ignore{{
Thus overall we have a partition $\partq(\alphaa)$ refining $\partp(\alphaa)$.
And then 
`signs' (the bar-unbar data in the edges labeled from $\{ \pfii,\mfii,\pai,\mai\}$)
induce a  
total order on equivalence
classes of the 0-equivalence relation 
(the counties in a nation).
\\
In particular
}}
And then
there is a total order on counties in each nation so that
if county $q_j$ is
before 
$q_k$ in the order and vertex $v$ is in $q_j$
and vertex $w$ in $q_k$ then the $vw$ edge is +-type (i.e. from $\pfii,\pai$) if $v<w$
and is $-$-type (from $\mfii,\mai$) otherwise.
Thus overall we have an order $\partr(\alphaa)$. 
\\
(II) The $\off$, 0 and $\pm$ data of $\alphaa$ agrees with what we obtain by
applying $\Rec$
(from \ref{pa:rec})
to a germ with
$p=\partp(\alphaa)$, $q=\partq(\alphaa)$ and $\rho=\partr(\alphaa)$.
\end{lemma}

\proof
(I)
We see that
edge 0's induce an equivalence by comparing allowed triangle
configurations with the \roo\ Lemma, \ref{le:chip}.
Taking the county partition data across all nations we have
a refinement of
$\partp(\alphaa)$ --- 
this is
$\partq(\alphaa)$.
A configuration induces a well-defined order
$\partr(\alphaa)$
on the quotient
--- the counties in each nation ---
since
triangles are never cyclic-ordered
(we will give  illustrative examples next)
and triangles with a 0 `collapse'
consistently.

\newcommand{\<}{\prec}
\newcommand{\succc}{\succ}

For example consider a triangle with vertices 2,5,9
(say, any triple will do --- note that our edge directions on $K_N$ will give
an `ambient' order on these vertices)
and suppose first in the same nation
but not the same county --- so decorated from $\{ \pfii,\mfii,\pai,\mai\}$.
For now we consider only the $\pm$ aspect of this.
The induced relation between 2,5 is
$2 \< 5$
if + and
$2 \succc 5$
if -. 
The relation between 5,9 is $5\< 9$ if + and
$5 \succc 9$
if -.
From the \ninerule\ (in particular Figure~\ref{fig:sga1})
we see that if $2\< 5$ and $5\< 9$ then $2\< 9$,
and so on, so this is an order.
Specifically, looking at the Figure, the top case is all-+ and so we use the
ambient order, which is of course acyclic.
Subsequent cases are generated by the
$\Sym_N$ action so, for example,
if we perm 2,5 then this edge becomes - while the others
remain +, thus we have $5\< 2\< 9$.
Meanwhile if instead we perm 2,9 then all signs are changed
{(from + to -)}
and we indeed have
$9\< 5\< 2$.
On the other hand suppose there is a 0 decoration.
Again from the \ninerule, which is given explicitly for ++0 cases, we see here that
if $5\sim 9$ then $2\< 5$ and $2\< 9$;
and then the $\Sym_N$ action generating other cases
respects this consistent order.
$\;$ (II) follows immediately from \ref{pa:rec}.
\qed


\begin{lemma}
Let $\alphaa$ be  an admissible configuration
of $K_N$.
\\ (I)
Consider a part $p_i$ of the partition
$\partp(\alphaa)$
induced by $\off$;
and the partition of the vertices of this part induced by 0,
as in Lemma~\ref{le:0s}. 
The $\pai\pfii$ edge decoration data induces a partition on the parts of $p_i$
into (at most) two parts, 
with $\pai$ in the separating role;
with  two vertices from the same nation only having different parameter
if in different such parts, i.e. with edge label $\pai$.
\ignore{{(To record this data as a single two-part partitioning across all nations,
we may order the two parts in each nation by the child-first order,
then simply combine all the first parts, and combine all the second parts.)
}}
Taken over all nations we thus have 
$\parts(\alphaa) \in \prod_{i} \PP_2(q|p_i)$.
\\
{(II)}
For each nation $p_i$ either all edges are 0 and there is a single parameter;
or there is at least one edge from $\{ \pfii,\mfii,\pai,\mai\}$
and any one of these determines the two parameters
appearing in all.
\\
(III)
Applying $\Rec$ to a germ with 
$(\partp(\alphaa),\partq(\alphaa),\partr(\alphaa),\parts(\alphaa))$
we obtain $\alphaa$ {up to X-symmetry}.
\ppmx{-now explain the vertex variables - all the same if 0-connected
and same if f-connected. ...}
\end{lemma}

\proof
This is  a matter of unpacking the definitions and then using
the \roo\ Lemma, \ref{le:chip}, for a third (!) time, this time using  
\ref{le:chip}(II),
because there is no
$\pai\pai\pai$ triangle.
$\;$
{
(II) comes from observing the configurations in Prop.\ref{pr:N3} with two or more
$\pai/\pfii$ edges analogously to (\ref{le:off1})(II),
noting that parameters are X-invariant
(cf. e.g. Lem.\ref{le:+++}).
}
The final point follows immediately from the $\Rec$ recipe
in \ref{pa:rec}, and appropriate choice of parameters
(noting \ref{le:0s}(II)).
\qed

\ignore{{
\mdef
In short then, a configuration $\alphaa$ should be organised firstly as a
partition according to $\off$ (we say a partition into `nations').
Then as a partition according to 0 (`counties').
Then 
...
\ppm{NEXT explain variables}
Now observe
from Proposition~\ref{pr:N3}
that any triangle between two nations, thus with exactly two $\off$s,
has the same variable on each of these edges.
Thus every edge between these two nations has the same variable.
\\
Similarly 
...
}}


\medskip


This concludes requirement (\bI). $\;$
For
(\bII) Let $(p,q,\rho,s) \in \SSS_N$ and consider $\Rec(-)$.
We require to check 
that every triangle is as in Prop.\ref{pr:N3}.

\begin{lemma} \label{le:IIi}
\ppmx{Paul, push this to a later Lemma?, since it belongs to stage (II)
of main proof:}
Suppose we have a configuration of $K_N$ constructed  
from
a partition of vertices 
$p$ by 
decorating each edge between the same two parts with a $\off$
carrying the same parameter (as in $\Rec$).
Then all the triangles involving $\off$ are
{admissibly configured}.
\end{lemma}

\proof  
As before we write $i<j<k$ for the vertices and order the edges
$ij,ik,jk$.
If all three vertices lie in different nations then $\Rec$ gives
edge configuration $\off\off\off$, which is admissible.
If two vertices lie in the same nation then we have $\off\off *$,
$\off *\off $ or $*\off\off $, for some *, with the $\mu$ parameter the same
on the two $\off$ edges.
From inspection of cases in Prop.\ref{pr:N3} we see that all such
are admissible. 
\qed

\begin{lemma} \label{le:IIii}
Suppose we have a configuration of $K_N$ constructed firstly as in
Lemma~\ref{le:IIi}
and then from a refinement of $p$ to $q$ by decorating
each edge within a county $q_i$ by 0;
and then from a
total ordering $\rho$ of the counties in each nation
by decorating each edge between counties in the nation as in $\Rec$
from \ref{pa:rec}.
Then
every edge is decorated from $\{\off,0,+,- \}$, and 
every triangle is  consistent at the level of orientation with
an admissible configuration.
\ppmx{...hmm, maybe better to already include the $s$ data now for this bit?
- perhaps then demoting this part of Lemma to next Lemma.
Or, sticking with $+-$, maybe we can get away with it here.}
\end{lemma}

\proof 
Triangles involving $\off$ involve more than one nation, and are
already admissible by Lemma~\ref{le:IIi},
so we need only check each triangle lying within a nation.
As before we write $i<j<k$ for the vertices and order the edges
$ij,ik,jk$.
If all three vertices lie in the same county then $\Rec$ gives
edge configuration 000, which is admissible.
If two of three vertices lie in the same county then $\Rec$ gives
a configuration like 0** or *0* or **0, with the *'s depending on $\rho$.
In the 0** case $i,j$ are in the same county.
If $\rho = (ij,k)$
(in the obvious shorthand) then $k$ is above both
in both the natural and the $\rho$ order, so we have 0++.
Comparing with  Prop.\ref{pr:N3} and the orbit table we see that
this is consistent with  admissibility.
(The cases shown in  Prop.\ref{pr:N3} are both $\rho=(i,jk)$ which is ++0.)
If $\rho = (k,ij)$ then we have $0--$ which is again consistent.
Finally for these cases note that $\Rec$ assigns the same parameters to
both edges as required.
It remains to consider the cases with three different counties,
such as $\rho = (i,j,k)$. This gives +++ with all edges assigned the
same parameters, which is again consistent.
The case  $\rho = (j,i,k)$ gives $-++$ (cf. Fig.\ref{fig:sga1})
and the rest of the orbit follows.
\qed


\begin{lemma} \label{lem:Rend}
Suppose we have a configuration of $K_N$ constructed firstly as in
Lemma~\ref{le:IIii}
and then
from a partition $s$ of the vertices into two parts
by assiging vertex parameters and replacing $\pm$ labels with
corresponding $\pai,\pfii,\mai,\mfii$ labels as in $\Rec$.
Then every triangle is admissible.
\end{lemma}
\proof
With regard to the $\pai\pfii$
labels, irrespective of signs, the admissibility conditions stipulate only
that configurations like $\pai\pfii 0$ and $\pai\pai\pai$ must not arise.
The former is impossible since $\Rec$ assigns the same edge label to
every edge between a given vertex and a county.
Since $s$ partitions counties in a nation into at most two parts
$\pai\pai\pai$ is not possible. \qed


\newcommand{\arr}{\begin{array}}
\newcommand{\rra}{\end{array}}
\newcommand{\guu}[1]{{}^{#1}}
\newcommand{\typ}[2]{
\arr{c}
\includegraphics[width=.41cm]{xfig/#1.eps}
\\ \guu{\alpha_1} \\ \guu{\beta_1}
\rra \!\!
\raisebox{5pt}{$
\stackrel{\rule{15pt}{.4pt}}{{}_{\mu_{12}}}
$}
\!\!
\arr{c}
\includegraphics[width=.41cm]{xfig/#2}
\\ {\guu{\alpha_2}} \\ \guu{\beta_2}
\rra}

\newcommand{\typp}[2]{
\arr{c}
\includegraphics[width=.41cm]{xfig/#1.eps}
\\ \guu{\alpha_{#2}} \\ \guu{\beta_{#2}}
\rra 
}
\newcommand{\typpp}[2]{
\arr{c}
\includegraphics[width=.41cm]{xfig/#1.eps}
\\ \guu{\alpha_{#2}} 
\rra 
}


\medskip

{By (\ref{de:admissible}), Lemma~\ref{lem:Rend} implies that
$\Rec: \SSS_N^\C \rightarrow \functor(\Bcat,\Match^N) $,
i.e. that $F(\sigma) = \Rec(\lambda)$ is a solution for all $\lambda$.
By (\bI) we know that it is surjective on X-orbits.}

Having established directions (\bI) and (\bII), 
 this concludes the proof of Theorem~\ref{th:mainx}. \qed




\newcommand{\three}{\{1,2,3\}} 

\mdef \label{pa:proofJJN}
{\em Proof} of Theorem~\ref{th:JJN}.
{
Transversality of $\TTT_N$ in $\SSS_N$
under the {\em direct} action is clear. 
Consider a diagram  $\lambda\in \SSS_N$, and hence a variety of solutions
$\Rec(\lambda)$
(meaning the output of $\Rec$ with the variables left indeterminate,
or the set of solutions as the parameters vary).
Write $f$ for the underlying shape of $\lambda$.
Observe that  $\bb(f) \in \Sym_N \lambda$
(under the direct action)
for all $\lambda$.
Since we consider the whole variety here the precise fate of the parameters
need not be tracked
(although see for example \S\ref{ss:flyp}).
The result is now a corollary of Thm.\ref{th:mainx}.
(NB Injectivity
of $\bb\circ\jjj :J_N(\three^*)\rightarrow \SSS_N$
follows routinely from the construction of $\jjj$.)
}
\qed
\ignore{{
\ppm{
[-not sure this is clear if we work with varieties, or true if we don't!]
\\
MAKE A MINIMAL FIX OF THIS AND BE DONE! Explain a bit what the variety is -
relate type-space to variety...
\\
maybe with example?...
\\
%
First observe that the $\Sym_N$ action
acts on a diagram $\lambda\in \SSS_N$ simply by permuting entries.
It acts on the solution variety $\Rec(\lambda)$ in a more complicate way in general,
also `changing signs'
(as well as making more trivial changes to the indices on variable names).
But by Lem.\ref{lem:no-} there is a \nominus\
representative.
From $\Rec$ we see that here if $i$ is in a county above $j$ in the same nation
then $i<j$.
Now noting the consecutive property from Lem.\ref{lem:no-} we see that our
representative will obey at least the `page' part of book order.
Next however, as already noted, nations can be permuted freely changing only the
indices on $\mu$ variables... [-this is not enough - somewhere we have a para that
says what we want here??]
\\
since the \nominus\ ...
FINISH ME! \qed
}
}}





\newcommand{\I}{\mathcal{I}}
\newcommand{\OO}{\mathcal{O}}
\newcommand{\V}{\mathcal{V}}
\newcommand{\om}{\omega}
\newcommand{\lan}{\langle}
\newcommand{\ra}{\rangle}
\newcommand{\End}{{\rm End}}
\newcommand{\ve}{\varepsilon}
\newcommand{\T}{\mathcal{T}}

\newcommand{\Dim}{{\rm Dim}}
\newcommand{\sspan}{{\rm span}}
\newcommand{\ga}{\gamma}
\newcommand{\Ga}{\Gamma}
\newcommand{\one}{\mathbf I}
\newcommand{\ka}{\kappa}
\newcommand{\Z}{\mathbb{Z}}
\newcommand{\g}{\mathfrak{g}}
\newcommand{\Vn}{V^{\otimes n}}

\newcommand{\la}{{\lambda}}
\newcommand{\La}{{\Lambda}}
\newcommand{\al}{\alpha}

\newcommand{\M}{\mathscr{M}}
\newcommand{\D}{\mathscr{D}}

\renewcommand{\Up}{+}
\newcommand{\Do}{-}

\ignore{{
\theoremstyle{plain}
\newtheorem{theorem}{Theorem}
\newtheorem{lemma}{Lemma}
\newtheorem*{prop}{Proposition}
\newtheorem*{cor}{Corollary}
\theoremstyle{definition}
\newtheorem*{definition}{Definition}
\newtheorem*{heur}{Heuristic}
\theoremstyle{remark}
\newtheorem*{prob}{Problem}
\newtheorem*{rmk}{Remark}
\newtheorem*{ex}{Example}
\newtheorem*{an}{Analogy}
}}


\module{future}

\section{Future Directions} \label{ss:future}

Several natural questions for future study spring to mind.
\begin{enumerate}
\item 
The Andruskiewitsch-Schneider classification program for pointed Hopf algebras takes as input a solution $(R,V)$ to the YBE to produce Nichols algebras which can potentially be lifted via bosonization. The $4\times 4$ solutions were addressed \cite{AndruskGiraldi2018}. Do our solutions lead to finite-dimensional Nichols algebras?
\item Turaev \cite{Turaev1988} developed a method for constructing link invariants from solutions to the Yang-Baxter equation. Of course the well-known $4\times 4$ $R$-matrix associated with the $U_q\mathfrak{sl}_2$ provide link invariants (essentially the Jones polynomial), and many of our solutions generalize these solutions so we certainly expect interesting invariants.
\item
It is noteworthy that the verification of a $N^2\times N^2$ solution comes down to verifying  constraints on $2$-indexed scalars for all $3$-element subsets.  This bears some similarity to the verification of associativity in a monoidal category: it is sufficient by Mac Lane coherence \cite{MacLane} to verify constraints involving $4$-indices on scalars with 3 indices.  This suggests a higher category connection, see \cite{KapVoev}.
\item
Also natural is the question of generalisations to other targets besides $\Match^N$,
for example generalised 8-vertex model
(in the sense of \cite{Hietarinta:1992}) instead of the 6-vertex model.
{
It is indeed interesting to 
ask what is the measure of our set of solutions in the space of all
solutions up to isomorphism.
This could in principle be addressed for example via a strategy as in \cite{Martin92}.}
\item
  Many interesting questions arise from the perspective of representation theory
  itself. We touch on this in \S\ref{ss:repthy}. 
\end{enumerate}

\medskip
\appendix  
\section{Applications and representation theory} \label{ss:repthy} 
\module{braid-transversalP}

\ignore{
In \S\ref{ss:flyp} we use the $\PPP$ symmetry and some examples to set up
our result as a paradigm in higher/extended representation theory,
from which several further intriguing questions arise.
}


One intended
application for this work is to provide paradigms in  
representation theory.
For example we can consider the
source and 
target dependent
notion of equivalence of representations.
A canonical treatment should not be expected, but some 
boundaries of the problem can be usefully delineated.
In the one extreme we have
manifest symmetries which certainly imply equivalence, but evidence to suggest,
as in \S\ref{ss:flyp}, that these are not exhaustive;
and in the other we have
computable relations between
representations 
that may not imply equivalence,
but whose {\em absence} implies non-equivalence,
as in \S\ref{ss:signature}.
%


\medskip

We have seen  
in \S\ref{ss:comb0}
that the $\Sym_N$-transversal $\TTT_N$ of $\SSS_N$
(indexing varieties of solutions)
has a beautiful and well-behaved combinatoric
- i.e. the Euler transform of the sequence $3^{N-1}$,
see also \S\ref{ss:signature}.
In this sense
the $\Sym_N$ symmetry 
is  ideal for
the purpose of classifying and constructing solutions.
But we know from (\ref{de:flip1}) that this $\Sym_N$
does not exhaust the available symmetries,
and hence not the possible
realisations or
notions of equivalence.
%

A broad notion of equivalence of representations
for single-generator structures is matrix similarity of generator image.
Not 
every similar matrix to a representation gives a representation
(this is special to ordinary representation theory),
but {\em if} two representations are
{\em not} similar then we say they are not equivalent.
On the other hand, while it is computationally relatively easy to determine if two
matrices are similar
(i.e. have equal Jordan form)
case by case, there is not yet a general theory here.
In our $\Match^N$ case, though, Jordan forms are relatively easy to compute
in practice.
Simpler still
is to address the necessary but not sufficient condition that the
spectrum is the same. For this see \S\ref{ss:spectrum}.

\newcommand{\locy}{locally similar}
\newcommand{\locyp}{locally permutation-similar}

From this perspective, one can say that representations are {\em \locy}
if they are similar by conjugation by a matrix of form $Q \otimes Q$
(in $\Mat^N$ any representation and any $Q$ gives a representation, but in
$\Match^N$ this is far from true). 
The $\Sym_N$ action can be seen as
locally {\em permutation}-similar.

Another possibility is to consider
the transversal of $\SSS_N$
with respect to  
the $\Zz_2$-orbits under the $\PPP$ symmetry
(in {\em tandem} with the $\Sigma_N$ symmetry).
We do this next.


\subsection{On the flip action on solution space  
} $\;$ \label{ss:flyp}

\begin{figure}
\[
\includegraphics[width=6.8\ccm]{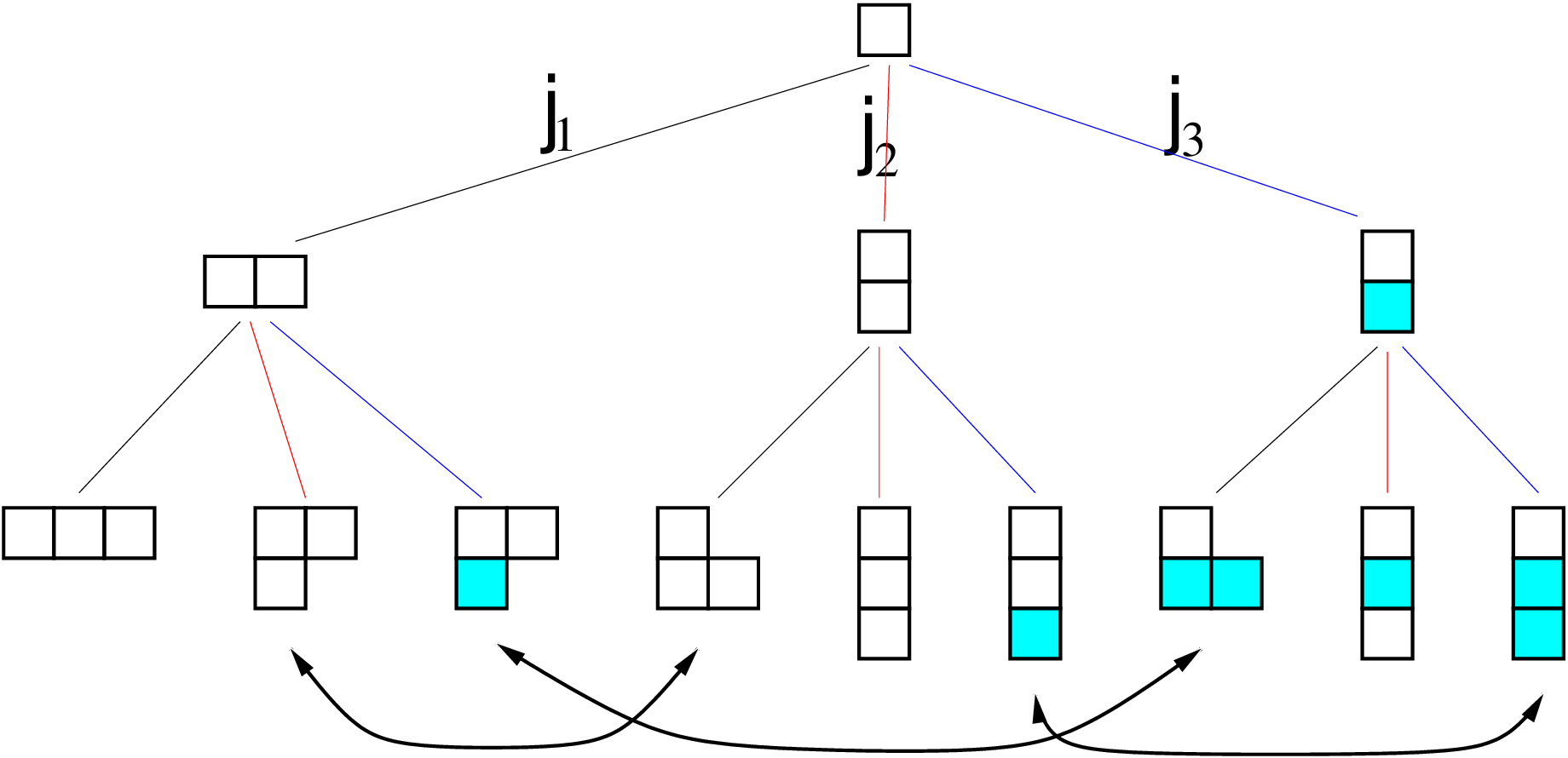}
\hspace{3cm}
\includegraphics[width=.8\ccm]{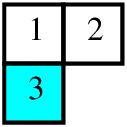}
\;\raisebox{.1in}{$\stackrel{\PPP}{\leftrightarrow}$}\;
\includegraphics[width=.8\ccm]{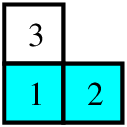}
\;\raisebox{.1in}{$\stackrel{(13)}{\leftrightarrow}$}\;
\includegraphics[width=.8\ccm]{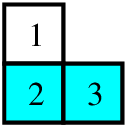}
\]
\caption{ All one-nation solution types up to rank-3,
and 
all non-trivial $\PPP$-orbits.
\label{fig:tree123}}
\end{figure}


\mdef \label{pa:Paction}
On a given solution $R \in \Match^N(2,2)$
the  
$\PPP$ action leaves the $a_i$s unchanged
but  
interchanges the diagonal and skew-diagonal entries on all $A(i,j)$ simultaneously.
On a germ this has the effect of simultaneously changing all  $\pfii$ edges
to $\underline{\pfii}$,  $\pai$ edges to $\underline{\pai}$ and vice versa.
It leaves the $/$ edges and the $0$ edges unchanged, so that the
\emph{memberships} of nations
{and counties}
are unchanged; 
but the order of counties in their nation is reversed.

For example consider the variety of solutions associated to the first diagram here:
\[
\includegraphics[width=.8\ccm]{xfig/tree123Pa13b.eps}
\;\raisebox{.1in}{$\stackrel{\PPP}{\leftrightarrow}$}\;
\includegraphics[width=.8\ccm]{xfig/tree123Pa13bx.eps}
\;\raisebox{.1in}{$\stackrel{(13)}{\leftrightarrow}$}\;
\includegraphics[width=.8\ccm]{xfig/tree123Pa31.eps}
\]
\ignore{{
Whether this new solution
\ppm{[-do you want to talk about individual solutions, or varieties?
probably better varieties?]}
is distinct from the original solution depends on the configuration data.
In general then,
the $\Zz_2$ action breaks $\SSS_N$ into ``palindrome'' singletons and pairs.
}}
%
%
From this diagram,
for each point $(\alpha_1, \beta_1)$ in type space
(a specific pair of complex numbers)
we 
get the solution with
$a_1 = a_2 = \alpha_1$, $a_3 = \beta_1$ and so on.
Applying $\PPP$ does not change the $a_i$s,
so we have the middle diagram and still $(\alpha_1, \beta_1)$.
Applying $(13)\in\Sym_3$ takes the $a_i$s to $a_1 = \beta_1$
and $a_2  = a_3 = \alpha_1$.
This
is the specific solution we would associate  to
$( \beta_1, \alpha_1)$ not
$(\alpha_1, \beta_1)$ via the last diagram, but these are points in the same variety. 

We see that starting
in $\TTT_N$ the action of $\PPP$  takes us out of the subset.
But
since $\PPP$ and $\Sym_N$ actions commute
we have an action of $\PPP$ on classes (represented by unfilled diagrams).
Up to rank 3
the orbits are  
indicated in Fig.\ref{fig:tree123}.



Consider the varieties of solutions indicated by the following
rank-4 one-nation diagrams:
\[
\includegraphics[width=.8\ccm]{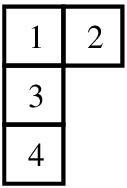} \;
\includegraphics[width=.8\ccm]{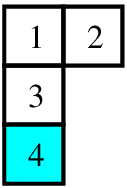} \;
\includegraphics[width=.8\ccm]{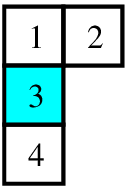} \;
\includegraphics[width=.8\ccm]{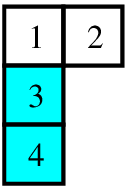}
\;\;\;\;\;
\includegraphics[width=.8\ccm]{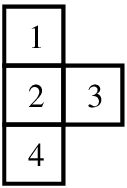} \;
\includegraphics[width=.8\ccm]{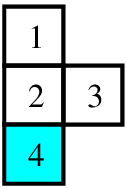} \;
\includegraphics[width=.8\ccm]{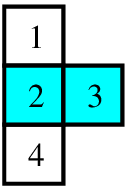} \;
\includegraphics[width=.8\ccm]{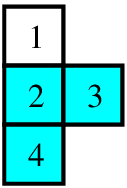}
\]
The
first is  $\bb(\jjj(122))$,
and the fifth for example is $\bb(\jjj(212))$.
The fifth (also the seventh) is fixed under
our
extended transformation group $\GG_N = \Z_2\times\Sym_N$.
However
the fifth 
gives solutions that are conjugate, as  matrices, to those from the
first.  
Conjugation here can be achieved by a
{\em non-local} permutation matrix in $\Match^N(2,2)$.
\ignore{{
However
the fifth (and seventh) is fixed under
our
extended transformation group $\GG_N = \Z_2\times\Sym_N$.
}}
\ignore{{
However the
first,  $\bb(\jjj(122))$
gives solutions that are
 conjugate, as  matrices, to those from the fifth,
$\bb(\jjj(212))$.
}}

\ignore{{

\noindent
\ppm{[DELETE the rest of this??...]}

Consider our example with $N=13$ from (\ref{eq:eg13}):
\\
\includegraphics[width=6.8\ccm]{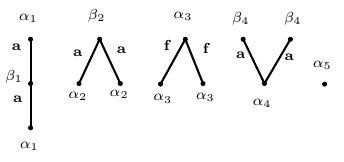}
%
\ppm{[...Not sure I understand what is going on here exactly.
(I know it develops an example from  (\ref{eq:eg13}) 
but I realise that I did not understand that either...)]}
The effect of the $\Zz_2$ action is to replace the edge decorations $\textbf{x}\rightarrow\underline{\textbf{x}}$.  
We first apply a permutation to change these back.  
Explicitly we should apply the permutations (in cycle notation) $(1\;3)$, $(4\; 5\; 6)$, $(7\; 8\; 9)$ and
$(10\; 12)$).
Next we apply a permutation to ensure the nations respect our monomial ordering, namely we interchange the 2nd and 4th nations.
The resulting diagram is on the left below, with the variables renamed on the right (with our top-to-bottom, left-to-right indexing convention):
\[
\includegraphics[width=7\ccm]{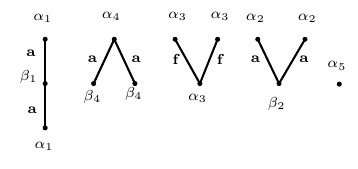}
\;\;\;\;\;\raisebox{1.8cm}{$\stackrel{}{\rightsquigarrow}$}\;\;\;\;\;
\includegraphics[width=7\ccm]{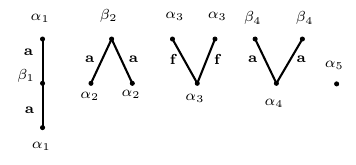}
\]
 The $\mu_{ij}$ with $i$ or $j$ in $\{2,4\}$ that we have suppressed would also be out of order due to the permutation.

It is worth pointing out a few subtleties. First it should be clear that the
$\PPP$ symmetry is distinct from that of $\Sigma_N$.
Second, the $\pfii$ and $\pai$ choices, i.e., the element $s$ responsible for the vertex labels, can break symmetry.  For example the sub-solution corresponding to the third and fourth nations in our example (abbreviated $\pfii\wedge,\pai\vee$)  has orbit size $2$ when resolved to elements of $\TTT_N$ (the other being $\pai\wedge,\pfii\vee$), whereas the sub-solution corresponding to the second and fourth nations (i.e., $\pai\wedge,\pai\vee$) has orbit size $1$ at this level.


}}




\module{reps-N4signature}

\newcommand{\maps}{\rightarrow}  

\subsection{Signatures:
Eigenvalues and isomorphism from a linear perspective} \label{ss:signature} \label{ss:spectrum}   $\;$ 



\newcommand{\varalpha}{\alphaa}
\newcommand{\sep}{:}  

A {\em signature}  is, roughly, the eigenvalue spectrum of the solution matrix
for each solution type.
Since the actual eigenvalues of a solution $\FF(\sigma)$ come from
the specific point on the variety,
what we focus on here is the generic
degeneracies
rather than the values.
Since the largest block matrix in $F(\sigma)$ is $2\times 2$, we can
compute the generic degeneracies for each variety
(indeed in principle we can compute a  Jordan form including all points).
We will give the degeneracies as an unordered list, and hence as a partition of $N^2$.
Notationally we
may
write
$\varalpha \maps \omega(\varalpha)$ if
$\omega(\varalpha)$ is the list of degeneracies for
configuration $\varalpha$.


It is straightforward to work out eigenvalues and multiplicities
for a specific class of solutions --- i.e. for $F(\sigma)$.
As noted,
if two eigenvalue spectrums do not agree, then the solutions
cannot be isomorphic. But even when they do agree, they may not be
isomorphic.

Rank-1 is trivial: we write $\square\maps (1)$ to denote that the trivial
variety $\square$ has solutions with a single eigenvalue with degeneracy 1.
Rank-2 is straightforward.
Our types $0,\pfii,\pai,\off$ are:
\\
$\square\!\square\maps (4)$, \quad
$\includegraphics[width=.251cm]{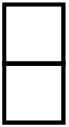} 
\maps (3,1)$, \quad
$\includegraphics[width=.251cm]{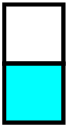} 
\maps (2,2)$, \quad
$ \square\;\square\maps (1^4)$.
\\
Let us note in passing: 
(1) Type $\pfii$ has two eigenvalues with degeneracies (3,1),
corresponding to the classical $U_qsl_2$ symmetry:
$2\otimes 2 = 3 \oplus 1$ (sometimes written
$\square\otimes\square=\square\!\square\oplus\emptyset$
\cite{drinfeld,Jimbo86}).
\\ (2) Type $\off$ can be specialised to $\alpha_1 = \alpha_2$ whereupon
(up to the trivial overall factor)
we
recover the representations given in \cite[\S12.1.1]{Martin91}
(cf. the Birman--Murakami--Wenzl algebra). 


In general $\omega(\varalpha)$ `factorises' into a contribution
$\phi(\alphaa)$
from $\off$ edges
and a contribution  $ \nu_i(\alphaa)$ from each nation:
$\omega(\varalpha) = \phi(\alphaa) \prod_{i} \nu_i(\alphaa)$.

\mdef \label{de:phinu}
The contribution of eigenvalues of the form $\pm \mu_{ij}$ from
$\off$ edges is given in general by the following.
Let nations
of a solution $\alphaa$ 
be indexed by some ordered set $I = I(\alphaa)$; and for $i \in I$ let
$N_i = |p_i|$ denote the size of nation $i$.
Then the part of the signature
of solution $\alphaa$
coming from edges between nations is
\[
\phi(\alphaa) \; =  \prod_{i<j \in I} (N_i N_j, N_i N_j)
\]
--- that is, there are $N_i N_j$ eigenalues $\pm\mu_{ij}$.
For example, If $\lambda = \square\!\square\;\;\square$ then we have just $(2,2)$
from these edges.
The number of free parameters $\mu_{ij}$
in the variety coming from this part of the matrix
is $ {{|I|}\choose{2}}  $.

\newcommand{\hash}{\#} 

Meanwhile
each nation $i$ contributes (one or) two eigenvalues, with multiplicities given
as follows. Let $\hash_i0$ denote the number of 0 edges.
Let set $J_i$ index the counties, i.e. the elements of $q|p_i$.
Then
$  \hash_i 0 \; = \;   
 \sum_{j\in J_i} {{|q_j|}\choose{2}}  $.
Let $s_2(i)$ denote the set of blue (colour 2) vertices.
We have
\beq \label{eq:nationsignature}
\nu_i(\alphaa) =
({{N_i}\choose{2}} + \hash_i 0 +|s_1(i)|,
 {{N_i}\choose{2}} - \hash_i 0 +|s_2(i)|) . 
\eq
Note: (i) That this does not depend on the order on counties.
(ii) That each one-nation representation is
(up to an elementary rescaling) a representation of some Hecke algebra.
\ignore{{
...by the number of vertices of each colour in that nation
plus \ppm{bla bla bla - write it out!
...Notation...}.
Thus (ignoring the single-county exceptions)
\[
\omega(\alphaa) \;  = 2...|...
\]
}}%

In practice we may notate $\omega(\alphaa)$ as
$(\nu_1 ; \nu_2 ; ... ;\nu_{|I|} \sep \phi)$.
\quad 
Now consider rank-3.


\mdef We start with
$\square\;\square\;\square =\off\off\off$ in $N=3$ as in Lemma~\ref{le:///}.
Here there are (up to) 9 distinct eigenvalues
($a_1, a_2, a_3, \pm\mu, \pm\nu,\pm\lambda$),
each of multiplicity 1;
given by 6 free parameters ($a_1, a_2, a_3, \mu, \nu,\lambda$), as shown.
We record this as the `signature'
$(1^9) = (1;1;1\sep 1,1,1,1,1,1)$.
\[
\ber{2}{3}{\;\;\;\;/_{\pm\mu} }{\;\;\;\; /_{\pm\nu} }{\;/_{\pm\lambda} \;}
\hspace{.1in}
\hspace{-.02471in} signature: \; (1^9),6
\]

\ppmx{[PAUL: Maybe demote this to an eigenvalue section later?]}

\mdef
For
$\includegraphics[width=.251cm]{xfig/tree123Pa2.eps}\;\square$,
so for example
$\off\off\uppfii$ as in Lem.\ref{le://+},
there is
1 eigenvalue of multiplicity 3;
2 eigenvalues of multiplicity 2;
and there 
are 2 eigenvalues of multiplicity 1, all given by 4 parameters.
Note that there is a variable $\beta$ associated to the $\pfii$ edge which does
not manifest elsewhere.
Note also that the object cannot be diagonalised over the ground ring,
or more specifically it would have a non-trivial Jordan block in the
specialisation $a_1 = \beta$.
\ignore{{
\hspace{.1in}
For $\off\off\uppai$
the $\uppai$ does not need such a decoration, since the
eigenvalues are given by the vertex eigenvalues.
}}



\newcommand{\zzz}{, \qquad}
\newcommand{\zzzz}{&}
\newcommand{\chomp}[3]{\!\!\!\begin{array}{c} #1 \\  #2 \\  #3 \end{array}\!\!\!}
\newcommand{\chmp}[3]{\!\!\!\begin{array}{c} #1  \\  #3 \end{array}\!\!\!}

\mdef
Altogether for $N=3$ (here giving just the degeneracies,
for a $\GG_N$-transversal cf. Prop.\ref{pr:N3}): 
\smallskip

\noindent
\begin{tabular}{|c|c|c|c|c|c|c|c|c|c|c|}
\hline
$\chomp{\square\!\square\!\square}{=000 }{\maps (9)}$
\zzzz
$\chomp{\includegraphics[width=.51cm]{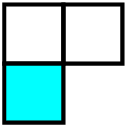}}{0\pai\pai}{(6,3)}$
\zzzz
$\chomp{\includegraphics[width=.51cm]{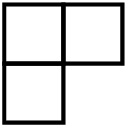}}{0ff}{(7,2)}$
\zzzz
$\chomp{\includegraphics[width=.251cm]{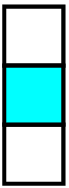}}{afa}{(5,4)}$
\zzzz
$\chomp{\includegraphics[width=.251cm]{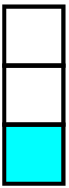}}{faa}{(5,4)}$
\zzzz
$\chomp{\includegraphics[width=.251cm]{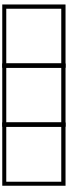}}{fff}{(6,3)}$
%
\zzzz
$\chomp{\square\!\square\;\square}{0\off\off}{(4;1:2,2)}$ 
\zzzz
$\chomp{\includegraphics[width=.251cm]{xfig/tree123Pa2.eps}\;\square}{f\off\off}{(3,1;1:2,2)}$ 
\zzzz
$\chomp{\includegraphics[width=.251cm]{xfig/tree123Pa3.eps}\;\square}{a\off\off}{(2,2;1:2,2)}$ 
\zzzz
$\chomp{\square\;\square\;\square}{\off\off\off}{  (1;1;1:1^6)}$  
\\ \hline
\end{tabular}

\begin{figure}
\includegraphics[width=15.6cm]{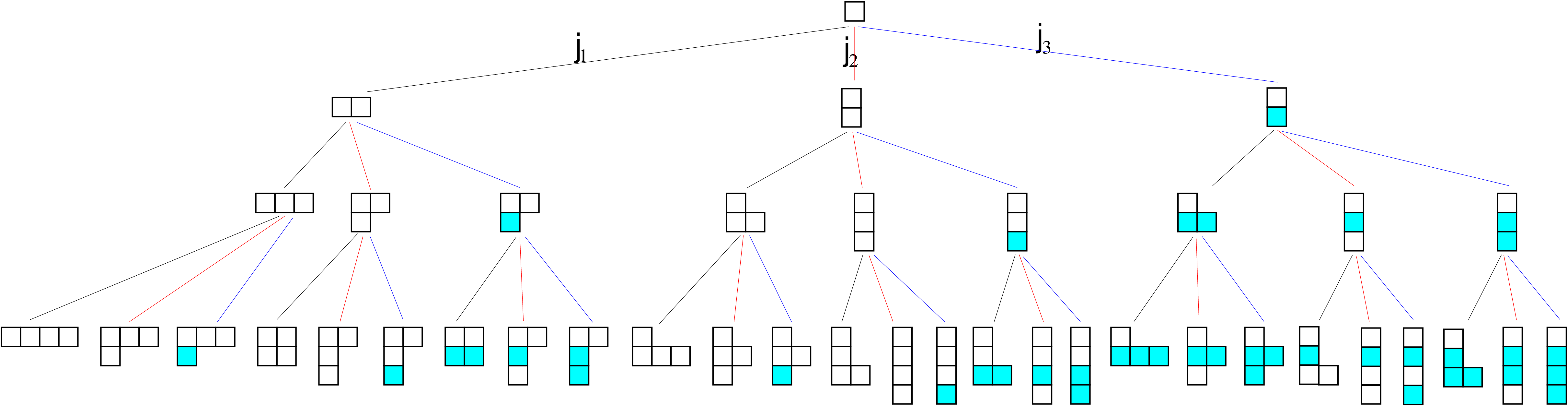}
\caption{One-nation configurations up to rank-4. \label{fig:N1234}}
\end{figure}

\smallskip \noindent
Notes:
(i) In general  
$|\TTT_N | = 1,4,13,46,154,...$
is the Euler `forest' transform of the one-nation sequence $3^{N-1}$
(cf. Fig.\ref{fig:N1234}).
This is A104460 in \cite{sloane1}. 
From the $\PPP$-orbits in 
Fig.\ref{fig:tree123}, 
the $\GG_N$-transversal for $N=3$  
here is 3 elements smaller than
the $\Sym_N$-transversal, thus ten entries.
\\
(ii) The two (5,4) cases are
related as per (\ref{de:phinu})(i), but are
not related by our symmetries,
but they are similar as matrices by a non-local conjugation.
\\
(iii) The two (6,3) cases have different Jordan structure over the polynomial ring
of indeterminates.

\medskip


\mdef
For $N=4$  
the one-nation cases (cf. Fig.\ref{fig:N1234},
and taking account of (\ref{de:phinu})(i))
and some indicative multination cases are: 

\noindent
\begin{tabular}{|c|c|c|c|c|c|c|c|c|c|c|c|c|c|}
\hline
$\chmp{\square\!\square\!\square\!\square}{000000}{\maps (16)}$
\zzzz
$\chmp{\includegraphics[width=.75cm]{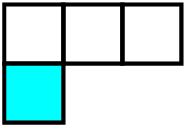}}{000aaa}{ (12,  4)}$
\zzzz
$\chmp{\includegraphics[width=.75cm]{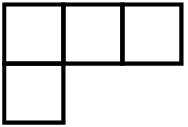}}{000fff}{ (13,  3)}$
\zzzz
$\chmp{\includegraphics[width=.5cm]{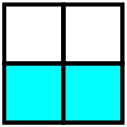}}{0aaaa0}{{(8,8)}}$ 
\zzzz
$\chmp{\includegraphics[width=.5cm]{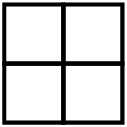}}{0ffff0}{ (12,  4)}$ 
\zzzz
$\chmp{\includegraphics[width=.5cm]{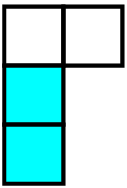}}{aa0faa}{(9,  7)}$ 
\zzzz
$\chmp{\includegraphics[width=.5cm]{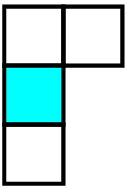}}{aa0aff}{(10,  6)}$ 
\zzzz
$\chmp{\includegraphics[width=.5cm]{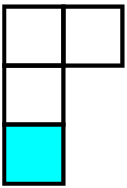}}{ff0aaa}{(10,  6)}$ 
\zzzz
$\chmp{\includegraphics[width=.5cm]{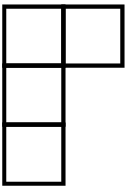}}{ff0fff}{{(11,5)}}$ 
\zzzz
$\chmp{\includegraphics[width=.27cm]{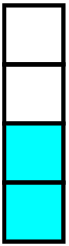}}{faaaaf}{(8, \; 8)}$
\zzzz
$\chmp{\includegraphics[width=.27cm]{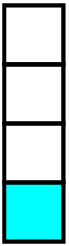}}{fffaaa}{(9, \; 7)}$
\zzzz
$\chmp{\includegraphics[width=.27cm]{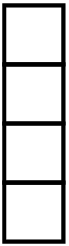}}{ffffff}{(10, \; 6)}$
\\ \hline \end{tabular}

\noindent
\begin{tabular}{|c|c|c|c|c|c|c|c|c|}
\hline
%
$\chmp{\square\!\square\!\square\;\square}{000\off\off\off}{ (9;1:3,3)} $
\zzzz
$\chmp{\includegraphics[width=.51cm]{xfig/tree123Pa12.eps}\;\square}{}{(7,2;1\sep 3,3)}$
\zzzz
$\chmp{\square\!\square\;\square\!\square}{0\off\off\off\off 0}{ (4;4:4,4)} $
\zzzz
$\chmp{\includegraphics[width=.2751cm]{xfig/tree123Pa2.eps}\;\square\!\square}{}{(3,1;4\sep 4,4)}$
\zzzz
$\chmp{\includegraphics[width=.2751cm]{xfig/tree123Pa3.eps}\;\square\!\square}{}{(2,2;4\sep 4,4)}$
\zzzz
...
\zzzz
$\chmp{\square\;\square\;\square\;\square}{\off\off\off\off\off\off}{
(1;1;1;1:1^{12})}$
\\ \hline \end{tabular}




\subsection{On $\Match^N$ representation theory} \label{section:repthy} $\;$ 
\module{braid-repthy}

\label{ss:RT101}

$\;$



In Artinian representation theory the natural notions of isomorphism and of
additivity `agree' with each other so that the classification is canonical.
In our case there are layers of naturality and the tie-up is more complex
(cf. \cite{MazorchukMiemietz14,RumyninWendland} etc).


Given a representation
$F$ on $\Match^N$
(we will just say a representation on $N$)
%
there is,
via Lemma~\ref{lem:chacon1}III,
a sense in which
each subset $S \subseteq \ul{N}$ gives a subrepresentation $F|_S$.
I.e. every subset
gives  
an invariant subset, with a `complementary' representation
on $\ul{N}\setminus S$.
However, there is also a novel form of additivity
--- an `adjoint' to this `restriction'.
It does not make immediate sense to talk of adding reprentations
(as in direct sums), or in particular of resumming the representations on $S$
and $\ul{N}\setminus S$.
However there is an additive property:
$\;$
Given representations on $N$ and $M$ we obtain a representation on
$N \sqcup M \hookrightarrow N+M$
by completing the $K_N \sqcup K_M$ subgraph of the complete graph
$K_{N \sqcup M}$ with $\off$ edges.
{
({\em Proof}:
every new triangle contains two $\off$s, and now use
the second row of configurations from Prop.\ref{pr:N3},
via Lemma~\ref{lem:3 suffice},
noting that the orientation of $\off$ edges
--- the relative placement of the singleton vertex of the three ---
does not here affect the validity of the solution.)
}
$\;$
This
addition
is substantially different from direct sum --- for one thing these edges
carry parameters so this gives not a single isomorphism class of representations
in the classical sense but a variety.
Let us write $F \oslash G$ for the
variety 
of representations constructed from
$F$ on $N$ and $G$ on $M$ say --- noting that many choices are needed to select
a concrete representation from this
variety. 
What is   {\em not}  true in general is that
$F|_S \oslash F|_{\ul{N}\setminus S} \ni F$.
For this to hold, we require that $S$ is a union of nations of $F$ ---
so that a nation is a kind of
irreducible representation.

In this sense
classifying representations  first by 
 classifying single-nation representations and then
 giving composition
multiplicities of general representations in terms of such
`irreducible' content
is roughly in analogy with Jordan--Holder.





\section{On rectanguloid/X-equivalence and Lemma~\ref{lem:Xeq}} \label{ss:Xeqa}

\module{XlemmaProof}



\ignore{{  
\newpage

\beff{}{{

\begin{lemma} \label{X}
Let $\FF:\Bcat\rightarrow\Match^N$ be a representation, hence given by
$\FF(\sigma)$.
Let $X$ be a diagonal invertible matrix in $\Match^N(2,2)$.
Then $\FF'(\sigma)= X \FF(\sigma) X^{-1}$ also gives a representation.
The relation $\FF\sim\FF'$
if there exists any such $X$ relating $\FF$ and $\FF'$
is an equivalence relation.
We call it X-equivalence.

Let $\FF(\sigma)$ be a solution,
with the edge submatrices
$  \mat{cc} a & b_{}\\ c & d \tam_{ij}
=
\mat{cc} a_{ij} & b_{ij}\\ c_{ij} & d_{ij} \tam_{}
$
($ij=12,13,23,...$).
Then changing these (for each edge $ij$ independently) to 
$  \mat{cc} a & xb_{}\\ c/x & d \tam_{ij} $
(any invertible $x_{ij}$)
and hence in particular to 
$  \mat{cc} a & cb_{}\\ 1 & d \tam_{ij} $
(`lower-1 form')
also gives
a solution.
\\
In general writing $\FF'(\sigma)$ for the changed solution we say that
$\FF$ is $X$-related to $\FF'$. This is an equivalence relation,
partitioning the set of all solutions into `$X$-equivalence classes'.
We write $\FF(\sigma)
\stackx  
\FF'(\sigma)$ for $X$-equivalent
solutions. 
\\
(Note that $\FF'$ is a solution not only in the same varietal branch, but with the
same spectrum as $\FF$.)
\end{lemma}
}}
}}


\beff{}{{

\mdef  \label{proof:X}  
{\bf Explicit proof of Lemma~\ref{lem:Xeq}.}
Observe from Theorem~\ref{th:constraints} that the
constraint equations
are satisfied by the substitutions
\[
\mat{cc} a_{ij} & b_{ij}\\ c_{ij} & d_{ij} \tam_{}
\leadsto \; 
\mat{cc} a_{ij} & x_{ij} b_{ij}\\ c_{ij} / x_{ij} & d_{ij} \tam_{}
\]
(the upshot of 
{a conjugation by $X\in\Match^N(2,2)$ as in the Lemma})
given that they are satisfied by
$\FF(\sigma)$, because the entries  
$b_{ij}, c_{ij}$,
if appearing with non-zero coefficient, appear as
$b_{ij} c_{ij}$,
or else as an overall factor.
\hspace{.1in} 
The equivalence relation property follows by construction.
{\qed}

}}

\ignore{{
\begin{remark}
\ppm{[delete this version and keep later one?]}
We observe that, for $N\geq 3$ the transformation $X\FF(\sigma)X^{-1}$ with $X$ an arbitrary diagonal cannot be achieved by $X$ of the form $Z\otimes W$.
Indeed, let $M\in \Match^N$ so that
$\alpha(M)=A=(a_1,\ldots,a_N,A(1,2),\ldots,A(N-1,N))$ as in eqn. (\ref{eq:alconv1}).
Suppose further that there are $i<j<k$ so that $b_{ij},b_{ik}$ and $b_{jk}$ are non-zero, without loss of generality these can be chosen to be $1<2<3$, so that by a judicious choice of $X$ we may transform the corresponding $b_{12},b_{13},b_{23}$ entries to be $1$.  Now denote the non-zero entries of diagonal matrices $Z$, $W$ by $z_i$ and $w_i$ respectively.  Then the corrresponding entries of $(Z\otimes W)M(Z\otimes W)^{-1}$ are of the form: $\frac{z_jw_ib_{ij}}{z_iw_j}$ for $i<j$ among $\{1,2,3\}$.  Setting each equal to $1$ we find that $b_{ij}=\frac{z_iw_j}{z_jw_i}$.
But this implies that $b_{12}b_{23}=b_{13}$, a relation that does not hold for general $M$.
\end{remark}
}}

\ignore{{
\whisper{
Remark.
Eric introduces a diagonal matrix that acts in rank-2 and
observes that the X change can be achieved there specifically by conjugating by this.
}
}}



\module{XlemmaNotes}


\newcommand{\mum}[1]{\mu^{(#1)}} 
\newcommand{\mumm}[1]{\overline{\mu}^{(#1)}}
\newcommand{\FFS}{F(\sigma)}  

We note that
neither
general solutions to the Yang-Baxter equation,
{nor spin-chain representations of other categories,}
are  
stable under conjugation by {non-local} diagonal matrices.
{--- Meaning that a solution
$\FF: \Bcat \rightarrow \Mat^N$
given in the form $\FF(\sigma)$ is not
necessarily taken to  
a solution unless $X$ of the form $X=\mu\otimes\mu$
(i.e. unless $\FF\leadsto \FF\circ\cC^\mu$).}

\ignore{{
However, since diagonal matrices are charge conserving, the set of
charge conserving matrices is stable under conjugation by diagonal matrices.
\ppm{[-so what? I think we need to say a bit more here?
the condition is something like the following, right?
(I am just going to temporarily give this a bit of room to breath.
I guess we'll delete it all again once we are done.)]}
}}


\mdef {
Formally a conjugating functor associates a suitable invertible matrix
$\mum{m}$ to each object $m$, and transforms morphism according to:
$f \mapsto \mum{s(f)} f \mumm{t(f)}$, where $s,t$ are the source and tail
functions.
This obviously behaves well under category composition (so gives a functor),
but not necessarily under monoidal composition (so does not give a monoidal
functor in general).
A necessary condition is that it acts on the whiskerings of
$\FFS \in \Match^N(2,2)$
like
$\mum{3} \FFS\otimes 1 \mumm{3} = \mum{2} \FFS \mumm{2} \otimes 1 $
and
$\mum{3} 1\otimes \FFS \mumm{3} = 1\otimes \mum{2} \FFS \mumm{2} $.
This of course holds if $\mum{m} = \mu^{\otimes m}$ as before - the `local' solution;
but we are looking for something more. }
{
An interesting
question is if this can be made to work with some degree of generality,
or if the only possibilities arise when we pass all the way to the constraints
coming from the YBE. 
\\
Observe that a conjugating functor, where defined, is determined by its action on
a generating set of morphisms,
and hence on the $\mumm{m}$'s for objects in their sources and tails
(so in our case by $\mumm{2}$
--- the $X$ of the X~Lemma).
It is an interesting question then, in our case,
given a valid $\mumm{2}$, what is $\mumm{m}$?
Our setup
is a nice  
playground for this.
}
Explicitly, taking the  in-line form
of $A \in \Match^N(2,2)$ 
as in (\ref{eq:alconv1})
\[
\alphain(A) =
(a_1, a_2, \ldots, a_N, \mat{cc} a_{12} & b_{12} \\ c_{12} & d_{12} \tam,\ldots,
\mat{cc} a_{ij} & b_{ij} \\ c_{ij} & d_{ij} \tam,\ldots),
\]
conjugation by an invertible diagonal matrix defined by 
$X(\mathbeef_i\otimes \mathbeef_j)=y_{ij}(\mathbeef_i\otimes \mathbeef_j)$ gives
\begin{equation} \label{eq:alphaX}
\alphain( XAX^{-1} ) =
(a_1, a_2, \ldots, a_N, \mat{cc} a_{12} &x_{12}b_{12} \\ c_{12}/x_{12} & d_{12} \tam,\ldots,
\mat{cc} a_{ij} &x_{ij}b_{ij} \\ c_{ij}/x_{ij} & d_{ij} \tam,\ldots)
\end{equation}
where $x_{ij}:=\frac{y_{ij}}{y_{ji}}$ for $i<j$.

The local case $X=\mu\otimes\mu$, with $\mu \mathbeef_i = y_i \mathbeef_i$ say, gives
$y_{ij} = y_i y_j$ and hence lies in the center
of $\Match^N(2,2)$
and acts trivially.
By virtue of the Lemma
(i.e. for the case of $\Braid$ and possibly not in general)
we have also
the decomposable case  $X=\mu\otimes\nu$, with
{$\nu \mathbeef_i = z_i \mathbeef_i$} say;
and the entangled case.
Just to see if it has interesting special features, let us consider
the case $X=\mu\otimes\nu$. This  gives
$y_{ij} = y_i z_j$ and hence
$x_{ij}=\frac{y_{i}z_{j}}{y_{j}z_{i}}$.
For $N=3$ we have $x_{12} = \frac{y_1}{z_1}\frac{z_2}{y_2}$, 
$x_{13} = \frac{y_1}{z_1}\frac{z_3}{y_3}$,
$x_{23} = \frac{y_2}{z_2}\frac{z_3}{y_3}$.
{This is clearly not generic.  
An interesting question is
whether X-equivalence can be realised `functorially'
(in particular via a conjugating functor)
in general.
}

\ignore{{
From  (\ref{eq:alphaX}) 
and the X~Lemma we see that
\ppm{given a solution $F(\sigma)$}
we may simultaneously and independently rescale the off-diagonal entries
$(b_{ij},c_{ij})$ of each of the $2\times 2$ matrices to any pair of the form
$(xb_{ij},c_{ij}/x)$.
\medskip
}}




\section{Details for proof of Theorem~\ref{th:JJN}} \label{ss:deets}
\label{section:repthy1}

\module{repsTransvers} 


\ignore{{
Recall that our Theorem~\ref{th:mainx} and \ref{th:JJN}
raises the question of the appropriate notion of
equivalence of representations in this setting.
As mentioned, in  part this comes down to the symmetries intrinsic to the target.
}}

\ignore{{
Here we  
work with the $\Sym_N$ symmetry.
(This  reduces the size of the set of classes
of solutions by a finite factor only, but  
in this way we are able to give a
{\em complete} characterisation
of the solution set - by Euler transforming the elementary sequence $3^{N-1}$.)
}}


\subsection{On further properties of $\Sym_N$-orbits of solutions} $\;$ 
\ppmx{What has happened to the all-+ Lemma and the 6-rule in this version, paul?!}


\newcommand{\Norb}{\Sigma_N \times \Zz_2}
\newcommand{\nominus}{$-$-free}  

\beff{}{{


\mdef  \label{6rule}
{\bf 6-rule}:
If an oriented chain of two edges is signed with the same sign
in
a solution
$\FF(\sigma)$
then the
`long' edge completing the triangle is signed with the same sign.
--- To see this check the \ninerule\ list
from
\ref{pr:N3}   
and
the orbit of +++ in Table~\ref{tab:1}. 
\\
{
(This is called 6-rule since it reduces the number of possible edge $\pm$-colourings
of a triangle from eight to six. It also reduces the number of $\pm0$-colourings.)
}
}}

\begin{lemma} \label{lem:no-}
Every $\Sigma_N$  
orbit of solutions in $\Rec(\germ_N)$ 
contains an element with no $-$-orientation.
\\
Furthermore, in such a
\; {`\nominus'} configuration the elements of a 0-part
(a county) 
are clustered with respect
to the natural order in their nation, i.e. they are consecutive. 
\end{lemma}

\newcommand{\dddots}{...}

\prooff\ 
First we note that the claim holds if it holds for 1-nation solutions,
since there are perms that act non-trivially only on a single nation and the $\off$-edges
out of that nation, restricting to a complete set of perms on that nation.
So now consider a single nation.
We work by induction on $N$. The claim is true for $N<4$ by inspection of our explicit
solution sets.
Suppose true at level-$N$ and consider $N+1$.
WLOG by inductive assumption consider a configuration with all (signed) edges $+$
between
vertices $1,2,\cdots,N$
(any non-signed edges are 0).
\ignore{{
Observe that configurations with 0s
on all these edges must have the same sign on all remaining edges, edges to $N+1$.
If all + then we are done. If all $-$ then
the transposition $(1 \; N\!+\! 1)$ \ppm{...} 
}}
Consider what configurations of edges to vertex $N+1$ are allowed here. 
Neglecting cases with 0s on the  edges to $N+1$ for a moment we have, say,
as on the left here:
\[
\includegraphics[width=5.4\ccm]{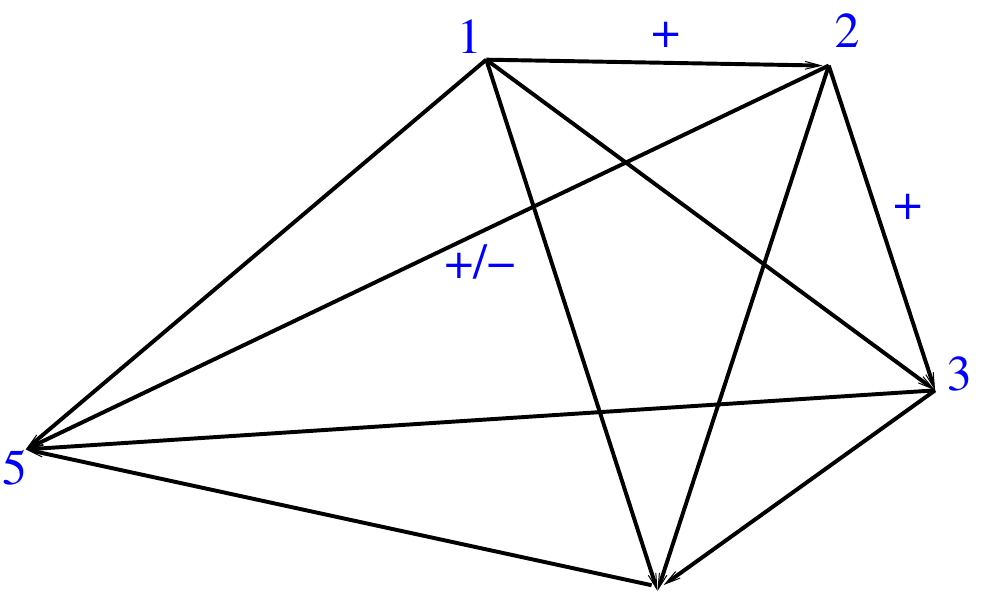}
\hspace{.5\ccm}
\includegraphics[width=5.4\ccm]{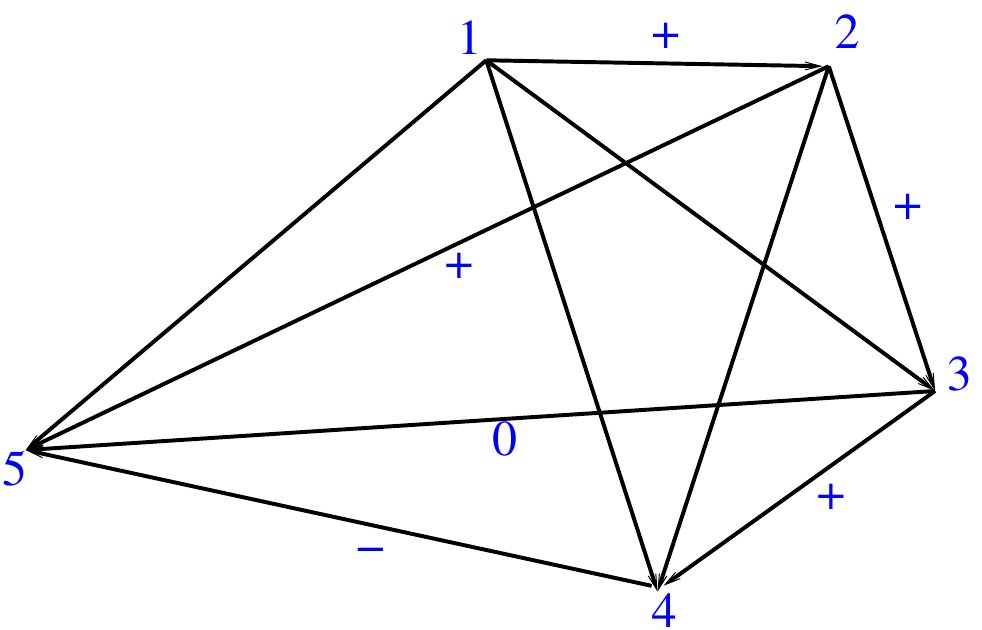}
\]
By the 6-rule, if the $i$ to $N+1$ edge is + then,
since the $i-1$ to $i$ edge is +,
the  $i-1$ to $N+1$ edge is also +.
Indeed this is also forced in the case where the $i-1$ to $i$ edge is 0.
Thus there are $N$ possible configurations of form $++...+--...-$
(NB all these arrows point to $N+1$).
The $i$-th of these configurations is taken to
$+++\dddots +$ by the perm
$(i\; i\!\!+\!\!1\cdots N)$,
so we are done in these cases.
In case of 0 on a long edge $i$ to $N+1$, the edges $j$ ($j<i$) to $i$ and
$j$ to $N+1$ must be the same, so these edges satisfy the no-$-$ claim.
Thus it remains only to address the edges after the last 0.
Suppose this last 0 is $i$ to $N+1$.
The edges $k$ to $N+1$ with $k>i$ will all be opposite to the corresponding
$i$ to $k$ edge (which is either + or 0), so either $-$ or 0.
Now apply the perm $(i\; N\!+\!1) $. This has the effect of making the $N+1$ vertex
the `new' $i$ vertex and incrementing $i$ et seq. Thus all the edges that were
$k$ to $N+1$ are reversed, so all the $-$s become +,
while the relative positions of all other vertices are preserved, so signs do not change
--- i.e. they remain \nominus. 

The second claim follows from the 6-rule: 
given $ij$ and $jk$ are + then $ik$ is +, so not 0.
\qed



\mdef \label{pa:order-pres}
{Note from (\ref{pa:rec})
that recipe $\Rec$ assigns a configuration to 
a $K_2$
subgraph dependent on the
{\em relative} rather than absolute values of the vertex labels.
\\
In this sense an  
inclusion
$\psi:\ul{M} \hookrightarrow \ul{N}$
as in \ref{lem:chacon1}III
that is  {\em order preserving}
induces
a `stronger' equivalence between a configuration on $M$ and
the restriction of a configuration on $N$ to
the 
image  $\psi(\ul{M})$
than for a generic inclusion. 
Such an order-preserved image is essentially identical
- changed only in being moved.
\\
Similarly, for a bijection $\psi : \ul{N} \stackrel{\sim}{\rightarrow} \ul{N}$
we may partition
the domain (not necessarily uniquely) into a set $p$ of subsets on
each of which $\psi$ is
order-preserving
(e.g. the set of singleton subsets is trivially preserved by any $\psi$)
and then $\psi$ moves these subsets `parallel-ly' around in $\ul{N}$
and the image of each under $\Rec$ gives an identical configuration.
Meanwhile, although the whole of
the image under
$\psi$
of a solution
is another solution, while it is isomorphic,
it is generally differently configured on edges {\em between} parts of $p$.
}


\vspace{.1cm}


\module{proofJJN} 



\mdef \label{pa:proofJJN0}
On the proof of Theorem~\ref{th:JJN}.
The Theorem asserts that every point on every variety of solutions may be
obtained by applying $\Rec$
from (\ref{pa:rec}) 
to the
image of some point in $\TTT_N^\C$ under some $w \in \Sym_N$.
This begs the question of how $\Sym_N$ acts on $\SSS_N^\C$.
%
Let $w \in \Sym_N$ act on $\SSS_N$ by permuting vertex labels.
Ignoring parameters for a moment it will be clear that $\TTT_N$ is a transversal
of $\Sym_N$-orbits of $\SSS_N$.
Thus we will be done if we can show that there is an action
on $\SSS_N^\C$
that
reduces to this,
and such that
the square:
\beq \label{eq:commsq}
\xymatrix{ \SSS_N^\C \ar[r]^-{\Rec} \ar[d]^{w\in\Sym_N}
                    &    \functor(\Bcat,\Match^N) \ar[d]^{w\in\Sym_N}
\\
\SSS_N^\C \ar[r]^-{\Rec}   &    \functor(\Bcat,\Match^N)
}
\eq
(where the action of $w$
on the right is via the functor in $\functor(\Match^N,\Match^N)$) 
commutes up to trivial movement on the variety.

\ignore{{ 

In addressing the action with parameters it will
be convenient to have names for them  
--- for any given solution they are
just given complex numbers, but 
in manipulating solutions in the form $\Rec(\lambda)$, $\lambda\in \SSS_N^\C$,
we will not generally want to specify the complex values of the parameters explicitly.
However there is not generally an easy way of {\em naming} them
- to treat specialisations collectively - that plays nicely with manipulations such as
$\Sym_N$ action.
This is because identical nations are ab initio unordered - while of course actual
asignments of complex parameters break this symmetry.
Once an element of $\jjj(J_N(\{1,2,3\}))$ has vertex numbers inserted this also breaks the symmetry.
But parameters `move' with the $\Sym_N$ action.
For example  
$\raisebox{-.1in}{
\includegraphics[width=.41cm]{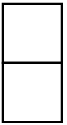}\;\;
\includegraphics[width=.41cm]{xfig/young11.eps}\;}
$
has type-space with parameters
$(\alpha_1, \beta_1, \alpha_2, \beta_2, \mu_{12})$
but there is no sense in which the nations are ordered.
With an insertion
$\raisebox{-.1in}{
\includegraphics[width=.41cm]{xfig/tree012ex2ab1b.eps}\;\;
\includegraphics[width=.41cm]{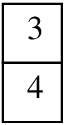}\;}
$
we can (say, child-first) order, but the action of $(13)\in \Sym_4$ would obfuscate
this choice.
There is no problem here - complex numbers do not need names.
But if we do not want to specify particular complex numbers we must name them, and then
we can only give a general characterisation of the $\Sym_N$ action on solutions by
describing {\em changes} rather than absolutes.
A resolution is then as follows.
Given $\raisebox{-.1in}{
\includegraphics[width=.41cm]{xfig/tree012ex2ab1b.eps}\;\;
\includegraphics[width=.41cm]{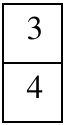}\;}
$,
the parameter names are ascribed in the obvious way - thus $\alpha_1 \in \C$
and so on:
\[
\typ{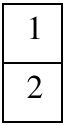}{young11-34}
\ignore{{ \;
\arr{c}
\includegraphics[width=.41cm]{xfig/tree012ex2ab1b.eps}
\\ \alpha_1 \\ \beta_1
\rra \;
\stackrel{---}{\mu_{12}}\;
\arr{c}
\includegraphics[width=.41cm]{xfig/young11-34}
\\ \alpha_2 \\ \beta_2
\rra }}
\]
- indicating $\alpha_1  , \beta_1$ ascribed to the first nation and so on.
The action of $(13)$ is:
\beq \label{eq:act13}
\typ{tree012ex2ab1b}{young11-34}
\ignore{{
\arr{c}
\includegraphics[width=.41cm]{xfig/tree012ex2ab1b.eps}
\\ \alpha_1 \\ \beta_1
\rra \!
\stackrel{\rule{25pt}{.4pt}}{\mu_{12}} \!
\arr{c}
\includegraphics[width=.41cm]{xfig/young11-34}
\\ \alpha_2 \\ \beta_2
\rra }}
\;\;\;\;
\stackrel{(13)}{\longmapsto} \;\;\;\;
\typ{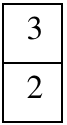}{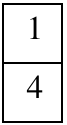}
\ignore{{ \;
\arr{c}
\includegraphics[width=.41cm]{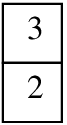}
\\ \alpha_1 \\ \beta_1
\rra 
\stackrel{---}{\mu_{12}}
\arr{c}
\includegraphics[width=.41cm]{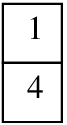}
\\ \alpha_2 \\ \beta_2
\rra }}
\eq
So the complex parameter that we (arbitrarily) named $\alpha_1$ is now associated to
vertices 3 and 2.
NB it was the use of symbol 1 to name it that was the anomaly.
As a complex number it is free to move under the action.


Recall we aim to show that
\eqref{eq:commsq} commutes.
\ignore{{
\[
\xymatrix{ \SSS_N^\C \ar[r]^-{\Rec} \ar[d]^{w\in\Sym_N}
                    &    \functor(\Bcat,\Match^N) \ar[d]^{w\in\Sym_N}
\\
\SSS_N^\C \ar[r]^-{\Rec}   &    \functor(\Bcat,\Match^N)
}
\]
\ppm{
commutes, where the action of $w$ 
on the right is via the functor in $\functor(\Match^N,\Match^N)$.}
}}
To see this note that it will be enough to check on triangles.

The action of $\Sym_N$ from $\functor(\Match^N,\Match^N)$
on  any matrix in $\Match^N(2,2)$, and hence in particular on a solution,
can be deduced from the action on a triangle $i<j<k$
(say  with  $ijk=123$)
as follows.
\ignore{{
The $ijk$-localised solution matrix is of form
\[
\mat{ccc|ccc|ccc}
a_{1}& & & & & & \\
 &a_{12}& &b_{12}& & & \\
 & &a_{13}& & & &b_{13}& \\ \hline
 &c_{12}& &d_{12}& & &  \\
 & & & &a_{2}& & \\
 & & & & &a_{23}& &b_{23} \\ \hline
 & &c_{13}& & & &d_{13}& \\
 & & & & &c_{23}& &d_{23} \\
 & & & & & & & &a_{3}
\tam
\]
(here with $ijk=123$).
}}
The action for example of $(13)\in\Sym_N$ is conjugation by
$\mat{ccc}0&0&1\\ 0&1&0\\1&0&0\tam \otimes \mat{ccc}0&0&1\\ 0&1&0\\1&0&0\tam $,
thus on any charge-conserving matrix:
\[
\mat{ccc|ccc|ccc}
a_{1}& & & & & & \\
 &a_{12}& &b_{12}& & & \\
 & &a_{13}& & & &b_{13}& \\ \hline
 &c_{12}& &d_{12}& & &  \\
 & & & &a_{2}& & \\
 & & & & &a_{23}& &b_{23} \\ \hline
 & &c_{13}& & & &d_{13}& \\
 & & & & &c_{23}& &d_{23} \\
 & & & & & & & &a_{3}
\tam
 \stackrel{(13)}{\leadsto}
\mat{ccc|ccc|ccc}
a_{3}& & & & & & \\
 &d_{23}& &c_{23}& & & \\
 & &d_{13}& & & &c_{13}& \\ \hline
 &b_{23}& &a_{23}& & &  \\
 & & & &a_{2}& & \\
 & & & & &d_{12}& &c_{12} \\ \hline
 & &b_{13}& & & &a_{13}& \\
 & & & & &b_{12}& &a_{12} \\
 & & & & & & & &a_{1}
\tam
\]
Note that this
verifies the action $w a_{ij} = a_{w(i) w(j)}$ as in
Theorem~\ref{th:constraints}.
\\
On an edge $(i,j)$ the action is (localised at this edge) as follows.
\\
Firstly the transposition $(i,j)$
gives 
$\mat{cc} a_{ij} & b_{ij}\\c_{ij}&d_{ij}\tam
 \stackrel{(ij)}{\leadsto}
 \mat{cc} d_{ij} & c_{ij}\\b_{ij}&a_{ij}\tam$.
\\
On a signed edge this is
$\mat{cc} a_{ij} & b_{ij}\\c_{ij}&0\tam
 \stackrel{(ij)}{\leadsto}
 \mat{cc} 0 & c_{ij}\\b_{ij}&a_{ij}\tam$
 which, up to X-equivalence, is simply a sign change.
So here the claim follows
from the use of the natural order in the construction of $\Rec$
in \ref{pa:rec}(iv).
On a 0 edge agreement is trivial. On a $\off$ edge we have agreement up to
X-equivalence as required.
\\
Secondly consider transposition $(i,k)$ on an $(i,j)$ edge:
%
$\mat{cc} a_{ij} & b_{ij}\\c_{ij}&d_{ij}\tam
  \stackrel{(ik)}{\leadsto}
 \mat{cc} a_{kj} & b_{kj}\\c_{kj}&d_{kj}\tam
 =  \mat{cc} d_{jk} & c_{jk}\\b_{jk}&a_{jk}\tam
$
 (noting $k>j$).

}}

\begin{figure}
\[
\includegraphics[width=7.05cm]{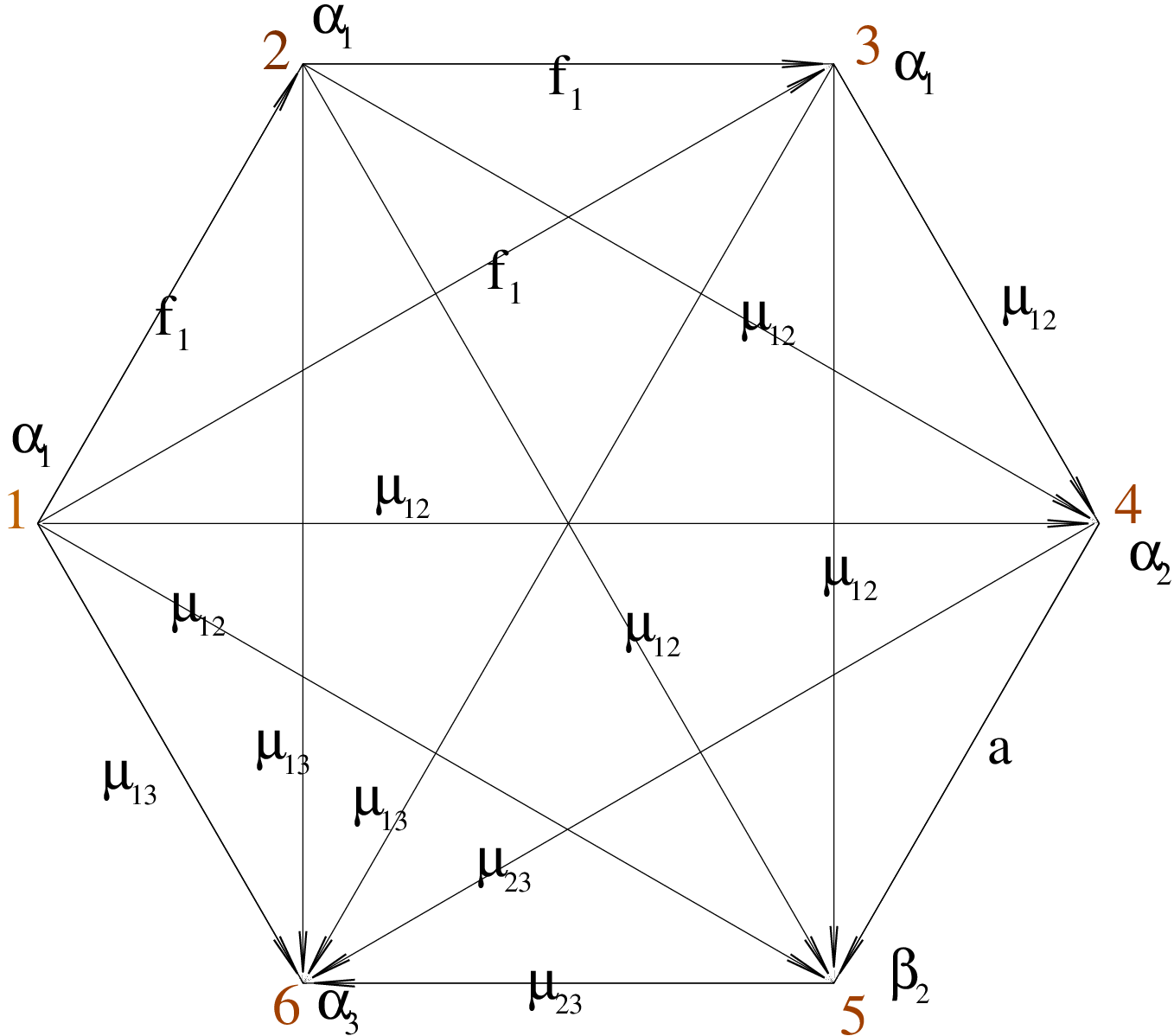} 
\raisebox{2.952cm}{$\;\;\;\;\stackrel{(34)}{\leadsto}\;\;\;\;$}
\includegraphics[width=7.05cm]{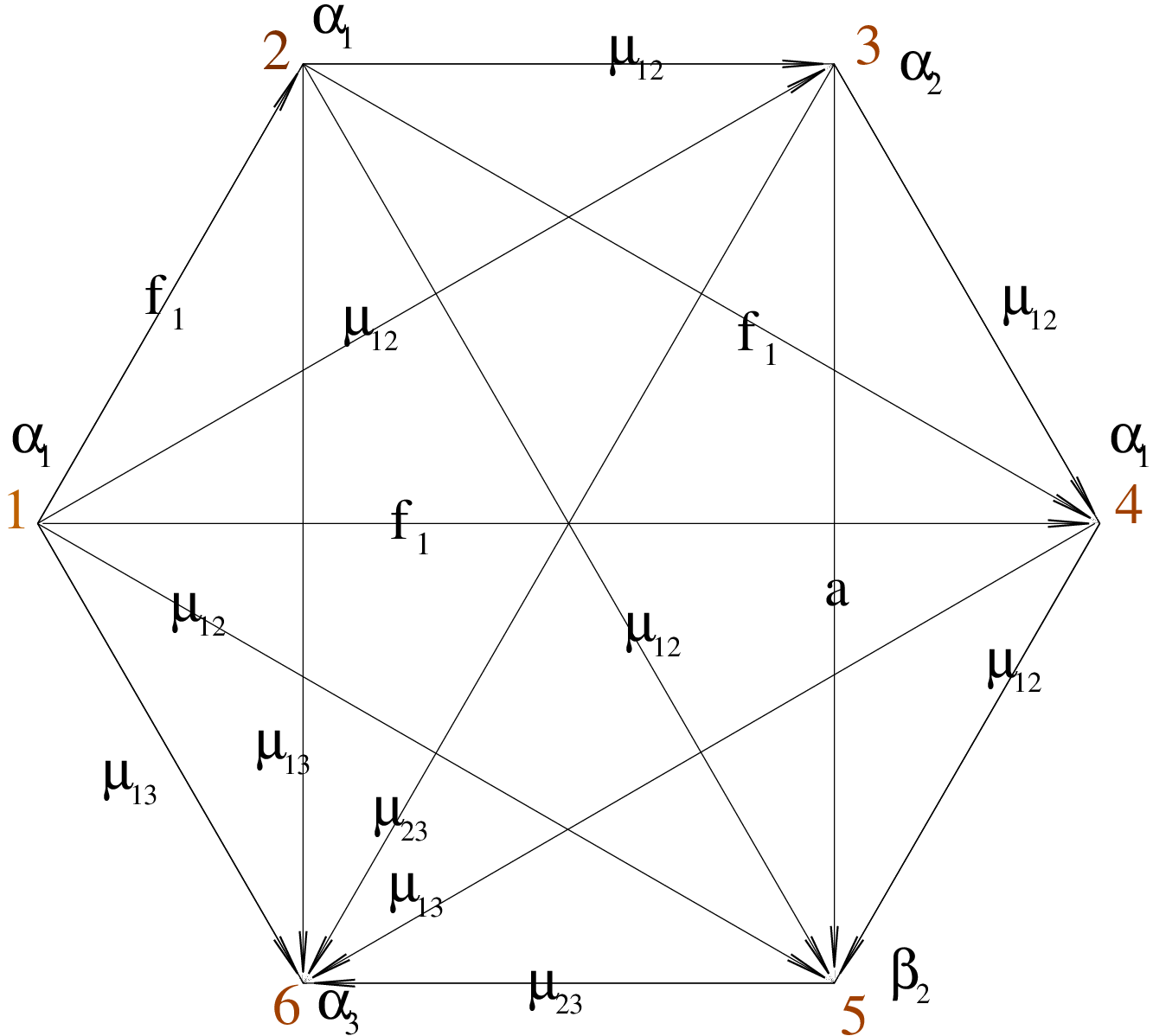} 
\]
\caption{Action of $(34)\in \Sym_N$ on a solution with $N=6$ and three nations.
On $\off$ edges we simply give (names of form $\mu_{ij}$ for) the parameters.
We write $\pfii_1$ for $\pfii_{\beta_1}$ to avoid clutter. 
\label{fig:hex34}}
\end{figure}

Figure~\ref{fig:hex34} illustrates the action of $\Sym_N$ on an example solution.
%
The change expressed at the level of $\SSS_N^\C$ is
\[
\typp{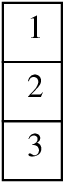}{1}
\hspace{.21cm}
\typp{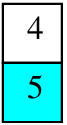}{2}
\hspace{.21cm}
\typpp{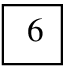}{3}
\hspace{.21in}
\raisebox{.652cm}{$\;\;\;\;\stackrel{(34)}{\leadsto}\;\;\;\;$}
\hspace{.21in}
\typp{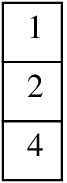}{1}
\hspace{.21cm}
\typp{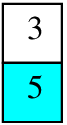}{2}
\hspace{.21cm}
\typpp{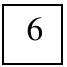}{3}
\]
(NB the specific {\em names} $\alpha_1$ etc used
here to represent parameter values are not
portable in general 
--- only which are `related' and which are independent.
To give a point in type-space we would need to order the component nations.)
Ignoring parameters, i.e. at the level of $\SSS_N$, this is of course
the same as the direct action.
%
Now
observe that up to `signs' the collection of triangle
configurations are simply rearranged on $K_N$
by the permutation.
For example the fff triangle moves from vertices 123 to 124
under the action of $(34)\in \Sym_N$.
(In the illustrated case there are no sign changes since
3,4 are adjacent and the 34-edge is $\off$-decorated.)
Since triangles
collectively
are at most rearranged up to signs, the relationships of $\alpha_i$
and $\beta_i$ parameters are preserved, so the open conditions defining the type-space
are preserved.


\ignore{{ 

\ppm{
 ????
If $i<k$ then this simply swaps-in the $ik$ edge:
$\mat{cc} a_{ij} & b_{ij}\\c_{ij}&d_{ij}\tam
 \leadsto \mat{cc} a_{ik} & b_{ik}\\c_{ik}&d_{ik}\tam$ ????}
Comparing with the other side of the square,
using \eqref{eq:act13} say,
we end up with 
$ \mat{cc} 0 & \mu_{12}\\ \mu_{12}&0\tam$
going either way.
The other cases are: 123 in same nation
(subcases ... as per rule-of-9?);
or 123 all in different nations... .
%
...
\ppm{-FINISH.
If $i>k$ then ...
Where $k$ is in the same nation, this is ...;
\\ finally where $k$ is in a different nation ... FINISH.
 [-could use example following \ref{de:TTTN}?]
Theorem~\ref{th:JJN} follows immediately. \qed}


A complete worked example can be found in \S\ref{sec:N4}. 
}}


%
\subsection{2-coloured composition combinatorics} $\;$
\module{notation-composition2}



{In light of Lemma~\ref{lem:no-} we can characterise orbits of representations in
terms of 
suitable partitioned integer compositions
--- describing nations in an element of the \nominus\ transversal of $\SSS_N$.
Next we develop 
useful facility with these objects,
and match to \S\ref{ss:main}.}

\beff{More composition notation}{{

\mdef \label{pa:modelss1}
  If $\lambda=\lambda_1 \lambda_2 \cdots \lambda_L$ is a composition then we write
  $\Gamma_\lambda$ for the set
  \[
\Gamma_\lambda = \times_{i=1}^{L} \Gamma_{\lambda_i}
\]
Example: $\Gamma_{22} = \{ (11,11), (11,2), (2,11), (2,2) \}$.
We may write $\mu \models \lambda$ for $\mu \in \Gamma_{\lambda}$. 
\\
If $\lambda$ is a partition
we write $\Gamma'_\lambda$ for the subset of $\Gamma_\lambda$ 
of multicompositions $\mu$
such that if there are parts of $\lambda$ of tied order then
their compositions appear in $\mu$ in (weak) lex order.
\\
Example $\Gamma'_{22} = \{  (11,11), (11,2), (2,2) \}$.
We may write $\mu \modelss \lambda$ for $\mu \in \Gamma'_{\lambda}$. 

}}

\medskip


\mdef \label{de:Gamma2}
We can regard
a composition
$\lambda\in\Gamma_N$ as a set --- of ordered parts ---
thus $\Gamma_5\ni\lambda = 1+3+1 = \{ (1,1),(3,2),(1,3)\}$, say.
Then $\PP_2(\lambda)$ is the set of two-part partitions of the set of parts.
Thus  
$
\Gamma^2_N = \bigcup_{\lambda\in\Gamma_M} \PP_2(\lambda) .
$
\\
We can think of an element of $\Gamma^2_N$ either just as the 2-part partition,
$s$ say,
or as the pair $(\lambda,s)$, where $\lambda$ is the underlying composition
(which can anyway be recovered from $s$).


\mdef \label{de:laml1}
The following is simply a diagram-free recasting of the $\jjj$-function machinery
from \ref{de:j}.
\\
Let $\lambda=\lambda_1 \lambda_2 \cdots \lambda_L \in \Gamma_N$
be a non-empty composition.
We define
$$
\gaml(\lambda) =\lambda_1 \lambda_2 \cdots (\lambda_L+1)
\hspace{1cm} \mbox{ and } \hspace{1cm}
 \gamr(\lambda) =\lambda_1 \lambda_2 \cdots \lambda_L 1
 .$$
Observe that both
$\gaml$ and $\gamr$
define injective functions from $\Gamma_N$ to  $\Gamma_{N+1}$.
\\
Similarly for
$(\lambda,s) \in\Gamma^2_N$ define $\gaml(\lambda,s)$ as above
with $(\lambda_L+1)$ in the same part of the 2-part partition as $\lambda_L$
was; and define $\gamr_0(\lambda,s)$ such that the added 1 is not in the same
part of $s$ as $\lambda_1$, and
$\gamr_1(\lambda,s)$ such that the added 1 is  in the same
part of $s$ as $\lambda_1$.
Again observe that these functions are injective.


\mdef
The first part of the following is well known; the second part not so.
\\
{\bf{Proposition}}. \label{pr:2N3N}
(a) We have that $\gaml(\Gamma_N) \cap \gamr(\Gamma_N) = \emptyset$
and
$\gaml(\Gamma_N) \cup \gamr(\Gamma_N) = \Gamma_{N+1}$.
Thus in particular
$
| \Gamma_N | = 2^{N-1}
$.
\\
(b) 
We have that $\gaml(\Gamma^2_N) , \gamr_0(\Gamma^2_N), \gamr_1(\Gamma^2_N)$
are disjoint,
and
$\gaml(\Gamma^2_N) \cup \gamr_0(\Gamma^2_N)\cup \gamr_1(\Gamma^2_N)
= \Gamma^2_{N+1}$.
Thus in particular
\[
| \Gamma^2_N | = 3^{N-1}
\]
\\
{\em Proof.} (a) Disjointness is clear. For completeness note that
every composition either ends in 1 or does not; and that in the former case
it lies in $\gamr(\Gamma_N)$ (indeed $\gamr$ has an inverse on the subset
ending in 1),
and in the latter in  $\gaml(\Gamma_N)$.
\\
(b) The argument for (a) extends straightforwardly.
The {\em three} subsets of $\Gamma^2_{N+1}$ to consider here are
elements whose composition ends in 1 and this is not in the same part of $s$
as $\lambda_1$;
elements whose composition ends in 1 and this is  in the same part of $s$
as $\lambda_1$;
and elements whose composition does not end in 1.
\qed



\mdef
Given an integer partition, or indeed a composition,
$\lambda = (\lambda_1, \lambda_2, ..., \lambda_m )$ of $N$
we have a partition of
$\ul{N} = \{1,2,...,N \}$ into
`nations':
$$
\ppp(\lambda) = \{
\{ 1,2,...,\lambda_1 \},
\{ \lambda_1 +1, \lambda_1 +2, ..., \lambda_1 + \lambda_2 \},
\hspace{1.65in}$$ $$\hspace{.65in}
\{ \lambda_1  + \lambda_2+1, \lambda_1  + \lambda_2+2, ..., \lambda_1  + \lambda_2+ \lambda_3 \},
...,
\{ 
\left(\sum_{i=1}^{m-1} \lambda_i \right)+1
, ..., N \}\}
$$

\ignore{{
\mdef
Given a set partition $p$ of a set $S$ of order $N$,  we
write $|p|$ for
the integer partition
of $N$ given by the orders of the parts of $p$
\ppm{$\;$--- a mild abuse of notation}.
\ppm{[-does this actually work?! ambiguous? needed!?]}
%
(NB  
If $S$ is ordered, then we may also extract
a {\em composition} $||p||$.) 
}}

\ignore{{

\mdef \label{de:Gammall2}
Define $\Gammall = \sqcup_{M \in \N} \Gamma_M$, the set of all compositions.
Define  $\Gammall^2 = \sqcup_{M \in \N} \Gamma^2_M$ similarly.

For $\chi\in \Gammall^2$ then $|\chi| = M$ such that $\chi\in\Gamma^2_M$.
And for $\chi\in hom(\Gammall^2,\N)$ then
$|\chi| = \sum_{\mu\in\Gammall^2} \chi(\mu) |\mu | $.
Similarly we can regard $\chi\in hom(\Gammall^2,\N)$ as a collection of
(augmented) compositions...
$||\chi||\vdash N$ is the integer partition corresponding to this collection
...
\\
Define \label{de:UUUN}
$$
\UUU_N = \{ \chi \in hom(\Gammall^2, \N) \; | \; |\chi|=N \}
.$$

\ppm{[Now state and prove the Theorem, as per ericpaulx.]}

}}



\module{braid-transversal}


\subsection{Symmetry orbits and transversals of solutions}
\label{ss:transversal} $\;$
%

\mdef
Our action of $\Sym_N$
can be used to 
take any solution $\FF$ to a representative with, say,
largest nations first; and then at fixed nation size to lex order on
compositions, so that all nations with composition $\lambda$ come before
all with composition $\lambda' > \lambda$.
\ignore{{
\sout{
However $P \in \ZZ_2$ reverses all compositions
(up to the $\Sym_n$ action) so takes $1+2$ to $2+1$
(but fixes $(1+2,2+1)$, although not necessarily the $\pai,\pfii$
part of a representation).
}}}


\mdef
\ignore{{
\sout{
For $M \in \N$ let  $\Gamma_M$ denote the set  
compositions of $M$ --- ordered sequences of natural numbers summing to $M$.
(The set $\Gamma_M$
has size $2^M$ and 
can be totally ordered for example lexicographically:
$1+2\; <\; 2+1$ etc.)
We write  $\Gamma_M^2$ for the set of all partitions into two parts of
the ordered summands of a composition, over all compositions in $\Gamma_M$.
This is  a set of size $3^{M-1}$ that is straightforward to enumerate explicitly.}
By  
{Lemma~\ref{lem:no-}} a partial classification of the $\Sigma_M$-orbits of varieties
of single-nation representations on $M$ say is given by $\Gamma_M$.
}}
{\bf{Proposition}}. 
A  
classification
of the $\Sigma_M$-orbits of 
varieties
of
single-nation representations on $M$ say
is 
given by
$\Gamma_M^2$
--- for $(\gamma,s) \in \Gamma_M^2$ the vertices in the two parts of $s$ are assigned
the $\alpha$ and $\beta$ parameter respectively
(e.g.
$\Gamma_2^2 =\{ (2,\{\{1,2\}\}),\; (1+1,\{\{1\},\{2\}\}),\; (1+1,\{\{1,2\}\}) \}$
corresponds to the set of edge labels $\{ 0, \ai, \fii \}$).
\\
{\em Proof}.
  This follows from {Lemma~\ref{lem:no-}}
  (or from our composition diagram construction and $\Rec$).
\qed


\mdef
A classification of orbits of varieties of representations on $N$ is then by
multiplicities $f \in \hom( \sqcup_{M \in \N} \Gamma_M^2 , \N_0 )$
such that the total degree is $N$,
as in the Theorem.
{
{ }  Combinatorics for this follows immediately from
our Proposition~\ref{pr:2N3N}(b) and 
the classical
construction of multiset sequences
(see for example \cite{SloanePlouffe95,Flajolet}
for the Euler transform).
}

\ignore{{
\ppm{[we have better way to do this now:]} \sout{
For this reason it is convenient to organise the nation content first
by nation sizes --- thus as an integer partition of $N$ ---
and then to give the multiplicities of nation types of each size $M$
--- compositions of $M$.}
}}

\ignore{{



\medskip


\mdef \label{pa:transvers1}
Observe that
$\Sigma_N$ acts on $\SSS_N$ via its action on $\ul{N}$, 
{and that this is
effectively
the same as the action induced via $\Rec$ and
Lemma~\ref{lem:chacon1}(II).}
For a transversal of $\Sigma_N$-orbits we proceed as follows.

First observe that in the $\Sigma_N$ orbit of any $(p,q,\rho,s)$ there are elements
in which the vertices in the largest nation
(or possibly one of the joint largest)
are relabelled to be consecutive from 1;
followed by the next largest,
and so on.
\ppm{Equivalently, in enumerating a transversal, the nations may be
considered to be unordered. Thus for such an enumeration it is sufficient
to enumerate a transversal of one-nation elements and then apply the multiset
combinatoric as in (\ref{de:multiset}). [-to be added from ericpaulx]}
Observe then that the orbits of nation partitions
\ppm{[-say this better?: 
that in each orbit
there is a representative
where $p=p_{\Lambda}(\lambda)$ for some integer partition $\lambda$
... hmmm, better just aim here for $\UUU_N$...]}
may effectively
be transversed by integer partitions.
{
If the nation sizes are not unique then we may order further so that the
tied nations are ordered according to their county content and ordering,
as we describe next.

Secondly,
in the $\Sigma_N$ orbit there is,
by Lemma~\ref{lem:no-},
an element in which the vertices
within a county are relabelled to be consecutive.
Furthermore they may be ordered so that the nominal order of these counties
agrees with $\rho$.
Note that configurations of this form are characterised by compositions.
Finally, in case of ties in nation size we may reorder nations by lex
order of compositions. The $s$ partition must simply be consistent
with the prescribed order.}

\mdef \label{de:TTTN}   
We have  
\[
\TTT_N \cong \{ (\lambda,\mu,s) \; | \;
              \lambda \vdash N, \mu \modelss \lambda,  s \in \PP_2 (\ppp(\mu))  \}
\]
Here  
$\mu \modelss \lambda$ denotes that $\mu$ consists of a composition of each part
of $\lambda$, with the compositions of order-tied parts written in lex order.
\\
For example $(1+1+1,1+2,1+2,2+1,1) \modelss (3,3,3,3,1) \vdash 13$.
}}
\ignore{{
For one particular choice of $s$ for this $(\lambda,\mu)$ pair we have the following abbreviated notation, where each connected component represents one nation, and the vertices at the same height within a nation correspond to a county:
\beq \label{eq:eg13}
\includegraphics[width=5.8\ccm]{tex-eric/diagram132.pdf}
\eq
We have suppressed the labels $1,\ldots, 13$ as the vertices are taken to be linearly ordered from top to bottom and left to right, and similarly suppressed the arrows as they all point downwards.

\mdef
The 
matrix $\FF(\sigma)$
corresponding to the example in (\ref{eq:eg13}) above
consists of $N=13$  scalar $1\times1$ matrices
$a_1,\ldots,a_{13}$ given by the vertex labels and
${{N}\choose{2}} = 78$ $2\times 2$ matrices $A(i,j)$ for $1\leq i< j\leq 13$,
that can be read off the diagram.
For example there are $66$ matrices of the form $\left(\begin{array}{cc} 0& \mu_{ij}\\\mu_{ij} &0\end{array}\right)$ for $1\leq i< j\leq 5$ corresponding to the $/$ edges between the vertices in distinct nations.  There are $3$ matrices $A(5,6),A(8,9)$ and $A(10,11)$ of the form $\alpha_i I$ for $i=2,3,4$ corresponding to the counties with population at least 2.
The matrices $A(1,2)=A(1,3)=A(2,3)=\left(\begin{array}{cc} \alpha_1+\beta_1& \beta_1\\-\alpha_1&0\end{array}\right)$ come from the two present and one implied edge of the first nation.  The remaining $2\times 2$ matrices corresponding to edges between vertices within a county are also of the $+$ form
$\left(\begin{array}{cc} \alpha_i+\beta_i& \beta_i\\-\alpha_i&0\end{array}\right)$ for $i=2,3,4$.

\newcommand{\Chin}{\chi} 

\sout{
\mdef
The recipe for constructing an explicit element of $\SSS_N$ from an element of the
formal transversal $\TTT_N$ (and hence a solution, via the above recipe)
is as follows.
\\
From $\lambda$ we construct the set partition $\ppp(\lambda)$.
From $\mu$ we construct the refinement $\ppp(\mu)$.
The order on the counties is simply the natural one,
which yields
$\rho_{\lambda\mu}\in Perm_{\ppp(\mu)}(\ppp(\lambda))$.
From $s$ we immediately have our element of $\PP_2 (\ppp(\mu))$.
Thus the map $\Chin: \TTT_N \hookrightarrow \SSS_N$ is
}
\[
\Chin : (\lambda, \mu, s) \mapsto (\ppp(\lambda),\ppp(\mu),\rho_{\lambda\mu},s)
\]


}}




\section{Useful examples and schematics} \label{ss:figs}

\module{notation-partition1}


\mdef
Each of our
2-coloured  composition
tableau  $\lambda\in \SSS_N$  
yields a decoration of the vertices and edges of $K_N$,
and hence a 
variety of functors via $\Rec$.  
To see the consistency of the construction, as in \ref{le:chip},
we may firstly visualise the clustering of vertices into counties and nations
as in Fig.\ref{fig:puppyears}.
\begin{figure}
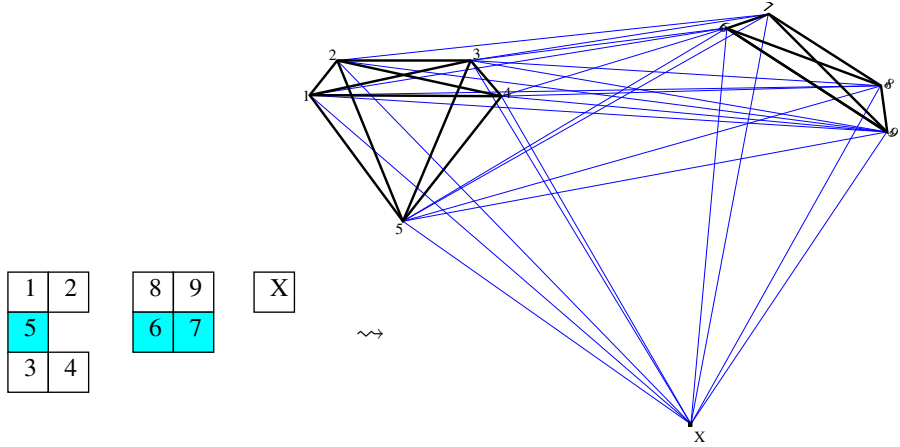

\[
\raisebox{.251in}{
\includegraphics[width=3.84cm]{xfig/tree012ex10bbb.eps}
}
\hspace{.2in}
\raisebox{.53in}{
$\leadsto$}
\hspace{-.431in}
\includegraphics[width=7.84cm]{xfig/yk212a22x5.eps}
\]
\caption{Geometric complete graph from 
a multiset of 2-coloured compositions. \label{fig:puppyears}}
\end{figure}
Then all `long' edges between nations $i,j$ get a $\mu_{ij}$ decoration;
and
`short' edges a 0
(parameter inherited from vertices).
(For remaining intranation edges see e.g. Fig.\ref{fig:hexa321a}.) 



\newcommand{\PPu}{\PP^u}
\newcommand{\PPm}{\PP^{-}}

\ignore{{
\mdef
Note that there  are natural bijections between the following sets:
\\
$\{ (p,q) \; | \; p,q \in \PP(S), p<q \}$,  $\;\;$
$\{ (q,p) \; | \; p,q \in \PP(S), p<q \}$,  $\;\;$
\hspace{.01in}
$ \{ (q,p') \; | \; q \in \PP(S), \; p' \in \PP(q) \}$
$\;$
and
$$
\PPu(S) =
 \bigcup_{q \in \PP(S)} \PP(q) = 
 \{ p' \; |  \; p' \in \PP(q), \; \mbox{for some } q \in \PP(S) \}
$$
 \\
 The map from first to second is simply to flip the pair.
 The map from second to third uses that $p<q$ has elements that are unions of
 parts from $q$, so we obtain a $p'$ by partitioning $q$ according to this data.
 Finally
 given $p' \in \PP(q)$ we can recover $q$ by taking the union over parts,
 so the $q$ is redundant.
 That is we can map $(q,p') \mapsto p'$.
 
 Note that $\PP$ has a kind of inverse $\PPm$ that takes an element of
any $\PP(S)$ and forms the union, thus restoring $S$
(caveat: this construction does not play well with the convention that
the result of a function applied to a subset of the domain is the set of images,
since here it may not be clear if an argument is a subset of the domain or an
element of it).
\\
Example:
$$
\PPu(\{1,2\})=
\PP(\{ \{1,2\}\}) \cup \PP(\{ \{1\},\{2\}\})
  = \{ \{\{\{1,2\}\}\}, \{\{ \{1\},\{2\}\}\}, \{ \{\{1\}\},\{\{2\}\}\} \}
$$
  
}}

\ignore{{
Sequences of form $(p,q,r,...)$ with $p<q<r<...$
(any length of sequence, $l$)
are sometimes called ramified partitions.
}}


\beff{Permutations}{{

\mdef
Here we consider the $(p,q,\rho,s)$ form of $\SSS_N$ as in \ref{de:typespace}.
  Given a  
  finite set $S$ then 
  $\Perm(S)$ 
  is the set of total orders of the  
  elements, i.e.
  $\Perm(S) = Hom^{iso}(S,\ul{|S|})$,
  the set of isomorphisms
  (although we express elements simply by giving the order).
\\
If $p<q$ as in \ref{pa:child1st} {\em et seq}
then $\Perm_p(q)$
denotes the set of total orderings of  the elements of $q$ in each $p_i$
(counties in each nation)
without imposing an overall total order:
$$
\Perm_p(q) = 
            \prod_{p_i \in p}  \Perm(q|p_i)
$$ 
Example: \hspace{.1in}
$
\Perm_{ \{ \{1\},\{2,3\}\}}( \{ \{1\},\{2\},\{3\}\}) =
\{  \{ (\{1\}),(\{2\},\{3\}) \} , \{ (\{1\}),(\{3\},\{2\}) \} \}
$. 
\\
Finally,
$\PP_2(q/p) = \prod_{p_i \in p} \PP_2(q|p_i) \;$
--- the set of partitionings of counties of each nation into at most two parts.
Example:
$ \PP_2(\{\{1\},\{2\},\{3\},\{4\}\}/\{\{1,2\},\{3,4\}\})
  = \PP_2(\{\{1\},\{2\}\}) \prod \PP_2(\{\{3\},\{4\}\})$.  
}}



\begin{figure}[h]
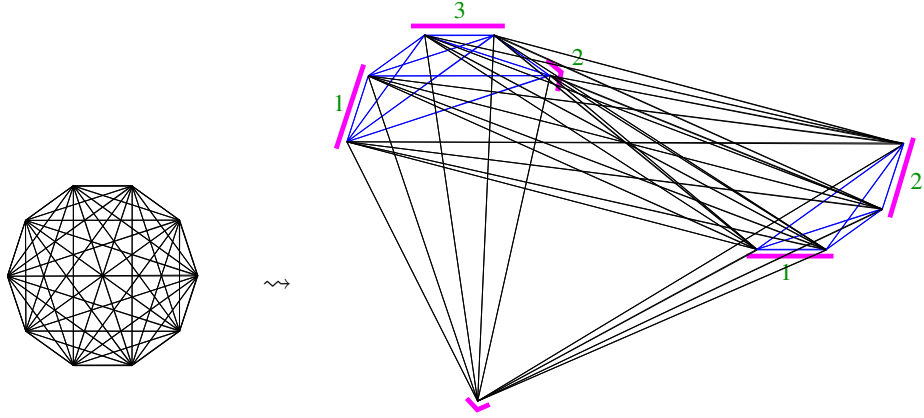

    \[
    \includegraphics[width=2.52\ccm]{xfig/reg10gon7.eps}
\hspace{.71\ccm}   \raisebox{1cm}{ $\leadsto$}
     \hspace{.41\ccm}
     \raisebox{-1\ccm}{   \includegraphics[width=7.9\ccm]{xfig/SN3cccc.eps} }
\]
\caption{Intermediate representation schematic with $N=10$. \label{fig:10gonx}}
\end{figure}


\ignore{{
\subsection{Constructions for solutions} $\;$ 

\sout{
Beside the geometrical 
encoding
we have a combinatorial encoding
of $(p,q,\rho,s)$
in terms of collections of Young tableaux
--- one tableau for each nation; one row for each county,
rows stacked  
according to $\rho$,
with each nation coloured so that counties  
that are in the same part in $s$ have the same colour.
The {\em same} element (but now  
complete with $s$ data)
becomes:}

}}


\beff{Classification: first schematic}{{

\mdef \label{de:SSSN}
For $N \in \N$, with the above notations, we have:
\[
\SSS_N =
\{ (p,q,\rho, s) \; | \; p<q \in \PP(N);\; \rho\in \Perm_p(q); \; s \in \PP_2(q/p)   \}
\]

One further schematic
(neglecting only $s$)
for elements
of $\SSS_N$ 
is 
as in Fig.\ref{fig:10gonx}.
%
%
In the schematic the well-separation of vertices again indicates $p$, thus
$\;$ 
$p= \{ \{\{ 1,2,3,4,5\},\{6,7,8,9\}, \{\barroman{X}\}\}$.
The coloured bars indicate $q$, thus
$\;$ 
$q=\{\{ 1,2 \},\{3,4\},\{5\}, \{6,7\}, \{8,9\}, \{\barroman{X}\}\}$.
And from the numbering
$\;$ 
$\Perm_p(q/p) = \{ (\{ 1,2 \},\{5\},\{3,4\}),( \{8,9\},\{6,7\}),( \{\barroman{X}\})
\}$.
Up to $s$ data we thus have
the same $\lambda\in\SSS_{10}$ as in Fig.\ref{fig:puppyears},

}}

\ignore{{

\newcommand{\ovb}{\overbrace}

\mdef
Example. For $N=2$ we have $\;$ 
 $
\SSS_2 = \{\;
\includegraphics[width=1cm]{xfig/tree012ex2aab.eps}\;,\;\;
\includegraphics[width=.51cm]{xfig/tree012ex2ab1b.eps}\;,\;\;
\includegraphics[width=.51cm]{xfig/tree012ex2ba1b.eps}\;,\;\;
\includegraphics[width=.51cm]{xfig/tree012ex2abb.eps}\;,\;\;
\includegraphics[width=.51cm]{xfig/tree012ex2bab.eps}\;,\;\;
\includegraphics[width=1.3cm]{xfig/tree012ex2acb.eps}\;
\}.
\ignore{{
\\ { } \hspace{.18in} =
\{ \;
(\; \ovb{\{\{1,2\}\}}^p, \; \ovb{\{\{ 1,2 \}\}}^q, \;\;\; \ovb{(\{1,2\})}^\rho ,
 \;\;\; \ovb{\{\{1,2\}\}}^s ),
$\\ $ { } \hspace{1.3\ccm}
(\; \{\{1,2\}\}, \{\{ 1\},\{2 \}\}, (\{1\},\{2\}), \{\{1,2\}\}), \;
$ \\ $ { } \hspace{1.3\ccm}
(\; \{\{1,2\}\}, \{\{ 1\},\{2 \}\}, (\{2\},\{1\}), \{\{1,2\}\}), \;
$ \\ $ { } \hspace{1.3\ccm}
(\; \{\{1,2\}\}, \{\{ 1\},\{2 \}\}, (\{1\},\{2\}), \{\{1\},\{2\}\}), \;
$ \\ $ { } \hspace{1.3\ccm}
(\; \{\{1,2\}\}, \{\{ 1\},\{2 \}\}, (\{2\},\{1\}), \{\{1\},\{2\}\}),
$ \\ $ { } \hspace{1.3\ccm}
(\{\{1\},\{2\}\}, \{\{ 1\},\{2 \}\}, (\{1\},\{2\}), \{\{1\},\{2\}\})
\}
}}
$
\\
To enumerate
we first vary over nation partitions
$\PP(2) = \{ \{ \{1,2\}\}, \{\{1\},\{2\}\} \}$. 
For the single-nation $p =  \{ \{1,2\}\} $
cases we  
have two possible refinements into counties,
again given by $\PP(2)$.
Here each county partition has two total orders; and there are two possible $s$
partitions.
Finally we have
an element with 
two nations, each a singleton so the other components are forced.

}}


\ignore{{
\input tex/reps-preamble03  
\section*{Attempt to integrate Eric-standalone}
\input tex/eric-summary
\input tex/eric-standalone-p
\section{Computation and Analysis for $(1,1,1)$} \label{ss:caa}
\input tex/L3  
\input tex/mat9macro   
\section{Computation and analysis for cases of form $(0,1,1)$} \label{ss:caa2}

We claim that cases of form $(0,1,1)$ give a transversal of the set of
cases of the forms $(0,1,1)$, (1,0,1) and (1,1,0).
Within $(0,1,1)$
observe that 
we can  exhaust all possibilities by checking the following table of
$(0,x,y)$ possibilities:
\[
\begin{array}{|c|c|c|c|c|c|c|}
  \hline
  (0,x,y) & 1.0 & 1.1i & 1.1ii & 1.2i & 1.2ii \\ \hline
  1.0     &  5  &  -   &   -   &  -   &  -    \\ \hline
  1.1i    &  -  &  3   &   -   &  -   &  -    \\ \hline
  1.1ii   &  -  &  -   &   4   &  -   &  -    \\ \hline
  1.2i    &  -  &  -   &   -   &  1   &  -    \\ \hline
  1.2ii   &  -  &  -   &   -   &  -   &  2    \\  \hline
\end{array}
\]
To explain the entries in the table,
we now turn (until we can include the tex source here) to Eric's file
\\ {\sf{somesolutionsLequal3011caseAND111case.pdf}}.
\\ (In the file first $x,y$ of form 1.2* are tried.
Then 1.1* -- which we would expect to be related to 1.2* by symmetry,
and this is the case.
Then 1.1 with 1.2 (no solutions).
Then $x=y=1.0$.
Then 1.1 or 1.2 with 1.0 (no solutions).)
The numbers in the table refer to the order in which solutions are
given in that file. Thus for example the first solution in the file is
of type $(0,1.2i,1.2i)$.
We write ``$-$'' to report finding {\em no} solutions of this type.

Now we should say something intelligent about the solutions 1-5
above. I.e. to characterize them intrinsically and extrinsically...
\\
Case 1: The first case seems to work as a kind of `cancelling meld' of type 1.2i
(recall i here means 1 boson in the Rittenberg classification)
in the 13 and 23 subspaces...
\\
Is this new? It certainly looks interesting...
How does the cancel work?
\\
We have, rescaling $a_{11}$ to 1 and applying similarity
(cf. \eqref{eq:cbnxx}):
\[
\acdc{0}{-x}{1}{1+x}{0}{-x}{1}{1+x}{x}
\]
\ppm{Hi Eric. In order to hit a `transversal'
(of some classification -- I have in mind the one from ordinary rep
  theory broadened to `Hamiltonian rep theory'...)
  the idea here is to apply
  overall renormalisation (its a group with generators with
  homogeneous relations!) to bring the $a_{11}$ eigenvalue to 1; and
  then to similarity transform so that as many below-the-diagonal
  entries as poss are 1. Does this work?!...}
\\
Case 2: ...
\[
\acdc{0}{1-x}{1}{x}{0}{1-x}{1}{x}{1}
\]
\\
Cases 3,4 are similar...
\\ ...
\\
Case 5: has the eigenvalue of the trivialised 12 sector;
then the
$\pm$ eigenvalues in the 13 and 23 sectors are forced to be the same
between the two sectors;
and the pure-3-sector eigenvalue can be chosen independently.
\[
\acdc{0}{x}{1}{0}{0}{x}{1}{0}{y}
\]
We might characterise this as a `symmetric meld'.

(OLDER:
We need to drag this output into the doc here and then apply
as much human analysis as possible.
We asked in particular if there might be cases included that are
essentially collapsing to $\NL = 2$ cases.
The answer is no, not exactly, but sort of...)

...

\section{A main Theorem for the $\NL = 3$ case} \label{ss:mainL3}

Tyeing the above together neatly we have
...
A natural notion of equivalence ...
Comparison of the nominal classification with the notion of
equivalence ...

\section{So how to we nicely characterise a transversal?} \label{ss:transvers}

What does the question mean? For example, a transversal when $\NL=2$
is largely provided by the Rittenberg classes with
$2=2+0$ (types 1.Xii) and
$2=1+1$ (types 1.Xi).
... together with the appropriate version of the reps discussed in  KMRW  (type 1.0).
...

\section{General $\NL$} \label{ss:genN}

A big (but not impossible-looking) question is how to lift to general
$\NL$. We can certainly give all the constraints generalising
\eqref{}-\eqref{} already.
Specifically we have the orbits under the action given in
\S\ref{ss:sym} of the transversal given in \eqref{}-\eqref{}.
\\
And the organisation scheme lifts in the obvious way.

For example, for $\NL = 4$ we have pairings 12, 13, 14, 23, 24, 34,
so a 6-fold nominal classification. 

Let us start with representations of `type' (1.0,1.0,...).
Can we show that they always exist?
...

...ERIC SAYS: yes. ...
And indeed Eric will populate a verification of this assertion.
(In our discussion P argued that this case should be `given respect' --
being perhaps intrinsically easy, but extrinsically significant... 

\medskip

\input tex/maundringseric
}}

\bibliography{
local.bib,../new31.bib,HGT2.bib}{}
\bibliographystyle{alpha}
\end{document}